\documentclass[preprint,12pt]{elsarticle}

\usepackage{amssymb}
\usepackage{amsthm}
\usepackage{subfig}
\usepackage[monochrome]{color}
\usepackage{graphicx}
\usepackage{natbib}
\usepackage{amsmath}
 \usepackage{booktabs}

	\date{2020}
	
\begin{document}

	\begin{frontmatter}



	\title{\textcolor{magenta}{Extrapolated Shock Tracking: bridging shock-fitting and embedded boundary methods}} 

		
	\author[lab1]{Mirco Ciallella }
	\author[lab1]{Mario Ricchiuto }
	\author[lab2]{Renato Paciorri }
	\author[lab3]{Aldo Bonfiglioli}
	\address[lab1]{Team CARDAMOM, INRIA Bordeaux Sud-Ouest, 33405 Talence, France}
	\address[lab2]{Dip. di Ingegneria Meccanica e Aerospaziale, Universit\`a di Roma ``La Sapienza'', Via Eudossiana 18, 00184 Rome, Italy}
	\address[lab3]{Scuola di Ingeneria  - Universit\`a degli Studi della Basilicata,Viale dell'Ateneo Lucano 10, 85100 Potenza, Italy}
	\begin{abstract}
	We propose a novel approach to approximate numerically shock waves. The method combines  the unstructured shock-fitting approach developed in the last decade by some of the authors, with ideas coming from embedded boundary techniques. 
	The numerical method obtained allows avoiding the re-meshing phase required by the unstructured fitting method, \textcolor{magenta}{while guaranteeing accuracy properties very close to those of the fitting approach. 
	 This new method has many similarities with front tracking approaches, and paves the way to shock-tracking techniques }  truly independent on the data and mesh structure used by the flow solver. 
	 \textcolor{magenta}{The approach} is tested on several problems showing accuracy properties very close to those of more expensive fitting methods, with a considerable gain in flexibility and generality.
	\end{abstract}


	\end{frontmatter}

	
\section{Introduction}

The numerical techniques used to simulate flows with shock-waves  are 
essentially two: the widely used shock-capturing (SC) methods, and 
the less common shock-fitting (SF) methods.
The former relies on the proven mathematical legitimacy of weak solutions: 
all types of flows, including flows with shocks, can be computed by using 
the same discretization of the equations in divergence form. Nevertheless, 
the shocks always appear smeared in a region whose thickness is of two or three cells 
rather than actual discontinuities.
In addition to this, but perhaps more correctly because of this, since the states of the cells
inside this region are unphysical~\cite{zaide2011shock}, the shock-capturing 
methods suffer from some numerical problems concerning the stability, the accuracy and 
the quality of the solutions that sometimes give anomalous results.
A {\it catalogue of} these {\it failings}  was made by Quirk in the early 90s of the last 
century~\cite{Quirk1994555}.
Despite the  great efforts made by numerous  researchers in 
the last decades to develop shock-capturing methods, 
these numerical problems are not entirely solved and still plague the numerical solutions obtained by 
shock-capturing solvers. 

The shock-fitting technique for compressible flow computations 
has been  developed by Gino Moretti~\cite{Moretti1,Moretti2} in the 1960s.
It consists in explicitly identifying the shock as a line (surface in 3D) within the flow-field and computing its motion 
and upstream and downstream states according to the Rankine-Hugoniot equations. 
However, historically, the techniques developed by Moretti and his collaborators were designed 
for solvers based on structured grids  and this made their development very difficult and complex, 
especially when extended to flows with shock interactions~\cite{Moretti2}.

Two different shock-fitting methodologies blossomed between the 60s and 80s: the boundary shock-fitting 
and floating shock-fitting. In the former approach, the shock is made to coincide with one of the 
boundaries of the computational domain so that the treatment of the jump relations across the shock is 
confined to the boundary points. Even though this method greatly simplified the coding, the treatment
of shocks appearing within the computational domain and of shock interactions became a major challenge. 
The floating shock-fitting approach was developed to be capable of dealing with more complex flow 
configurations. In the floating version, discontinuities can freely move over a background structured mesh: 
a shock front is described by its intersections with the grid-lines, which give rise to \emph{x} and \emph{y}
shock points, meaning that they are allowed to move onto grid-lines. Even though  floating shock-fitting codes 
have been used with success in the past to compute steady and un-steady two- and three-dimensional flows
involving shock reflections and shock interactions~\cite{Nasuti1,Nasuti2,Zhong1}, they are very
complex to code and require extensive changes in the computational kernel of the gas-dynamic solver. 

In the 80s and early 90s the CFD community has shown increasing interest in unstructured 
meshes. This is mainly due to the features that characterize this kind of grids: the ability to easily mesh 
complex geometries and the possibility of locally adapting the mesh size to follow the flow features. 
The latter advantage makes them well-suited to simulate compressible flows with shock waves and contact discontinuities. 
Exploiting this flexibility, Paciorri and Bonfiglioli developed a new unstructured shock-fitting technique for unstructured
vertex-centered solvers, described in~\cite{Paciorri1}. This approach has alleviated many of the 
difficulties of the shock-fitting techniques in the structured-grid framework.  
In recent years, the unstructured shock-fitting technique was improved to deal with interactions among discontinuities in two-dimensional 
flows~\cite{Paciorri2}, three-dimensional flows \cite{Paciorri3} and un-steady compressible 
flows \cite{Paciorri4,Campoli2017} opening a new route in simulating flow-field with shock waves. In particular, not only 
shocks and contact discontinuities are fitted, but also the interaction points, for example the triple points arising in Mach reflections~\cite{Paciorri5}.
A limitation of this technique is that it heavily relies on the flexibility of triangular and tetrahedral grids to locally produce a fitted unstructured grid around the discontinuities.
This limits its application to unstructured vertex-centered codes.

Recently, the research  group headed by Prof.\ J. Liu proposed and developed a shock-fitting technique for unstructured 
cell-centered solvers~\cite{zou2017shock,Chang2019}. However, even this technique has an important  
limitation: it uses a deforming grid whose topology cannot be changed during the computation, 
unless an expensive re-meshing (and, consequently, interpolation of the solution) is carried out. This is an important
restriction, especially when shock-waves move throughout the flow-field  or whenever
new shocks appear during the computation.

Both these shock-fitting formulations currently available heavily rely on the data structure of the flow solver, and more particularly on the mesh. Indeed, in all these techniques, the jump conditions are attached to some  mesh entity (edge, face, or node). This very often makes the methods better suited for one or another family of flow solvers (node-centered, cell-centered, finite volume, finite element etc.), thus limiting its use. Another complication is that both these methods require the mesh to follow exactly the evolution of the shock wave, which puts additional requirements on the meshing/re-meshing techniques used.

In this work we aim at proposing a new approach, which is in some way more general and flexible. The initial idea  
 comes from the similarity between the constraints arising from shock-fitting, 
and those related to the construction of boundary-fitted grids for simulating flows around complex geometries. 
In this context, immersed and embedded boundary methods  have been developed since many years to allow a flexible management of complex geometries. The two approaches rely on a slightly different philosophy.

Immersed methods  are based on an extension of the flow equations outside the physical domain (typically within solid bodies). This extension is formulated using some smooth
approximation of the Dirac delta function to localize the
boundary, as well as to impose the boundary condition.
These methods are relatively old, and based on the original
ideas of Peskin \cite{IB0}. Finite element and unstructured mesh
extensions for elliptic PDEs as well as for incompressible, and compressible flows have been discussed in \cite{IB1,IB2,IB3}.

Embedded methods, on the other hand, solve the PDEs only in the physical domain, while
replacing the exact boundary with some more or less accurate approximation, combined
with some weak enforcement of the boundary condition. There is a certain number of techniques
to perform this task, which go from the combination of XFEM-type methods with
penalization or Nitsche's type approaches \cite{EM0},
to several types of cut finite element methods with improved stability \cite{EM1,EM2}, to
approximate domain methods such as the well known ghost-fluid method \cite{EM3,EM4},
and the more recent shifted boundary method (SBM) \cite{Scovazzi1,Scovazzi2}.

In this work we borrow ideas from approximate domain methods, and in particular from the SBM.
As in the latter, \textcolor{magenta}{we impose modified conditions on  surrogate
shock-manifolds,}
\textcolor{magenta}{acting  as boundaries} between the shock-upstream and shock-downstream regions.  
\textcolor{magenta}{These surrogate boundaries are composed of two sets of mesh faces enclosing the cavity of elements crossed by the shock.}
  \textcolor{magenta}{The values of  the flow variables imposed on these surrogate boundaries}
\textcolor{magenta}{ are extrapolated from the tracked shock front accounting for  the non-linear jump and wave propagation conditions, as done in the unstructured  shock-fitting  approach.  
As in the SBM, the extrapolation is based on a truncated Taylor series expansion from the surrogate boundaries to
 the front, allowing to preserve}
 the overall accuracy of the discretization.
 This paper, in particular, only deals with second-order piecewise linear approximations, but all the ideas
can be extended to higher order. Note however that differently from e.g.\ the extension of the SBM to hyperbolic problems \cite{Scovazzi3}, 
the approach proposed here requires the solution of \textcolor{magenta}{three  coupled  problems: the CFD  upstream of the shock, the CFD downstream of the shock,
the coupled algebraic system obtained from the Rankine-Hugoniot relations augmented with the characteristic information traveling toward the shock front}.
\textcolor{magenta}{As in  shock-fitting and  front tracking methods \cite{Glimm2016},     the shock front is explicitly discretized by an independent lower-dimensional mesh, and its position, as 
well as the position of the two    surrogate boundaries}, are themselves part of the computational result. These elements  make the present work  not only original w.r.t.\ previous shock-fitting methods, but also with respect to previous work in embedded methods and in particular the SBM approach. 
Indeed, the most recent work on the use of similar ideas, only considers interfaces independent on the solution, and linear elliptic partial differential equations~\cite{lishifted}. 
Moreover, the approach proposed in the reference is based on a single surrogate interface, while the approach proposed here uses a symmetric formulation with two surrogates. 
The resulting \textcolor{magenta}{method bears some similarities to front tracking approaches, and for this reason  is referred to as extrapolated Shock Tracking (eST) to differentiate it from previous unstructured shock-fitting methods in which the  faces of the shock mesh are part of the CFD meshes. This new method constitutes a}
 bridge between shock-fitting and embedded boundary methods.
It removes some of the constraints of the approach by Paciorri and Bonfiglioli, while keeping its flexibility. The method proposed is  actually  even more general as it constitutes a shock-fitting/tracking technique
virtually independent on the data structure of the underlying gas-dynamic solver. This paper focuses on the formulation
of the method in two-space dimensions, and on its validation on classical problems involving strong shocks, as well as on
problems with shock interactions, where a capability for   hybrid fitting-capturing computations is shown.

\section{Generalities}

We consider the numerical approximation of solutions of the  {\color{blue}steady} limit of the Euler equations  reading:
\begin{equation}\label{eq:euler0}
\partial_t \mathbf{U}+\nabla\cdot\mathbf{F} =0 \quad \text{in}\quad\Omega\subset\mathbb{R}^d
\end{equation}
with conserved variables and fluxes given by:
\begin{equation}\label{eq:euler0a}
\mathbf{U}=\left[
\begin{array}{c}
\rho \\ \rho \mathbf{u}  \\ \rho E
\end{array}
\right]\;,\;\;
\mathbf{F}= \left[
\begin{array}{c}
\rho \mathbf{u}\\ \rho \mathbf{u} \otimes \mathbf{u} + p \mathbb{I} \\ \rho H \mathbf{u}
\end{array}
\right] 
\end{equation}
having denoted by $\rho$ the mass density, by $\mathbf{u}$ the velocity, by $p$ the pressure, and with $E= e +  \mathbf{u}\cdot\mathbf{u}/2$ the specific total energy, $e$ being the specific internal energy. 
Finally, the total specific enthalpy is  $H=h+ \mathbf{u}\cdot\mathbf{u}/2$, with  $h=e+p/\rho$ the specific enthalpy. For simplicity in this paper we  work with the classical perfect gas 
equation of state: 
\begin{equation}\label{eq:EOS}
p=(\gamma-1)\rho e
\end{equation}
with $\gamma$ the constant (for a perfect gas) ratio of specific heats. However, 
note that the method discussed allows in principle to handle any other type of gas, {\color{green} see e.g.~\cite{pepe2015unstructured}}.

In all applications involving high-speed flows, solutions of~(\ref{eq:euler0})  are  only piecewise continuous.
In $d$ space dimensions, discontinuities are represented by $d-1$ manifolds governed by the  well known Rankine-Hugoniot jump conditions reading:
\begin{equation}\label{eq:euler1}
[\![\mathbf{F}\cdot\mathbf{n}]\!]=w[\![\mathbf{U}]\!]
\end{equation}
having denoted by $\mathbf{n}$ the local normal  vector to the shock, by $[\![\cdot]\!]$ the corresponding jump of a quantity across the discontinuity,
and with $w$ the normal component of the shock speed.\\
\textcolor{magenta}{As discussed in the introduction, the method proposed exploits ideas from two different approaches:
the unstructured shock fitting method \cite{Paciorri1} and subsequent works; the shifted boundary method by \cite{Scovazzi1} and subsequent works.
In the following sections we recall the main ingredients of these two techniques.}

{\color{blue}
\section{Unstructured shock-fitting algorithm}
We shall first briefly describe the unstructured shock-fitting technique developed 
by Paciorri and Bonfiglioli~\cite{Paciorri1,Paciorri2,Paciorri3}, }\textcolor{magenta}{in the following referred to with the acronym SF}.

\begin{figure}[!htb]%
\centering
\subfloat[Interaction of two shocks of the same family]{\includegraphics[draft=false,width=0.475\linewidth]{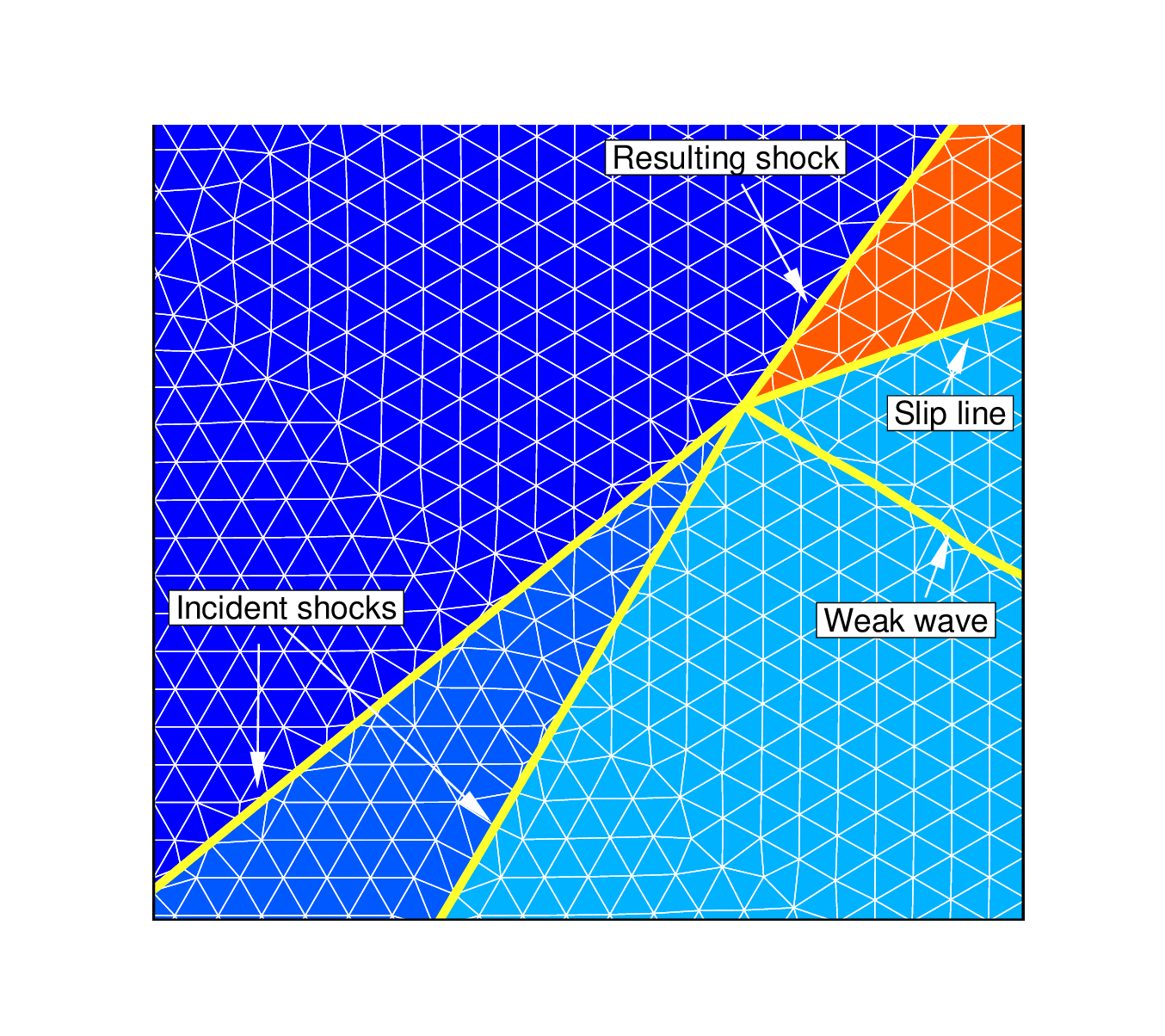}\label{fig:ss3-f5cc}}
\qquad
\subfloat[Pseudo-temporal evolution of the grid and the fitted discontinuities]{\includegraphics[draft=false,width=0.45\linewidth]{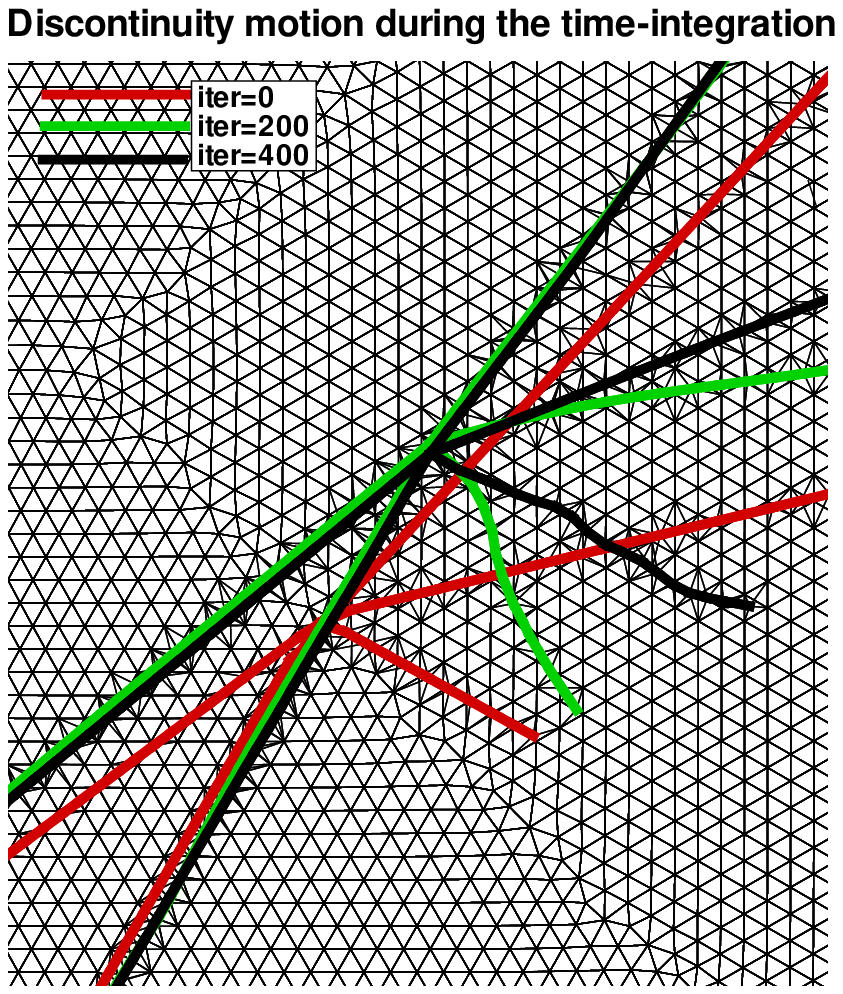}\label{fig:ss3-f3}} \\
\subfloat[Supersonic flow over a blunt-nosed body.]{\includegraphics[draft=false,width=0.475\linewidth]{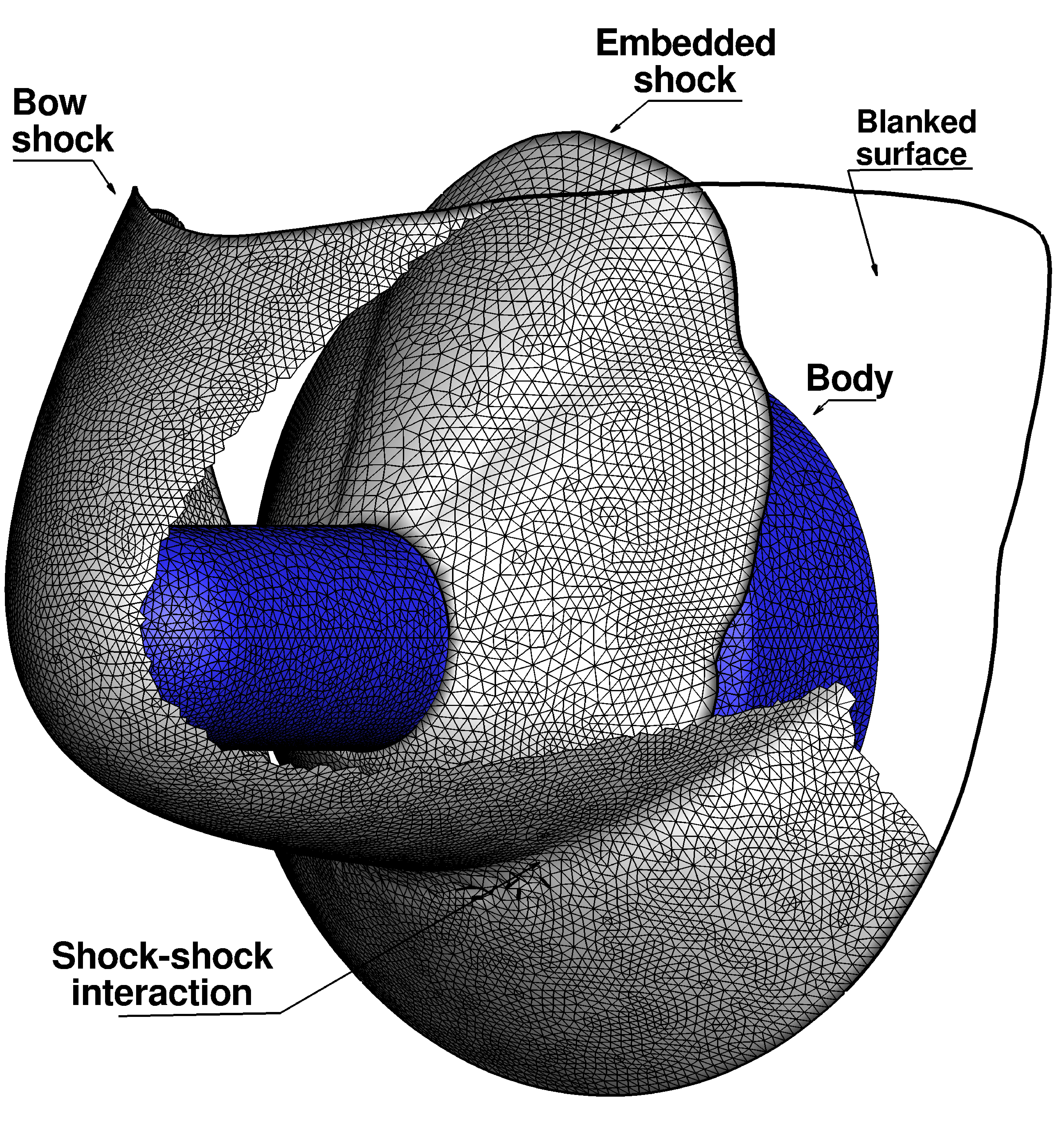}\label{fig:shk1_TC3_new}}
\caption{Examples of fitted discontinuities on unstructured meshes.}%
\label{fig:cont}%
\end{figure}

\textcolor{magenta}{In this approach} {\color{blue}the set of dependent variables is available within all grid-points of a tessellation 
(made of triangles in 2D and tetrahedra in 3D) that covers the entire
computational domain; this is what we call the {\em background} mesh.
In addition to the background mesh, the fitted discontinuities (either shocks or slip-streams) 
are discretised using a collection of points which are mutually joined to 
form a connected series of line segments, as shown in Fig.~\ref{fig:ss3-f5cc} for the 2D case, or a triangulated surface in 3D, 
as shown in Fig.~\ref{fig:shk1_TC3_new}. This is what we call the {\em shock} mesh.}
{\color{blue}For example, the thick solid (yellow) line in Fig.~\ref{fig:ss3-f5cc}
marks the various fitted discontinuities that arise due to the interaction between two shocks of the same family:
the two incident shocks, the resulting shock, a weak compression wave\footnote{could be an expansion wave instead,
depending on the upstream boundary conditions} and the slip-stream located between the former two. 
Figure~\ref{fig:shk1_TC3_new}, which refers to the three-dimensional, supersonic flow past a blunt-nosed object, 
shows the triangulated surfaces used to fit the bow shock and the imbedded shock that arises at the cylinder-flare junction.\\
Although it is not evident from Fig.~\ref{fig:cont}, each fitted discontinuity is a double-sided internal boundary of zero thickness. 
Being the width of the discontinuity negligible, its two sides
are discretised using the same polygonal curve or triangulated surface, so that
each pair of nodes that face each other on the two sides of the discontinuity share the same geometrical location,
but store different values of the dependent variables, one corresponding to the upstream state and
the other to the downstream one. 
Moreover, a velocity vector normal to the discontinuity is assigned to each pair of grid-points on the fitted discontinuity: it represents
the displacement velocity of the discontinuity.}
{\color{blue} The initial condition for a shock-fitting calculation is typically (see~\cite{Chang2019} for a different approach) 
supplied by running a shock-capturing calculation
on the background mesh; then, a feature extraction algorithm, such as the those described in~\cite{SFB:Andrea,Chang2019,PACIORRI2020109196}, 
is used to provide the initial (though approximate) location of the discontinuities.
Even when dealing with steady flows, the approach is inherently time-dependent, because
both the solution and the grid change with time, due to the displacement of the fitted discontinuities.
Whenever a steady solution exists, the shock speed asymptotically vanishes
and the tessellation of the flow domain does not any longer change.
This is illustrated in Fig.~\ref{fig:ss3-f3} which shows the pseudo-temporal evolution
of the various discontinuities involved in the shock-interaction of Fig.~\ref{fig:ss3-f5cc}.
Moreover, Fig.~\ref{fig:ss3-f3} reveals that
the spatial location of the fitted discontinuities is independent of the location of the grid-points that make up the background grid
and that local re-meshing only takes place in the immediate neighborhood of the moving discontinuities.}

\section{Shifted-boundary method}

\begin{figure}[!tb]%
\centering
\subfloat[The surrogate boundary $\tilde{\Gamma}$ and true boundary $\Gamma$]{\includegraphics[width=0.5\textwidth]{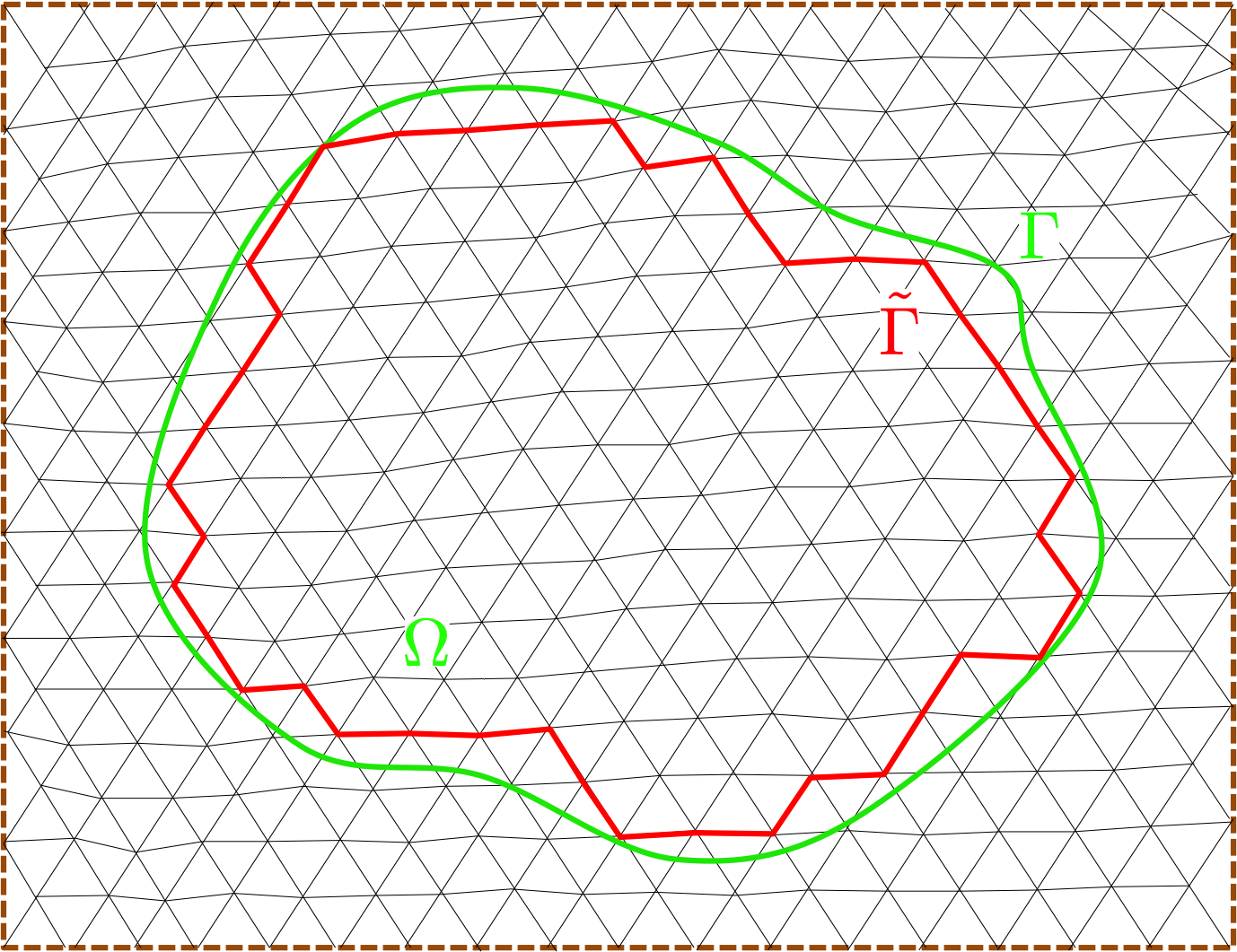}\label{fig:SBM1}}\quad
\subfloat[The distance vector $\vec{d}$ and unit normal and tangent vectors to the true boundary]{\includegraphics[width=0.45\textwidth]{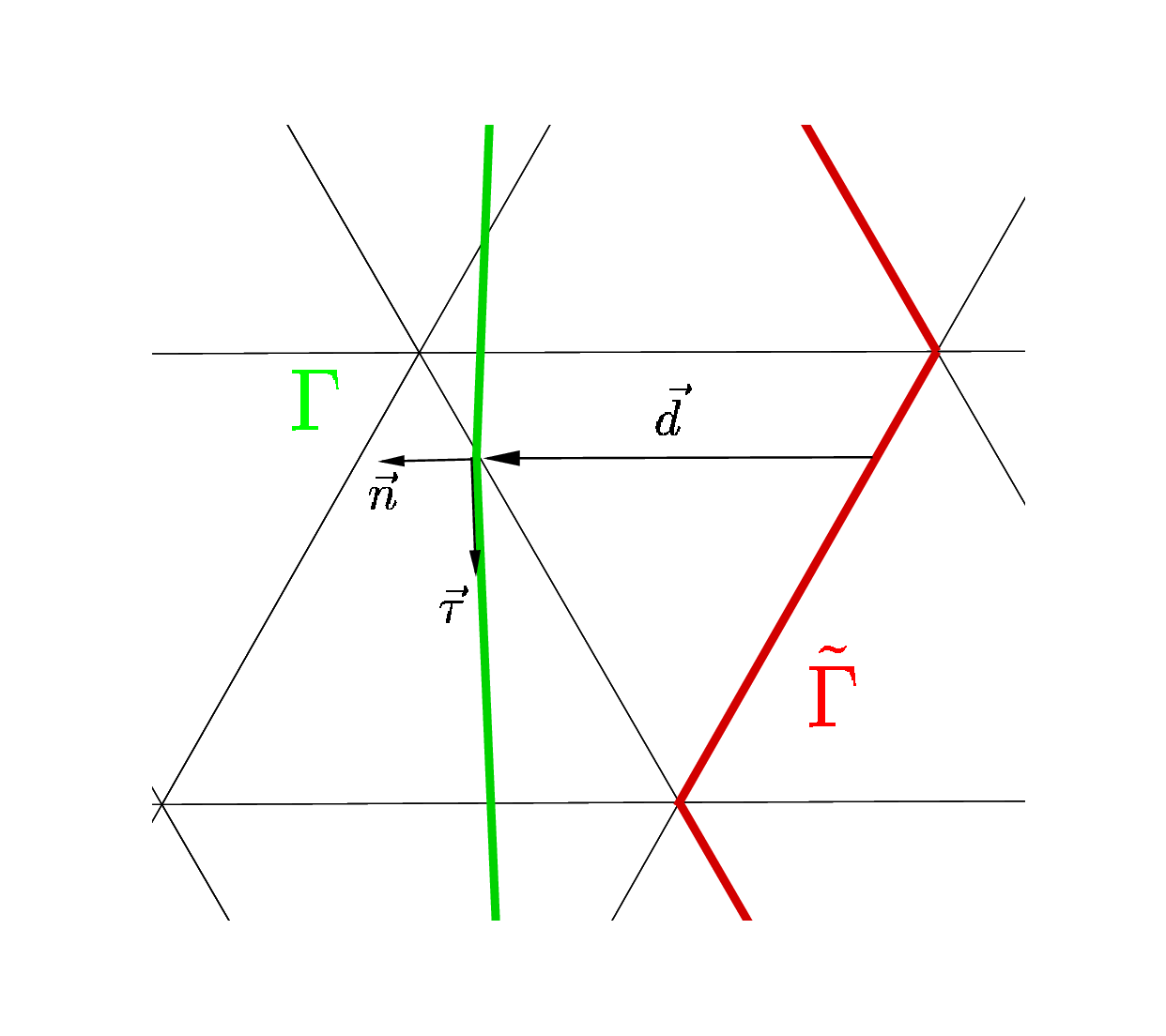}\label{fig:SBM2}} \\
\caption{The SBM: the surrogate and actual boundaries, and the distance vector $\vec{d}$.}%
\end{figure}

{\color{blue} The main advantage of embedded boundary methods, and among them the SBM \cite{Scovazzi1,Scovazzi2}, 
is the ease of mesh generation with respect to the classical body-fitted methods. It has been pointed out how trivial this task 
might be even when complicated geometries are taken into account. With all the benefits that characterize these approaches,
some shortcomings arose in the standpoint of the enforcement of boundary conditions. The originality of the SBM lies in the
idea of $shifting$ the location where the boundary conditions are applied.
In order to \textcolor{magenta}{guarantee consistency, and retain the mesh convergence rates of the original method,} 
the boundary conditions have to be modified.}\\
\textcolor{magenta}{The main steps of the method are the following. Given a mesh including the physical domain $\Omega$, not conformal w.r.t.\ the domain boundary $\Gamma$,
one must first define a surrogate boundary $\tilde\Gamma$. As shown in figure \ref{fig:SBM1}, $\tilde\Gamma$ is essentially built from the 
mesh faces and mesh points in $\Omega$ closest to the true boundary $\Gamma$. Next, for any point of the surrogate boundary $\tilde\Gamma$,
one needs to be able to define a map to a unique point of the true boundary $\Gamma$:}
{\color{blue}
\begin{equation}
\mathbf{M}\,:\,\tilde{\Gamma}\,\rightarrow\,\Gamma
\end{equation}
\begin{equation}
\mathbf{\tilde{x}}\,\rightarrow\,\mathbf{x}
\end{equation}
which maps $\mathbf{\tilde{x}}\in\tilde{\Gamma}$  on the surrogate boundary to 
 $\mathbf{x}\in\Gamma$  on the true boundary. The map \textbf{M} can be built} \textcolor{magenta}{in several ways, for example using a closest point projection, or  using level sets, or equivalently  using distances along directions normals to the true boundary $\Gamma$,} 
{\color{blue}as shown in Fig.~\ref{fig:SBM1}.
Since the gap between $\tilde{\Gamma}$ and $\Gamma$ is going to be of crucial importance, in terms of accuracy of the solution,
the map \textbf{M} will be characterized through a distance vector function:
\begin{equation}\label{dmap}
\mathbf{d_M(\tilde{x})}\,=\,\mathbf{x}\,-\,\mathbf{\tilde{x}}\,=\,[\mathbf{M\,-\,I}](\mathbf{\tilde{x}})
\end{equation}}
\textcolor{magenta}{If $\mathbf{M}$ is built using distances along  normals to   $\Gamma$, the vector  $\mathbf{d_M(\tilde{x})}$ is parallel to the normal to $\Gamma$ in $\mathbf{\tilde{x}}$.
} 
{\color{blue}Finally, the boundary conditions have to be modified to provide high-order (at least second-order) convergence rate of the solution.
This can be accomplished by writing a Taylor expansion formula centered at $\mathbf{\tilde{x}}\in\tilde{\Gamma}$, recalling Eq.~\eqref{dmap}:}
\begin{align}
u(\mathbf{x}) & = u(\mathbf{\tilde{x}})+\nabla u(\mathbf{\tilde{x}})\cdot(\mathbf{x}-\mathbf{\tilde{x}})+O(\|\mathbf{x}-\mathbf{\tilde{x}}\|^2) \nonumber\\
              & = u(\mathbf{\tilde{x}})+\nabla u(\mathbf{\tilde{x}})\cdot(\mathbf{M(\tilde{x})}-\mathbf{\tilde{x}})+O(\|\mathbf{M(\tilde{x})}-\mathbf{\tilde{x}}\|^2)\nonumber\\
              & = u(\mathbf{\tilde{x}})+\nabla u(\mathbf{\tilde{x}})\cdot\mathbf{d_M(\tilde{x})}+O(\|\mathbf{d_M(\tilde{x})}\|^2) \label{eq:sbm1}
\end{align}
{\color{blue}Equation~\eqref{eq:sbm1} is at most second-order accurate, unless additional terms in the Taylor expansion are included, as explained in \cite{Scovazzi4}}.
\textcolor{magenta}{Now, if on $\Gamma$ the prescribed boundary condition is $u(\mathbf{x})=g(\mathbf{x})$, the main idea of the SBM is to
deduce from Eq.~\eqref{eq:sbm1} that the boundary condition to be imposed on $\tilde\Gamma$ to allow for
 second order of accuracy w.r.t.\ $\|\mathbf{d_M(\tilde{x})}\|$ is  
 $$u(\mathbf{\tilde x})=g(\mathbf{M(\tilde x)})-\nabla u(\mathbf{\tilde{x}})\cdot\mathbf{d_M(\tilde{x})}.$$
 This extrapolation constitutes the main idea exploited in the following.
 }

\textcolor{magenta}{\section{ Extrapolated Shock-Tracking}
%
We discuss here the    extrapolated Shock-Tracking (eST) method we propose.
{\color{red}We focus on {\color{blue}steady state flows} in at most two space dimensions}, however most of the ideas 
discussed can be generalized to three space dimensions.
The  eST algorithm can be summarized in three main steps, allowing to
update the computational domains and solution values, 
 that lead from the available mesh and solution at pseudo-time $t$ 
to an updated mesh and solution at pseudo-time $t+\Delta t$:} 
\begin{description}
\item[1. (Shock/background-mesh coupling)] Geometrical coupling of the shock-mesh with the background-mesh, and definition of separate shock-upstream and shock-downstream computational domains;

\vspace{0.2cm}

\item[2. (Computational domain update)]  Iteration evolving in (pseudo-)time the flow variables in each computational domain independently; 

\vspace{0.2cm}

\item[3. (Shock update)] Evolution of the position of the shock and of the flow variables values at the shock, {\color{green} using
	the jump relations~(\ref{eq:euler1})}.
\end{description}
 
These three steps, 
are applied  iteratively until a steady state is obtained, {\color{blue}or, when dealing with un-steady flows,
in a time-accurate manner~\cite{Paciorri4,Campoli2017}}.
The most specific ingredients  of the method are those of steps 1 and 3.
Indeed, step 2 essentially relies on the use of an accurate multidimensional upwind unstructured grid solver to compute the smooth flows
upstream and downstream of the shock, and additional discontinuities not  being fitted by the above method. 
We will briefly recall in Sect.~\ref{sect:step of eulf} the cell-vertex solution method used here. 

The main difference between the technique described here and the one proposed by Paciorri and Bonfiglioli~\cite{Paciorri1,Paciorri2,Paciorri3}  
is in step 1. 
Indeed, the present technique removes the need to insert the shock-mesh in the background mesh, which can be a critical aspect, 
especially when different shock-surfaces mutually interact in the three dimensional space~\cite{Paciorri3}. 
To this end, we exploit ideas coming from embedded boundary methods. In particular, we propose to use an extrapolation from the background mesh 
to the shock mesh in the spirit of the {\color{green}SBM} initially proposed in~\cite{Scovazzi1}  for elliptic problems 
and extended to hyperbolic problems in~\cite{Scovazzi3}.
In other words, the method proposed consists in replacing re-meshing with the definition of  sufficiently accurate extrapolation functions, which allow the transfer of information between the background and shock meshes. This allows to completely remove the need of re-meshing.

\textcolor{magenta}{As already mentioned in the introduction, we refer to this new method
as to extrapolated Shock-Tracking (eST) to differentiate it from unstructured shock fitting, in which a conformal mesh 
fitting the shock front is generated,
and to differentiate it from the SBM, in which  the  true boundary is replaced by
a unique surrogate  with extrapolated boundary values. The eST method actually has similarities with high-order front tracking
approaches~\cite{Glimm2016}, and also  for this reason we prefer referring to it as shock-tracking}.
{\color{green}It may also be viewed as
some sort of elaborate
solution optimization procedure in which, starting from a captured result, one iteratively places the shock front and modifies the 
flow solution by
solving the nonlinear jump conditions.  In this respect there are  similarities with  approaches 
optimization based coupled with mesh adaptation as those  recently proposed e.g.\ by~\cite{persson,corrigan} in the 
framework of Discontinuous Galerkin methods.
}
{\color{magenta}
 The details of the method are discussed in the next sections, highlighting}
{\color{blue}
the major changes and differences w.r.t.\ the {\color{green}SF} approach of~\cite{Paciorri1,Paciorri2,Paciorri3}. }

	\subsection{Geometrical setting}

To illustrate the algorithmic features of the eST method, let us consider a two-dimensional domain and a shock front crossing the domain at a given time $t$ (see Fig.~\ref{pic1}). 
The shock front is described by a {\color{green}collection of shock-edges} whose endpoints are the shock-points, marked by squares
{\color{green}in Fig.~\ref{pic1}.
Shock-edges and shock-points make up the shock-mesh}. 
A background triangular mesh, whose {\color{green}grid-points} are denoted by circles
{\color{green}in Fig.~\ref{pic1}},
covers the entire computational domain. It is noted that the position of the shock-{\color{green}points} is completely 
independent of the location of the grid-points of the background mesh. 
While each grid-point of the background mesh is characterized by a single set of dependent variables, two sets of values, 
corresponding to the upstream and downstream states, are assigned to each shock-point. 
We assume that at time \emph{t} the solution is known at all grid- and shock-points. 
The computation of the subsequent time level \emph{t+$\Delta$t} can be split into several steps that will be described in detail in the following sub-sections.

\subsection{Cell removal around the shock front}\label{sect:cellremoval}

\begin{figure}
\begin{center}
	\subfloat[{\color{blue} Shock-mesh laid on top of the background-mesh}]{%
		\includegraphics[height=0.36\textwidth]{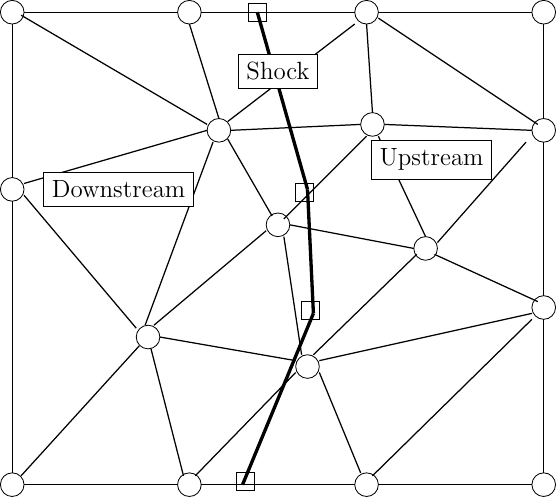}\label{pic1}}\qquad
	\subfloat[{\color{blue} Shock-mesh, computational-mesh and surrogate boundaries}]{%
		\includegraphics[height=0.38\textwidth]{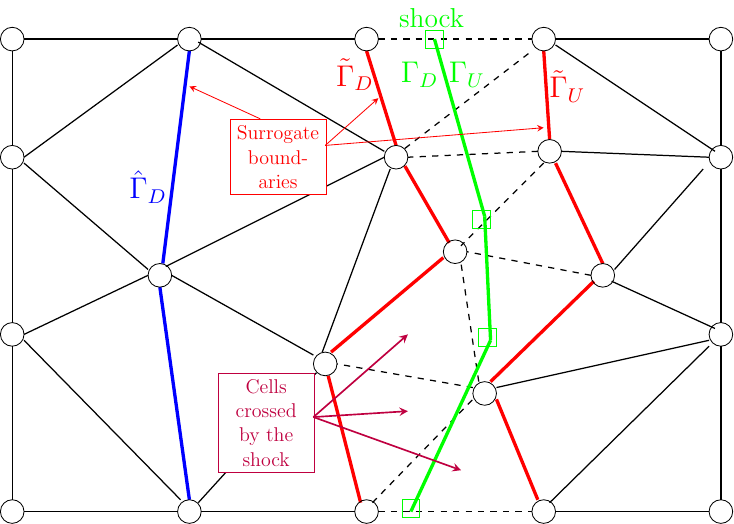}\label{pic2}}
\caption{{\color{blue}The computational-mesh is obtained by removing those cells of the background mesh
that are crossed by the shock-mesh}.}\label{fig:pic2}
\end{center}
\end{figure}

The first step consists in the removal of the triangles crossed by the shock, see Fig.~\ref{pic2}. By doing so, a hole is dug within the background mesh that, contrary to the technique proposed in~\cite{Paciorri1}, is not re-meshed. The creation of the hole splits the background mesh into two disjoint sub-domains which do not include the shock. Instead, we label certain boundaries as ``surrogate'' shock-boundaries which will be used to couple the flow domains, via the shock relations. We shall hereafter call ``computational mesh'' the background mesh with the triangles within the hole being removed. It is worth noting that the number of grid-points of the computational mesh is the same as that of the background grid, whereas the number of triangular cells is less, due to the cell removal. Hereafter, the upstream and downstream surrogate boundaries, drawn using red lines in Fig.~\ref{pic2}, will be called $\tilde{\Gamma}_U$ and $\tilde{\Gamma}_D$. Furthermore, the shock-boundary, which represents the actual shock position, will be referred to as $\Gamma$ and its upstream and downstream sides as $\Gamma_U$ and $\Gamma_D$, respectively. Finally, a second surrogate boundary located within the shock-downstream sub-domain (the blue line in Fig.~\ref{pic2}) will be called $\hat{\Gamma}_D$. This second surrogate boundary is obtained by removing all cells that have one, or more, nodes on $\tilde{\Gamma}_D$.

\subsection{Computation of the tangent and normal unit vectors}\label{sect:vectors}

\begin{figure}[h]
\begin{center}
		\includegraphics[width=0.4\textwidth]{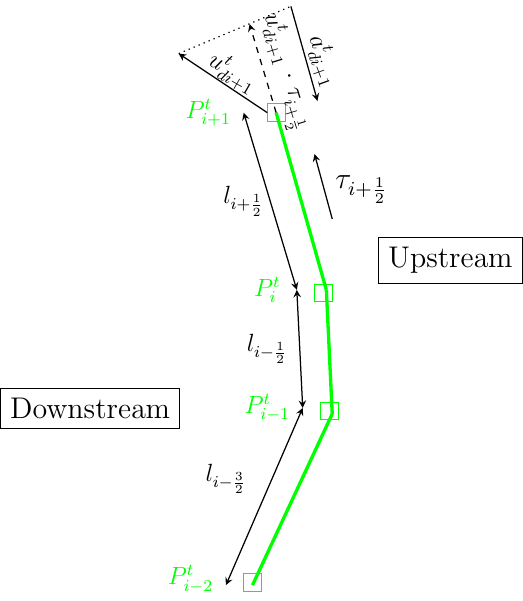}
\end{center}
\caption{Test needed to check whether point $P_{i+1}$ belongs to the domain of dependence of $P_i$.}\label{pic3}
\end{figure}

In order to {\color{green}apply} the Rankine-Hugoniot jump relations, {\color{green}Eq.~(\ref{eq:euler1})}, 
the tangent and normal unit vectors along the shock-front have to be calculated {\color{green}within each pair of shock-points}.
The tangent unit vector $\boldsymbol{\tau}_i$ {\color{green}in shock-point $i$} is obtained from:
\begin{equation}\label{eq2}
	\boldsymbol{\tau}_i\,=\,\frac{\mathbf{v}_{\tau_i}}{\mid \mathbf{v}_{\tau_i} \mid}
\end{equation}
where $\mathbf{v}_{\tau_i}$ is the vector tangent to the shock-front {\color{green}in shock-point $i$}.
The normal unit vector $\mathbf{n}_i$ is 
{\color{green} perpendicular to $\boldsymbol{\tau}_i$} and such that
it points from the shock-downstream towards the shock-upstream region.
The computation of $\mathbf{v}_{\tau_i}$ relies on finite difference formulae which involve the coordinates of 
the shock-point itself and {\color{green}those} of its neighboring shock-points. 
{\color{green}By reference to Fig.~\ref{pic3}, $\mathbf{x}\left(P_i^t\right)$ denotes the position
of shock-point $i$ at time level $t$.}
Shock-points $i-1$ and $i+1$ are located on both sides of shock-point $i$ and their position {\color{green}$\mathbf{x}\left(P_{i-1}^t\right)$
and $\mathbf{x}\left(P_{i+1}^t\right)$} at time level $t$ 
can be used to compute the tangent and normal unit vectors in shock-point $i$. 
A preliminary test is required to verify whether these adjacent shock-points belong to the domain of dependence of shock-point $i$.
This is easily checked using the following inequality:
\begin{equation}\label{eq1}
\mathbf{u}_{d,i+1}^t\,\cdot\,\boldsymbol{\tau}_{i+\frac{1}{2}}\,-\,a_{d,i+1}^t\,<\,0
\end{equation}
where:
\begin{equation}
{\color{green}
\boldsymbol{\tau}_{i+\frac{1}{2}} =
	\frac{\mathbf{x}\left(P_{i+1}^t\right)-\mathbf{x}\left(P_i^t\right)}{l_{i+\frac{1}{2}}} \quad\quad 
	l_{i+\frac{i}{2}} = |\mathbf{x}\left(P_{i+1}^t\right)-\mathbf{x}\left(P_i^t\right)|
	}
\end{equation}
and $\mathbf{u}_{d,i+1}^t$ and $a_{d,i+1}^t$ are the shock-downstream flow and acoustic velocity in shock-point $i+1$ at time level \emph{t}. 
If Eq.~\eqref{eq1} is verified, shock-point $i+1$ falls within the domain of dependence of shock-point \emph{i}.
Once this test has been repeated in shock-point $i-1$, three different situations may arise:
\begin{enumerate}
\item both shock-points $i-1$ and $i+1$ are in the domain of dependence of shock-point \emph{i};
\item only shock-point $i-1$ is in the domain of dependence of shock-point \emph{i};
\item only shock-point $i+1$ is in the domain of dependence of shock-point \emph{i};
\end{enumerate}
When case 1 applies, the computation of $\mathbf{v}_{\tau_i}$ must involve the shock-points on both sides; therefore:
\begin{equation}\label{vtau}
{\color{green}
	\mathbf{v}_{\tau_i}\,=\,\boldsymbol{\tau}_{i+\frac{1}{2}}\,l_{i-\frac{1}{2}}^2\,
	+\,\boldsymbol{\tau}_{i-\frac{1}{2}}\,l_{i+\frac{1}{2}}^2
	}
\end{equation}

When case 2 applies, shock-point $i+1$ must not be used in the computation of the tangent vector $\mathbf{v}_{\tau_i}$,
and the following upwind-biased formula,
{\color{green}
which involves shock-point $i-2$, instead of $i+1$, is used:}
\begin{equation}\label{eq3}
{\color{green}
	\mathbf{v}_{\tau_i}\,=\,\boldsymbol{\tau}_{i-\frac{1}{2}}\,\left(l_{i-\frac{1}{2}}+l_{i-\frac{3}{2}}\right)^2\,+\,
	\left(\boldsymbol{\tau}_{i-\frac{1}{2}}+\boldsymbol{\tau}_{i-\frac{3}{2}}\right)\,l_{i-\frac{1}{2}}^2
	}
\end{equation}
Finally, the third case is specular to the second one, but the corresponding formula involves shock-points \emph{i}, $i+1$ and $i+2$.

{\color{green}
The finite difference approximations~(\ref{vtau}) and~(\ref{eq3})
are both second-order-accurate even if the shock-points are un-evenly spaced along the shock-front.
}

%

\subsection{Solution update using the CFD solver}\label{sect:step of eulf}

\begin{figure}
	\begin{center}
\includegraphics[width=0.6\textwidth]{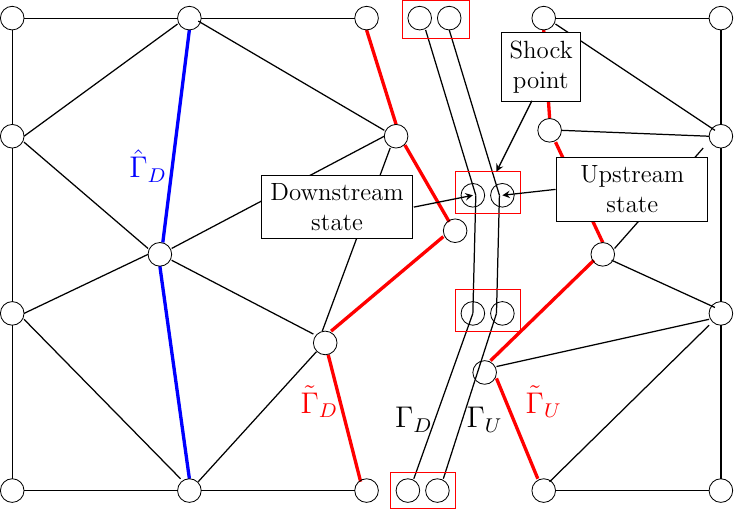}
	\end{center}
	\caption{The solution update is performed {\color{green}using} the computational mesh.}\label{pic4}
	\end{figure}

The solution is updated to time level \emph{t+$\Delta$t} using an unstructured shock-capturing code. The flow solver uses the computational mesh built in step~\ref{sect:cellremoval}, which includes  the  surrogate shocks $\tilde{\Gamma}_U$ and $\tilde{\Gamma}_D$ as part of its boundary. In particular, the flow computations are performed on two non-communicating domains separated by the hole bounded by the surrogate shock-boundaries (see Fig.~\ref{pic4}). As already mentioned, each shock-point {\color{green}consists in} 
two superimposed points  of the shock-mesh: one of these represents the {\color{green}shock-}downstream state and the other 
the {\color{green}shock-upstream one}.
Even though these points are not part of the flow domain, they will play a central role in the coupling of the surrogate shock-boundaries, as we will see in the next sections.\\
Concerning the solver used in this paper, it is based on a  \emph{Residual Distribution  (RD)}  method evolving in time approximations of the values of the flow variables in mesh nodes. The method has several appealing characteristics, including the possibility of defining genuine multidimensional upwind strategies for Euler flows,  by means of a wave decoupling exploiting  appropriately preconditioned forms of the equations \cite{Bonfiglioli1}. By combining ideas from both the stabilized finite element and finite volume methods, these schemes allow to achieve second order of accuracy and monotonicity preservation with a compact stencil of nearest neighbors.
The interested reader can refer to \cite{dr2017,ar2017}  and references therein for an in-depth review of this family of methods, as well as to \cite{Bonfiglioli1,Bonfiglioli2} and references therein for some specific choices of the implementation used here.

Note that the choice of the flow solver is somewhat independent on the rest of the method object of this paper. Concerning the presentation in the following sections,  the main  impact of our choice   is  on the structure of the  solver which is assumed to be evolving nodal values of the unknowns. Cell based  discretization  methods can be easily accommodated by minor modifications of the transfer operators discussed later in the paper.

\subsection{Solution transfer from/to the shock to/from the surrogates}

The flow solver provides updated nodal values within all grid-points of the {\color{green}computational mesh} 
at time level \emph{t+$\Delta$t}. 
The shock-upstream surrogate boundary, $\tilde{\Gamma}_U$ behaves like a supersonic outflow and, therefore, no boundary conditions should 
be applied. The situation along the shock-downstream surrogate $\tilde{\Gamma}_D$ is however different, 
since the flow is subsonic in the shock-normal direction and, therefore, boundary conditions corresponding 
to the downstream-running waves (the `fast' acoustic, entropy and vorticity waves) are missing. 
Moreover, the upstream and downstream states of the shock-points have not been updated, since {\color{green}the shock-mesh is} 
not part of the computational mesh. 
To perform this update, one needs to define appropriate transfer operators from the surrogate shock-boundaries to the shock.  
Due to the use of an upwind discretization in the CFD solver, we assume that the only variable that has been correctly computed 
along the shock-downstream surrogate boundary is 
the Riemann variable associated with the acoustic wave that moves upstream towards the shock:
\begin{equation}\label{eq4}
R_D^{t+\Delta t}\,=\,\tilde{a}_d^{t+\Delta t}\,+\,\frac{\gamma-1}{2}\,\tilde{\mathbf{u}}_d^{t+\Delta t}\,\cdot\,\mathbf{n}
\end{equation}
In Eq.~\eqref{eq4} $\mathbf{n}$ is the shock normal, $\tilde{a}_d^{t+\Delta t}$ is the speed of sound and $\tilde{\mathbf{u}}_d^{t+\Delta t}$ is 
the flow velocity on the shock-downstream side of the shock. It is noted that $R_D^{t+\Delta t}$ is assumed to be correctly computed by 
the CFD solver even if the values $\tilde{a}_d^{t+\Delta t}$ and $\tilde{\mathbf{u}}_d^{t+\Delta t}$ may each be incorrect.

These transfers operators need to be applied twice, once to and once from the shock, and are discussed  below.

\subsection{First transfer: from the surrogate boundaries to the shock}\label{first transfer}

Since the CFD solver uses Roe's parameter vector $\mathbf{Z}$~\cite{Roe1} as the dependent variable,
this is the set of variables used to transfer data between the shock and the surrogate boundaries.
The first transfer is required to update the shock-upstream points on $\Gamma_U$ and to {\color{green}transfer} $R_D^{t+\Delta t}$ from $\tilde{\Gamma}_D$ to $\Gamma_D$.\\
For both transfers, a Taylor series expansion truncated to the second order is used for the extrapolation:
\begin{equation}\label{eq5}
Z_i(\emph{\textbf{x}})\,=\,Z_i(\emph{$\tilde{\textbf{x}}$})\,+\,\nabla Z_i(\emph{$\tilde{\textbf{x}}$})\cdot(\emph{\textbf{x}}\,-\,\emph{$\tilde{\textbf{x}}$})\,+\,o(\|\emph{\textbf{x}}\,-\,\emph{$\tilde{\textbf{x}}$}\|^2)
\end{equation}
where {\color{green}$Z_i$ is any of the four components 
of $\mathbf{Z} = \sqrt{\rho}\left(1,H,u,v\right)^t$, \emph{$\textbf{x}$} and \emph{$\tilde{\textbf{x}}$} are the
coordinates of two different points that belong to \emph{$\Gamma$}, resp.\ \emph{$\tilde{\Gamma}$}},
and \emph{$\nabla Z_i$}(\emph{$\tilde{\textbf{x}}$}) is the gradient computed on the surrogate boundary, {\color{green} using Eq.~\eqref{eq6}}.
Note that, in order to achieve an overall second order of accuracy {\color{green} in the calculation
of $Z_i(\emph{\textbf{x}})$, the approximation of the gradient in Eq.~\eqref{eq5} only needs to be consistent, i.e.\ first-order-accurate}.

The first transfer consists in two phases.
\begin{enumerate}
\item Upstream: from $\tilde{\Gamma}_U$ to $\Gamma_U$\\

The first phase consists in extrapolating from $\tilde{\Gamma}_U$ to $\Gamma_U$.

\begin{figure}
\begin{center}
\includegraphics[width=0.7\textwidth]{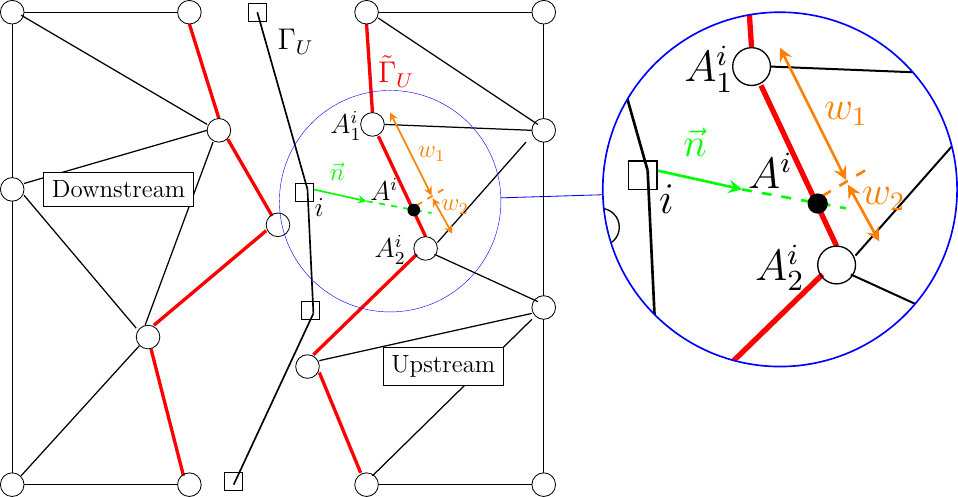}
\end{center}
\caption{{\color{green}Transfer} of variables from the surrogate boundary $\tilde{\Gamma}_U$ to the upstream shock state.}\label{pic8}
\end{figure}

In order to be consistent with the physics of the problem, the transfer of variables takes place along the direction of the shock-normal in the shock-point that has to be updated.
{\color{red}As shown in Fig.~\ref{pic8}, point $A^i$ is the intersection between the surrogate boundary $\tilde{\Gamma}_U$ and the straight-line parallel to the shock-normal in shock-point $i$ which passes through shock-point $i$.}
The  value of the dependent variables (and their gradients) in point $A^i$ is computed using the solution in grid-points $A_1^i$ and $A_2^i$ by means of the following formula:
\begin{equation}\label{eq6}
\phi(A^i)\,=\,\phi(A_1^i) \, w_2 \,+\,\phi(A_2^i) \, w_1 
\end{equation}
where $\phi$ {\color{green}is either $Z_i$ or $\nabla Z_i$}, and $w_1$ and $w_2$ are the weights, 
equal to the normalized distances between $A^i$ and grid-points $A_1^i$ and $A_2^i$. 
{\color{green}
When Eq.~(\ref{eq6}) is used to compute the gradient,
the evaluation of the gradient in the grid-points of the surrogate boundaries}
is performed here using an area-weighted formula, which is reported in Appendix 1. 

{\color{green}
Once the value of $\mathbf{Z}$ in the intersection point $A^i$ has been computed using Eq.~(\ref{eq6})
the value of $\mathbf{Z}$ in shock-point $i$ is computed
by means of Eq.~\eqref{eq5}, having set \emph{$\textbf{x}$} equal to the coordinates of shock-point \emph{i} 
and \emph{$\tilde{\textbf{x}}$} to those of $A^i$}.\\

\item Downstream: from $\tilde{\Gamma}_D$ to $\Gamma_D$:\\

	The Riemann variable defined by Eq.~\eqref{eq4}, which is the only quantity that has been correctly computed on $\tilde{\Gamma}_D$ by the unstructured shock-capturing solver, has to be {\color{green}transferred} from $\tilde{\Gamma}_D$ to $\Gamma_D$ using Eq.~\eqref{eq5}. The procedure is identical to that used for the upstream boundaries: starting from the shock-point to be updated and moving forward along its normal vector as far as an edge of $\tilde{\Gamma}_D$ is intersected in $B^i$ (see Fig.~\ref{pic11}). 

\begin{figure}
\begin{center}
\includegraphics[width=0.7\textwidth]{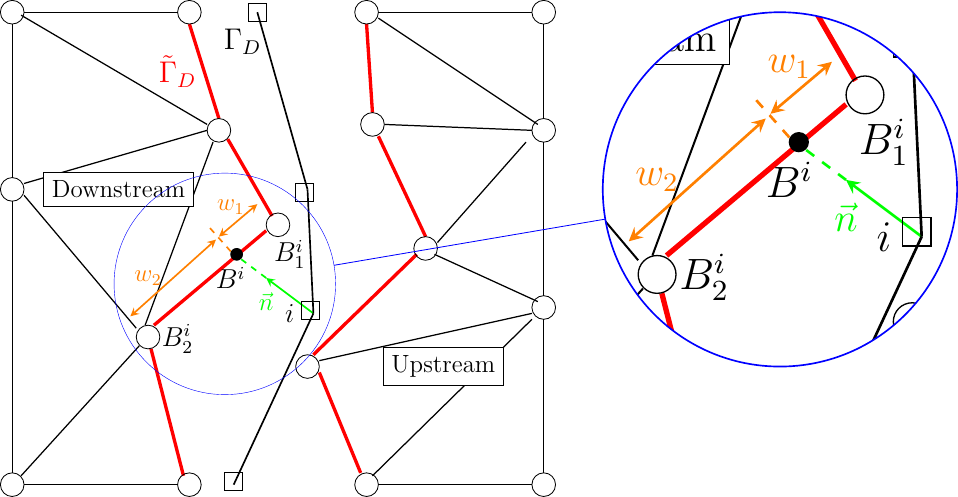}
\end{center}
\caption{{\color{green}Transfer} of the Riemann's variable from the surrogate boundary $\tilde{\Gamma}_D$ to the shock point.}\label{pic11}
\end{figure}
\end{enumerate}

\subsection{Shock calculation}\label{sect:shock}

As already mentioned and also schematically shown in Fig.~\ref{pic4}, each shock-point consists in two superimposed grid-points, which represent the shock-upstream and the shock-downstream states. The velocity component in the shock-normal direction, $w$, is also stored within each shock-point. For the reasons explained in step~5.5, the shock-upstream state and the Riemann variable $R_D$, Eq.~\eqref{eq4}, on the shock-downstream side of the shock have been correctly updated at time level $t+\Delta t$. 
The shock-downstream state ($\rho_d$, $p_d$ and $\mathbf{u}_d$) and the shock-speed at time $t+\Delta t$, which are yet unknown at this stage, 
are found by solving a system of five algebraic non-linear equations. The first four equations are the Rankine-Hugoniot jump relations
and the fifth is Eq.~\eqref{eq4}: 
\begin{equation}\label{RH}
	{\color{green}
	\begin{array}{ccl}
	\rho_d^{t+\Delta t}(\mathbf{u}_d^{t+\Delta t}\cdot\mathbf{n}-\,w)&=&\rho_u^{t+\Delta t}(\mathbf{u}_u^{t+\Delta t}\cdot\mathbf{n}-w) \\
	\rho_d^{t+\Delta t}(\mathbf{u}_d^{t+\Delta t}\cdot\mathbf{n}-\,w)^2+p_d^{t+\Delta t}&=&\rho_u^{t+\Delta t}(\mathbf{u}_u^{t+\Delta t}\cdot\mathbf{n}-w)^2+p_u^{t+\Delta t}\\
	\frac{\gamma}{\gamma-1}\frac{p_d^{t+\Delta t}}{\rho_d^{t+\Delta t}} + \frac{1}{2}(\mathbf{u}_d^{t+\Delta t}\cdot\mathbf{n}-w)^2&=&\frac{\gamma}{\gamma-1}\frac{p_u^{t+\Delta t}}{\rho_u^{t+\Delta t}}+ \frac{1}{2}(\mathbf{u}_u^{t+\Delta t}\cdot\mathbf{n}-w)^2 \\
	\mathbf{u}_d^{t+\Delta t}\cdot\boldsymbol{\tau}&=&\mathbf{u}_u^{t+\Delta t}\cdot\boldsymbol{\tau}\\
        R_D^{t+\Delta t}&=&\tilde{a}_d^{t+\Delta t}\,+\,\frac{\gamma-1}{2}\,\tilde{\mathbf{u}}_d^{t+\Delta t}\,\cdot\,\mathbf{n}
	\end{array}
	}
\end{equation}
Hence, for system~\eqref{RH}, the vector of known variables ($\rho_u$, $p_u$, $\mathbf{u}_u$, $R_D$) is then used to 
find the updated values of the unknowns ones ($\rho_d$, $p_d$, $\mathbf{u}_d$, $w$).
The system~\eqref{RH} is solved within each shock-point using the Newton-Raphson root-finding algorithm, 
thus providing the correct downstream state and shock-speed at time level \emph{t+$\Delta$t}.

\subsection{Second transfer: from the shock to the surrogate boundaries}\label{second transf}

Once the shock-downstream states along the shock-boundary $\Gamma_D$ have been updated as described in step~\ref{sect:shock}, the grid-points on the downstream surrogate boundary $\tilde{\Gamma}_D$ need also to be updated.	\\

Downstream: from $\Gamma_D$ to $\tilde{\Gamma}_D$\\

\begin{figure}
\begin{center}
\includegraphics[width=0.7\textwidth]{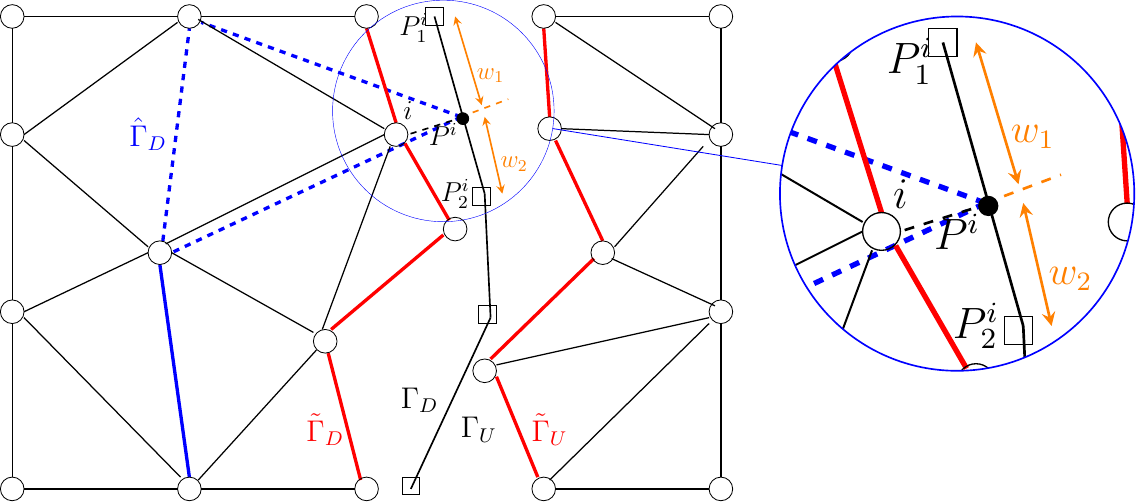}
\end{center}
\caption{Search of the auxiliary point $P^i$ on the shock poly-line used to interpolate the dependent variables in the grid-points of the surrogate boundary $\tilde{\Gamma}_D$.}\label{pic9}
\end{figure}

The first step needed to update grid-point \emph{i} on $\tilde{\Gamma}_D$ consists in finding the projection $P^i$ of grid-point $i$ 
on the shock-poly-line. 
{\color{red}This is accomplished by first locating the closest shock-edge to grid-point $i$
and then projecting along the direction which is
the weighted average\footnote{the weights depend upon the normalized distance between the two shock-points 
$P^i_1$ and $P^i_2$ and grid-point $i$.}
of the two vectors normal to the shock in $P^i_1$ and $P^i_2$}.
Then, the dependent variables in $P^i$ are computed using Eq.~\eqref{eq6}, 
\textcolor{red}{the weights $w_1$ and $w_2$ being the normalized distances of $P^i$ from
shock-points $P^i_1$ and $P^i_2$.}

{\color{green}
The second step consists in using point $P^i$ and two grid-points that belong to the second
surrogate boundary $\hat{\Gamma}_D$ to build a triangle 
(shown using a dashed blue line  in Fig.~\ref{pic9}),which contains grid-point $i$. 

Finally, the dependent variables in grid-point $i$ are linearly interpolated within that triangle}.

The reason for using a second surrogate boundary on the downstream side of the shock lies in fact that, 
whenever the shock-downstream flow is subsonic, the acoustic waves spread in all directions.
Under this circumstance, only grid-points (such as those on $\hat{\Gamma}_D$) that are surrounded {\color{green}on all sides}
by cells have been correctly updated by the CFD solver. 

\subsection{Shock displacement}

\begin{figure}[t]
\begin{center}
\includegraphics[width=0.7\textwidth]{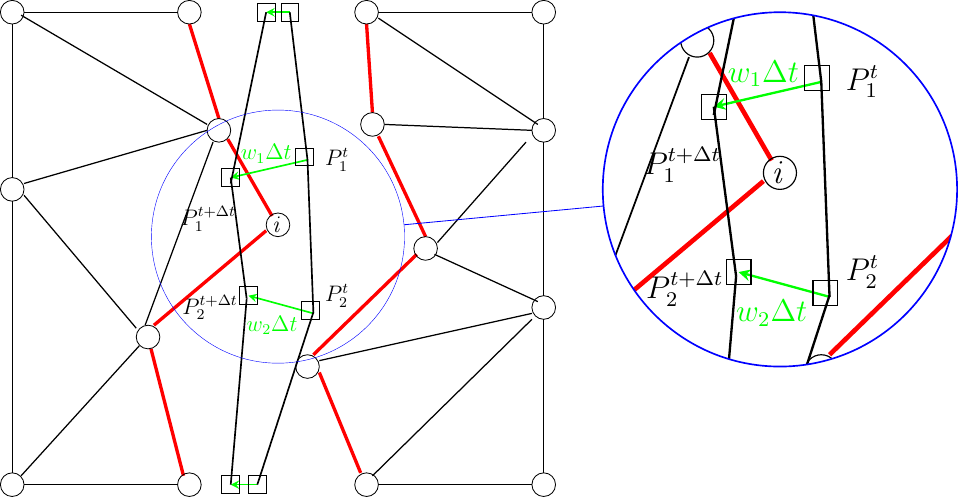}
\end{center}
\caption{The shock-front overtakes a grid-point of the background mesh during its motion.}\label{pic6}
\end{figure}

The new position of the shock-front at time level \emph{t+$\Delta$t} is computed by displacing all shock-points using the following first-order-accurate (in time) formula:	
\begin{equation}\label{displ}
	{\color{green}
	\mathbf{x}\left(P^{t+\Delta t}\right)\,=\,\mathbf{x}\left(P^t\right)\,+\,w^{t+\Delta t}\,\mathbf{n}\,\Delta t
	}
\end{equation}
where $\mathbf{x}\left(P\right)$ denotes the geometrical location of the shock-points.
The use of a first-order-accurate temporal integration formula in~\eqref{displ} is immaterial as long as steady flows are of interest. In this case, a second-order accurate representation of the shock-shape is guaranteed by the use of second-order-accurate formulae to compute the shock-normal, as described in step~\ref{sect:vectors}. For unsteady flows, second-order-accurate time integration
formulae should be used, as done for example in \cite{Paciorri4,Campoli2017}.
As can be seen from Fig.~\ref{pic6}, the shock-front can freely float over the background triangulation and, while doing so, it may cross the downstream surrogate boundary. This is the situation sketched in Fig.~\ref{pic6}, where the shock-fronts at time level $t$ and $t+\Delta t$ have both been drawn. In the sketch of Fig.~\ref{pic6}, grid-point $i$ has been overtaken by the moving shock-front. Whenever this happens, the flow state within grid-point $i$ should be changed accordingly. This is the task performed in the next step.

\subsection{Re-interpolation of nodes crossed by the shock}\label{re-interp}

This step of the algorithm consists in the interpolation of those grid-points of the background mesh that have been overtaken by the  shock-front, thus passing from one region to the other. In order to understand whether grid-point $i$ has been overtaken or not by the shock, the position of the closest shock-edge, before and after the displacement, i.e.\ at time $t$ and $t+\Delta t$, is used to build a quadrilateral, as shown in Fig.~\ref{pic6}. If grid-point \emph{i} falls inside the quadrilateral of vertices $P_1^t$, $P_2^t$, $P_1^{t+\Delta t}$ and $P_2^{t+\Delta t}$, grid-point $i$ has been overtaken and its state has to be updated. 
The state of grid-point $i$ is updated using an interpolation procedure  similar to that illustrated in step~\ref{second transf} .

\subsection{Shock-point re-distribution}

\begin{figure}[t]
\begin{center}
\includegraphics[width=0.4\textwidth]{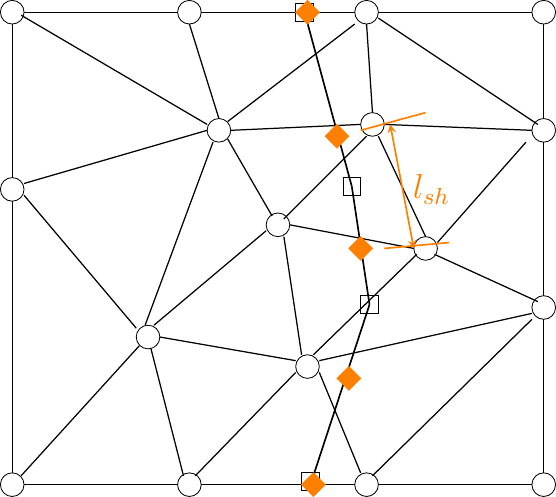}
\caption{Shock-point location along the  shock front: before (squares) and after (diamond) the re-distribution step.}\label{pic7}
\end{center}
\end{figure}

During the shock displacement step, the shock-edges may stretch or shorten, depending on the relative motion of the various shock-point that make up the shock-front. This might lead to a shock-poly-line made of shock-edges whose length is considerably different from the local size of the background mesh. To avoid such a risk, a shock-point re-distribution can be performed as the last step of the algorithm. Doing so, it is possible to ensure that the shock-edge lengths are approximately equal to the edges of the underlying background mesh. A naive shock-point re-distribution procedure is done by imposing that all shock-edges have the same fixed length \emph{$l_{sh}$}, preset by the user. 
Whenever the shock-points are re-located along the  shock-front, both the shock-upstream and shock-downstream state within each shock-point have to be re-computed, a task which is easily accomplished using linear interpolation along the shock-front.
Figure~\ref{pic7} shows the location of the shock-points along the  shock-front both before and after the re-distribution. \\
At this stage, the numerical solution has been correctly updated at time level \emph{t+$\Delta$t}.


	\section{Numerical results}
	\label{sec:tcases}
%
{\color{green}All physical quantities displayed in this section have been made dimensionless
using the following set of reference variables: $L$, $\rho_{\infty}$, $\mathbf{u}_{\infty}$, where $L$
is a length scale, $\rho_{\infty}$ and $\mathbf{u}_{\infty}$ the free-stream density and
flow speed. Using the aforementioned set of reference variables, the reference pressure
is twice the free-stream dynamic pressure: $\rho_{\infty} \left(\mathbf{u}_{\infty}\cdot\mathbf{u}_{\infty}\right)$.
}
%
%
%
\subsection{Quasi-one-dimensional nozzle flow}\label{Q1D}
The quasi-one-dimensional (Q1D) steady flow through a variable area duct (converging-diverging nozzle) turns 
out to be particularly well suited as a validation case because 
the flow is non-uniform both upstream and downstream of the shock and
an analytical solution is available, which allows to compute
the discretization error, $\epsilon$, i.e.\ the difference between the exact and computed solutions. 
Moreover, a similar study reported in~\cite{bonhaus1998higher}, showed that the discretization error
within the entire shock-downstream region exhibits
first-order convergence as the grid is refined even if high-order-accurate schemes are used. This is a known deficiency
of shock-capturing schemes which we will show does not affect the eST method.

The nozzle geometry has been taken from~\cite{bonhaus1998higher}:
\begin{equation}
	A/A_{*} = 1 +(A_e/A_{*}-1)(x/L)^2 \quad\quad\mbox{where}\quad -1/2 \le x/L \le 1
\end{equation}
and the exit-to-sonic area ratio is equal to $A_e/A_{*} = 2$.
Having set the ratio between the exit-static to inlet-total pressures equal to $p_{out}/p_{in}^0 = 0.7362$, a steady normal shock
occurs in the diverging part of the nozzle at about {\color{green}$x_{sh}/L = 0.75$}.

In order to simplify the treatment of the boundary conditions,
the left boundary of the computational domain
has been set at {\color{green}$x_{left}/L = 0.05$}, just downstream of the troath, where a supersonic inflow boundary condition applies.

The main advantage of the present test-case is the fact that it has an analytical solution. In particular, 
the Mach number distribution follows from the so-called area-rule:
\begin{equation}\label{eq14}
	\frac{1}{M}\bigg[\frac{2}{\gamma+1}\bigg(1\,+\,\frac{\gamma-1}{2}M^2\bigg)\bigg]^{\frac{\gamma+1}{2(\gamma-1)}}
\,=\,\frac{A}{A^*}
\end{equation}
\begin{figure}[t]	
	\begin{center}
		\includegraphics[width=0.65\textwidth]{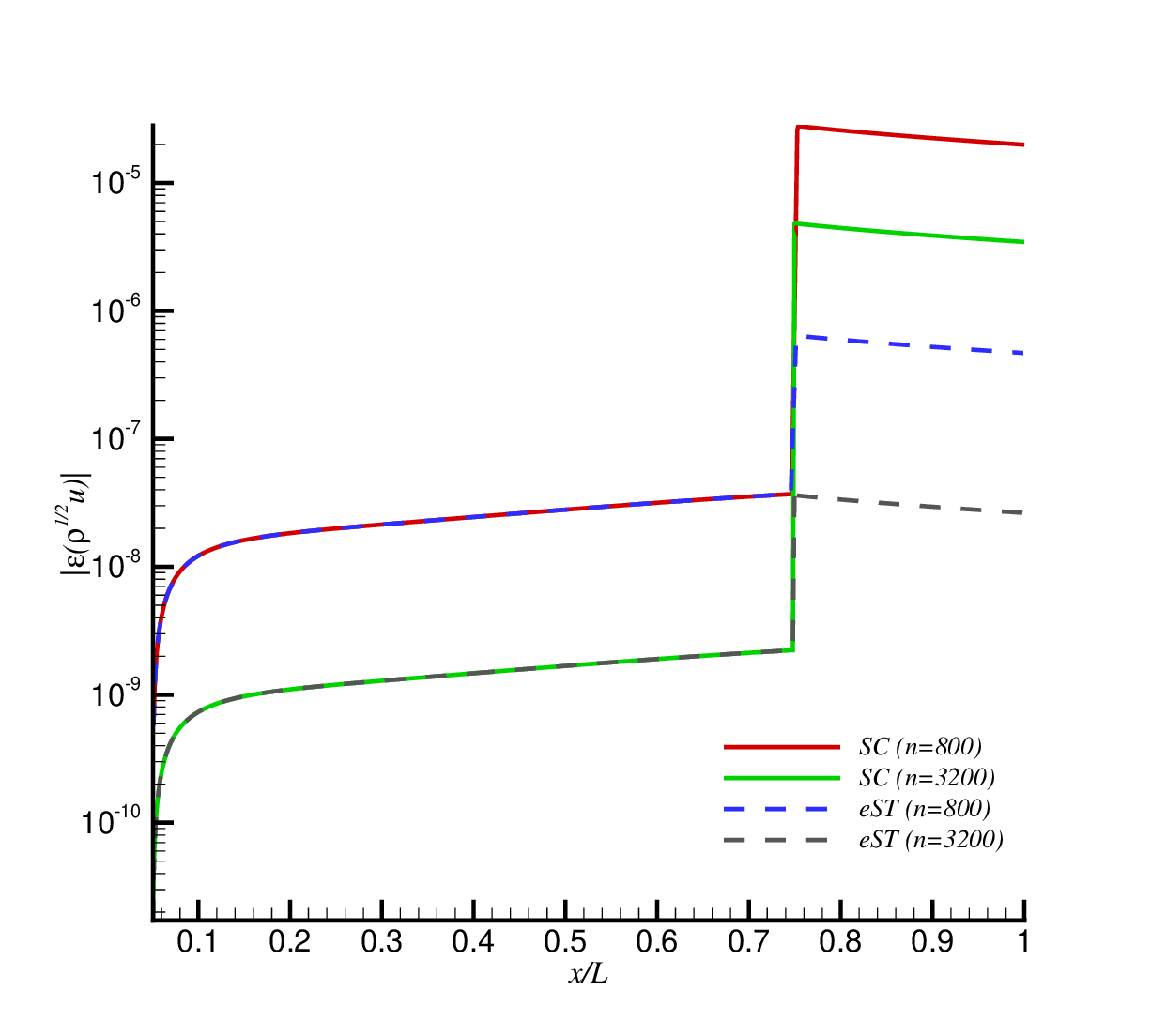}
		\caption{Q1D nozzle flow: pointwise errors analysis for SC and eST.} \label{pic17}
	\end{center}
\end{figure}
A comparison has been made between the SC and eST simulations using a sequence of uniformly spaced grids, with
grid densities ranging between 800 and 6400 cells, see Tab.~\ref{tab1}.
{\color{green}
\begin{table}
\centering
\caption{Q1D nozzle flow: characteristics of the background meshes used to perform the grid-convergence tests.}\label{tab1}
\begin{tabular}{crc} \hline\hline
	Grid level & Cells & $h$ \\ \hline
	0 & 800  & $1.875\,10^{-3}$ \\
	1 & 1600 & $9.375\,10^{-4}$ \\
	2 & 3200 & $4.688\,10^{-4}$ \\ 
	3 & 4800 & $3.125\,10^{-4}$ \\
	4 & 6400 & $2.344\,10^{-5}$ \\ \hline\hline
\end{tabular}
\end{table}
}

Knowledge of the exact solution allows to compute (rather than estimate) the discretization 
error and, therefore, to perform reliable convergence tests.
Figure~\ref{pic17} shows the pointwise distribution of the discretization error for the third component
of Roe's parameter vector. Note that the y-axis of Fig.~\ref{pic17} is in logarithmic scale.
As expected, the SC and eST simulations feature the same discretization error in the entire supersonic, shock-upstream region.
Downstream of the shock, however, SC incurs in a discretization error which is about two orders of magnitude larger
than that of the eST.
\begin{figure}	
	\begin{center}
		\includegraphics[width=0.65\textwidth]{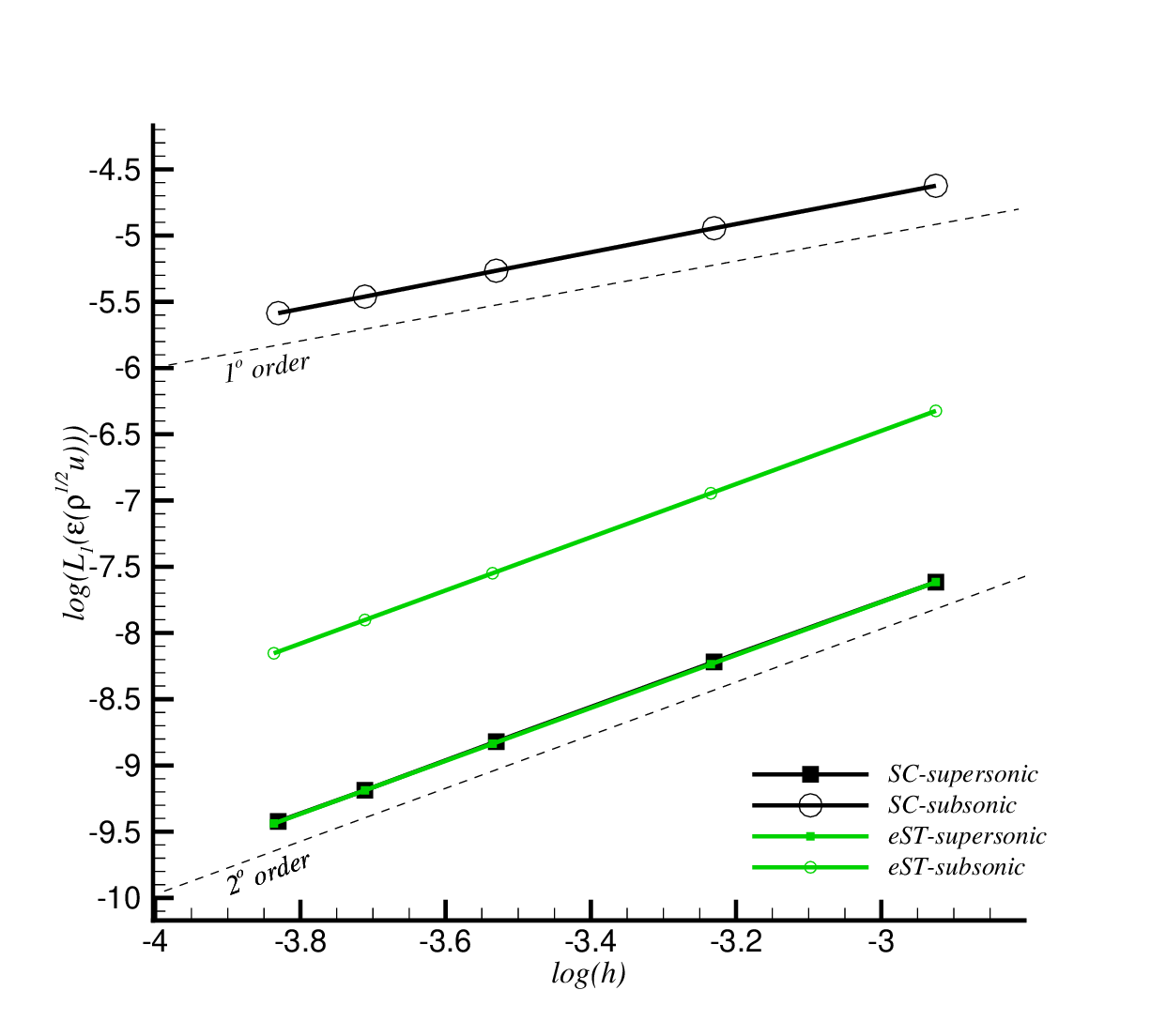}
		\caption{Q1D nozzle flow: global measures of the discretization error for SC and eST.} \label{pic19}
	\end{center}
\end{figure}
The accuracy degradation incurred by SC within the entire downstream region is further confirmed in Fig.~\ref{pic19},
which shows, in a log-log scale, the $L_1$ norm of the discretization error plotted against the mesh spacing.
In contrast to the pointwise measure displayed in Fig.~\ref{pic17}, Fig.~\ref{pic19} shows a
global measure, which has been separately computed within the shock-upstream ($x_L \le x < x_{sh}$) and 
shock-downstream ($x_{sh} < x \le L$)  sub-domains.
The difference between the error-reduction trends exhibited by SC and eST is striking:
the two different shock-modeling practices behave identically within the shock-upstream sub-domain, where both
exhibit second-order convergence as the mesh is refined; within the shock-downstream sub-domain, however,
eST retains second-order convergence, whereas SC drops to first-order.
%
%
\subsection{Planar source flow}\label{planarflow}

\begin{figure*}[h!]
	\subfloat[Sketch of the computational domain.]{%
		\includegraphics[width=0.5\textwidth,trim={3.2cm 2.9cm 3cm 2cm},clip]{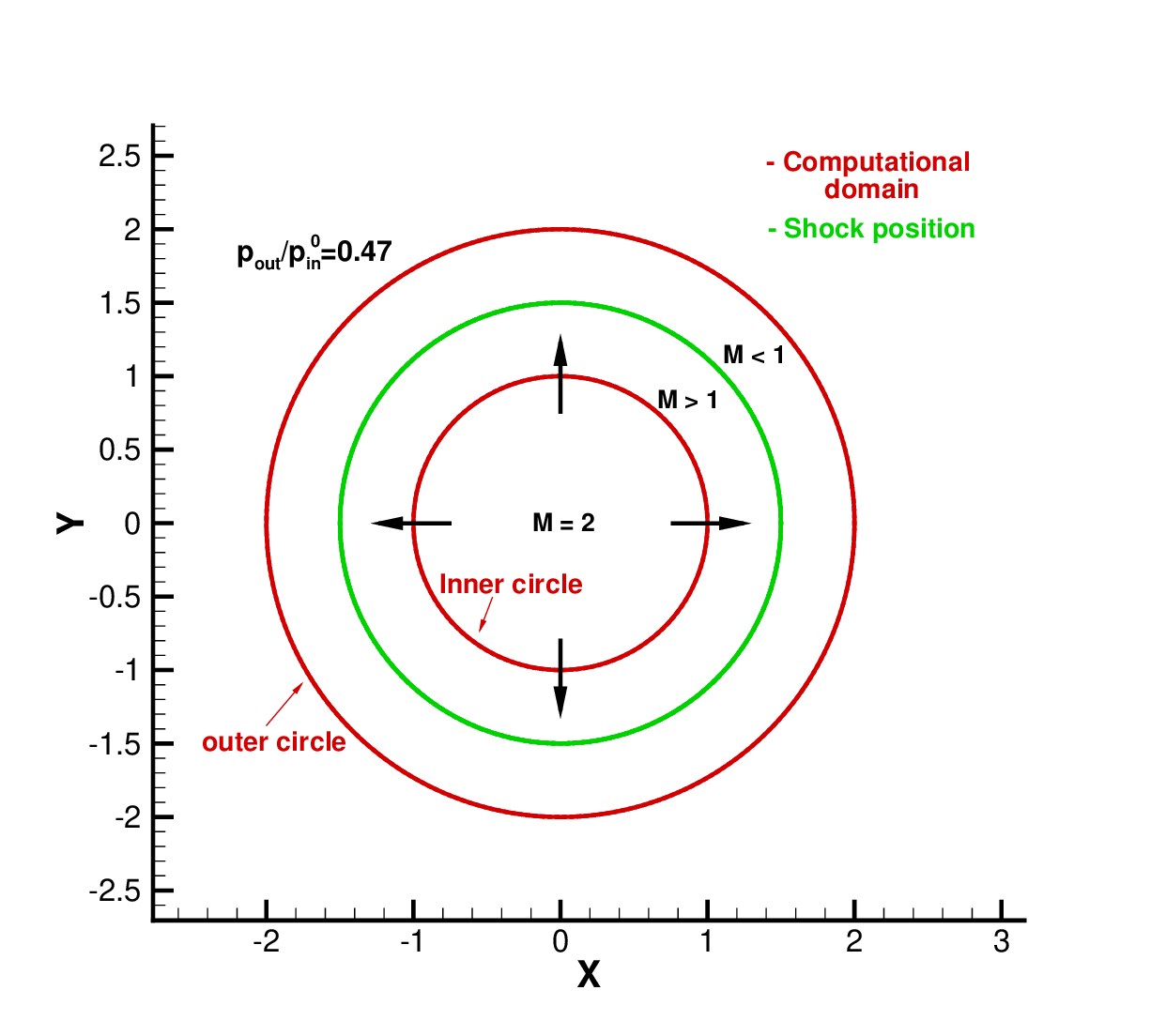}\label{pic20a}}
	\subfloat[Detail of the unstructured mesh employed for the simulations.]{%
		\includegraphics[width=0.5\textwidth]{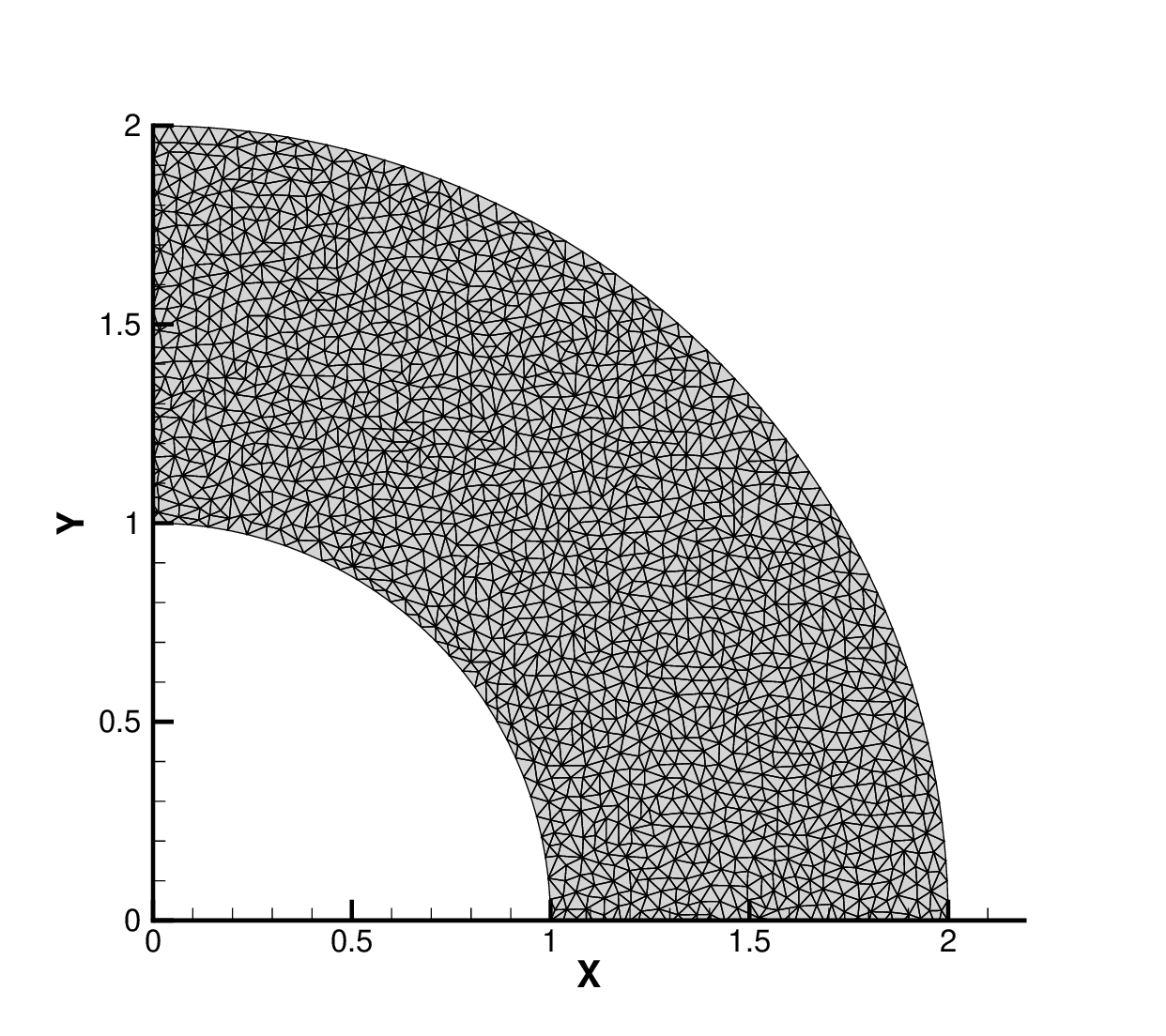}\label{pic20b}}
	\caption{Planar, transonic source flow.}\label{pic20}
\end{figure*}


This test case consists in a compressible, planar source flow 
that has already been studied in~\cite{Paciorri6,Campobasso2011} as a validation case, due to the availability
of an analytical solution.
Indeed, assuming that the analytical velocity field has a purely radial velocity component, it may be easily verified that the 
governing PDEs, written in a polar coordinate system, become identical to those governing a 
compressible quasi-one-dimensional variable-area flow~\eqref{eq14}, provided that the nozzle area varies linearly with the radial distance, $r$,
from the pole of the reference frame.
The computational domain consists in the annulus sketched in Fig.~\ref{pic20a}:
the ratio between the radii of the outer and inner circles {\color{green}($L = r_{in}$)} has been set equal to $r_{out}$/$r_{in}$ = 2. 
A transonic (shocked) flow has been studied by imposing a supersonic inlet flow at $M = 2$ on the inner circle and 
a ratio between the outlet static and inlet total pressures $p_{out}/p_{in}^0 = 0.47$ such that the shock forms at $r_{sh}$/$r_{in}$ = 1.5. 
The Delaunay mesh shown in Fig.~\ref{pic20b}, 
which contains 6916 grid-points and 13456 triangles, has been
generated using \texttt{\ttfamily triangle}~\cite{Triangle1,Triangle2} in such a way that no general alignment is present between 
triangle edges and shock, thus making the discrete problem truly two-dimensional. 
\begin{figure}[tb!]
	\subfloat[]{%
		\includegraphics[width=0.5\textwidth]{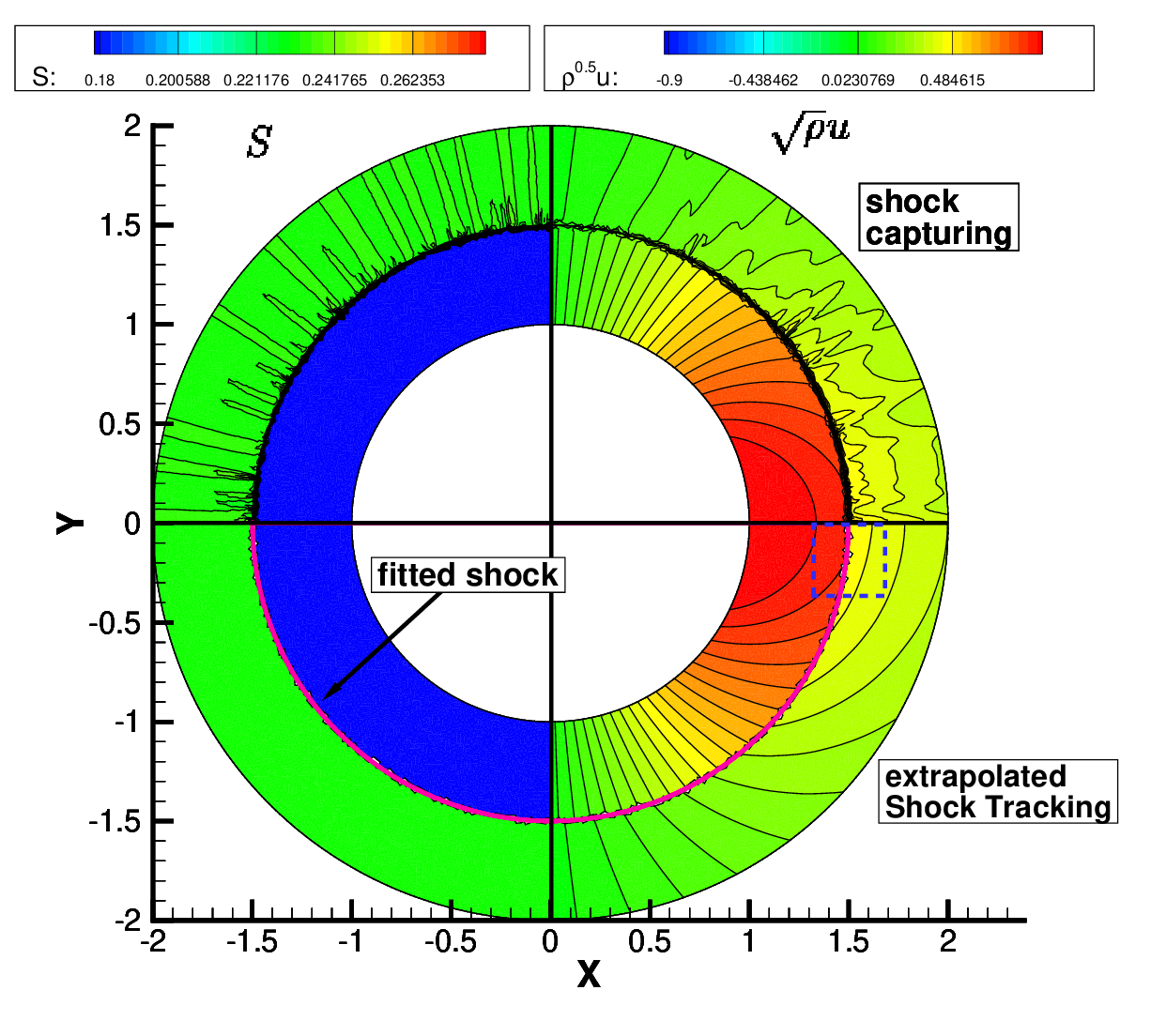}\label{pic22a}}
  \hspace{\fill}
	\subfloat[]{%
		\includegraphics[width=0.5\textwidth]{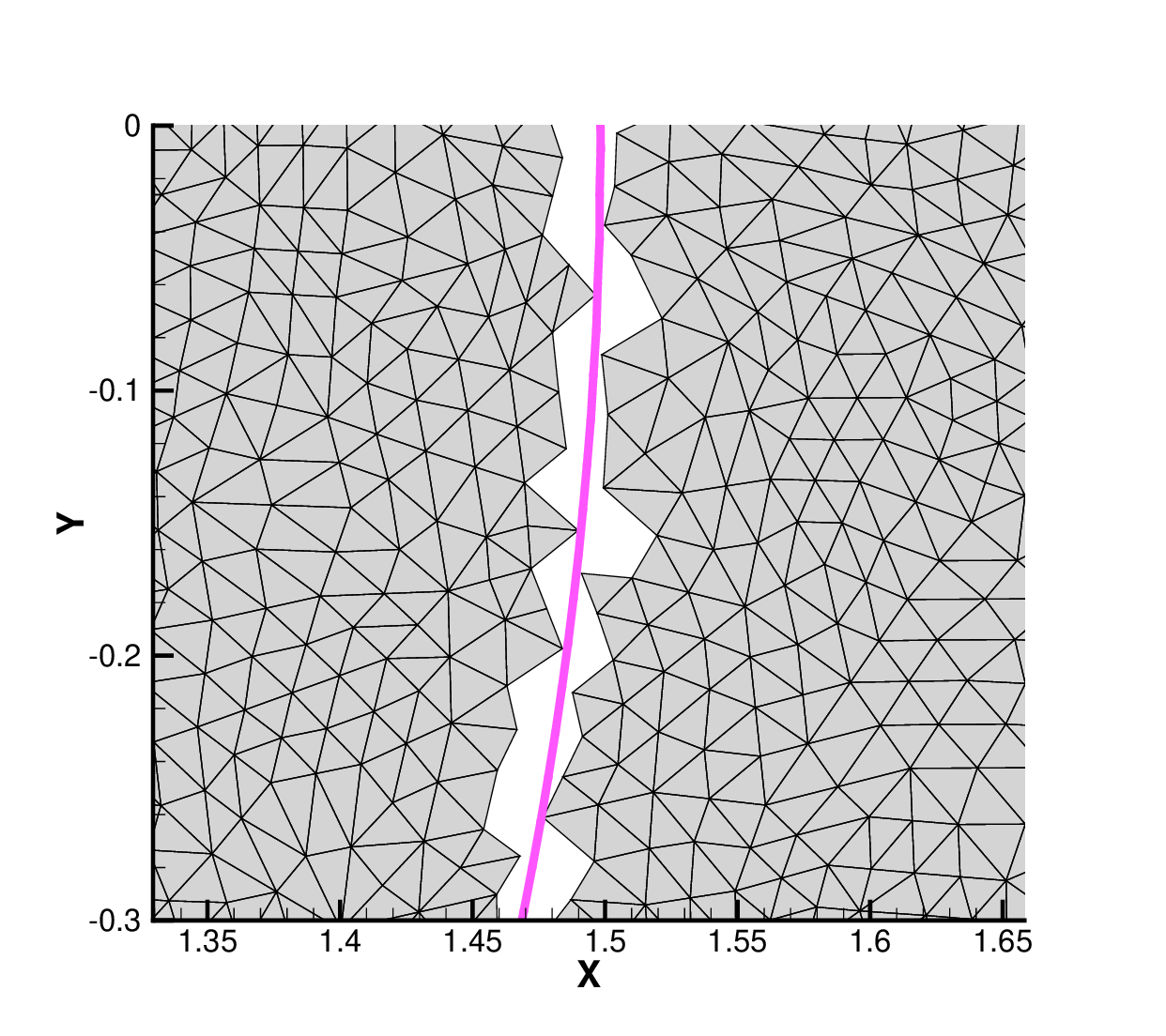}\label{pic22c}}
	\caption{Planar source flow. (a): entropy (left half of the frame) and third component ($\sqrt{\rho}u$) of $\mathbf{Z}$ (right half of the frame). SC result on the top, eST result on the bottom. (b): close-up of the blue square drawn in frame (a).}\label{pic22}
\end{figure}
Figure~\ref{pic22} shows a comparison between the SC and eST solutions, both in terms of entropy, $S = p\rho^{-\gamma}$, and
$\sqrt{\rho}u$ iso-lines.
Both flow variables clearly reveal that the SC solution is
plagued by severe spurious errors due to the misalignment between the edges of the mesh and the captured shock.
These errors propagate in the shock-downstream region, as it is evident from the entropy field of the SC calculation, compromising  the quality of the solution. Note that across the numerical shock layer obtained with the SC approach, the direction  
of the velocity vector is undefined and largely dependent on the  mesh topology. 
This explains the perturbations observed in the shock-downstream region.

\begin{figure}	
	\begin{center}
		\includegraphics[width=0.6\textwidth]{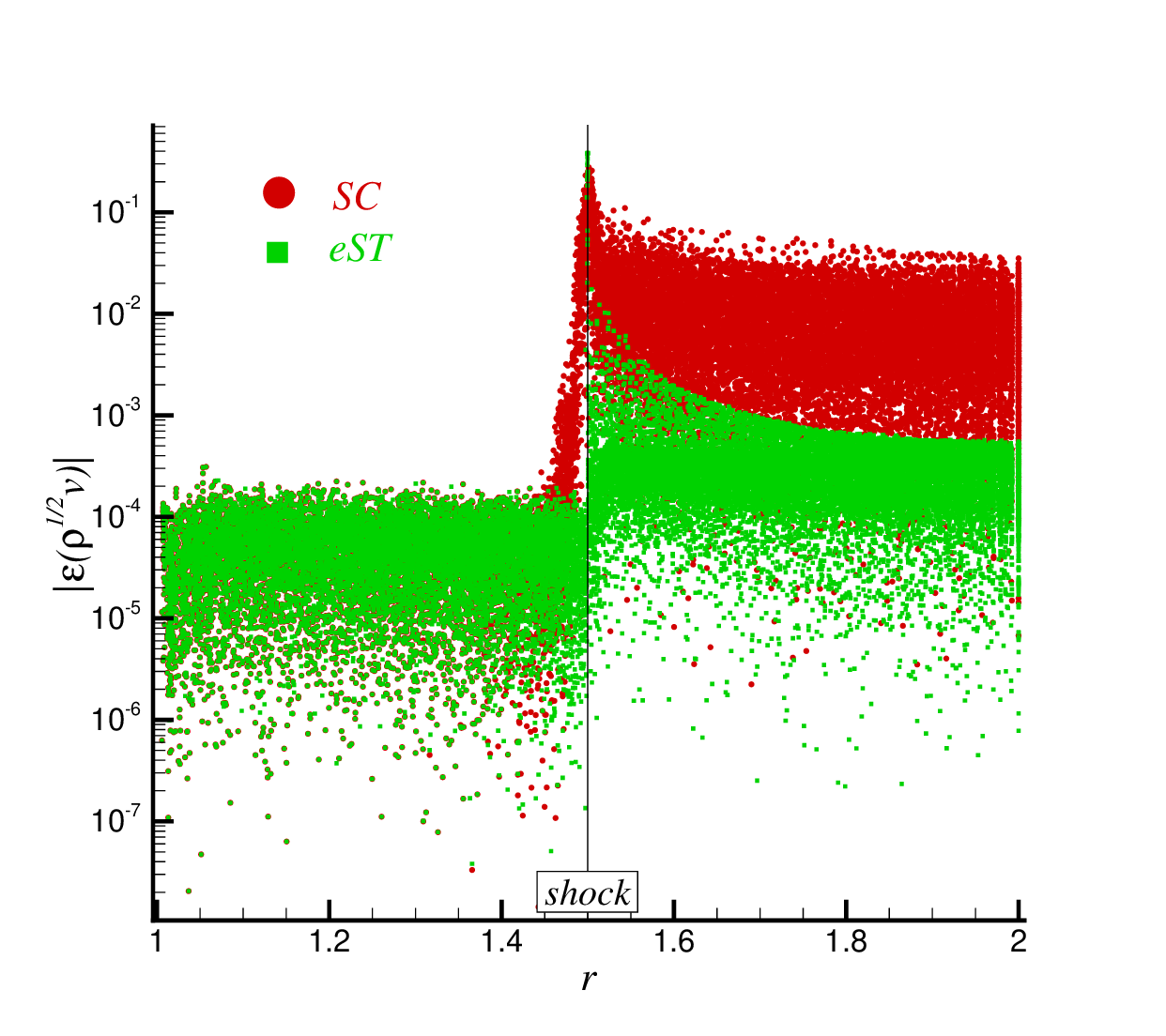}
		\caption{Planar, transonic source flow: pointwise error analysis (SC vs. eST).} \label{pic28}
	\end{center}
\end{figure}

Thanks to the availability of the analytical solution,
a point-wise error analysis has been carried out
by computing the discretization error.
Figure~\ref{pic28} shows the behavior of the local discretization error in all points of the mesh,
plotted against the radial distance from the center of the circle. 
\textcolor{red}{The vertical line drawn in Fig.~\ref{pic28} points to the position ($r=1.5$) where the shock-wave takes place.}  
It can be seen that upstream of the shock 
the error of the SC and eST solutions is equal.
Downstream of the shock ($r > 1.5$), however,
the eST solution exhibits an error which is one or two orders of magnitude lower 
than that obtained using SC.

\begin{figure}	
	\begin{center}
		\subfloat[$\sqrt{\rho}H$]{\includegraphics[width=0.5\textwidth]{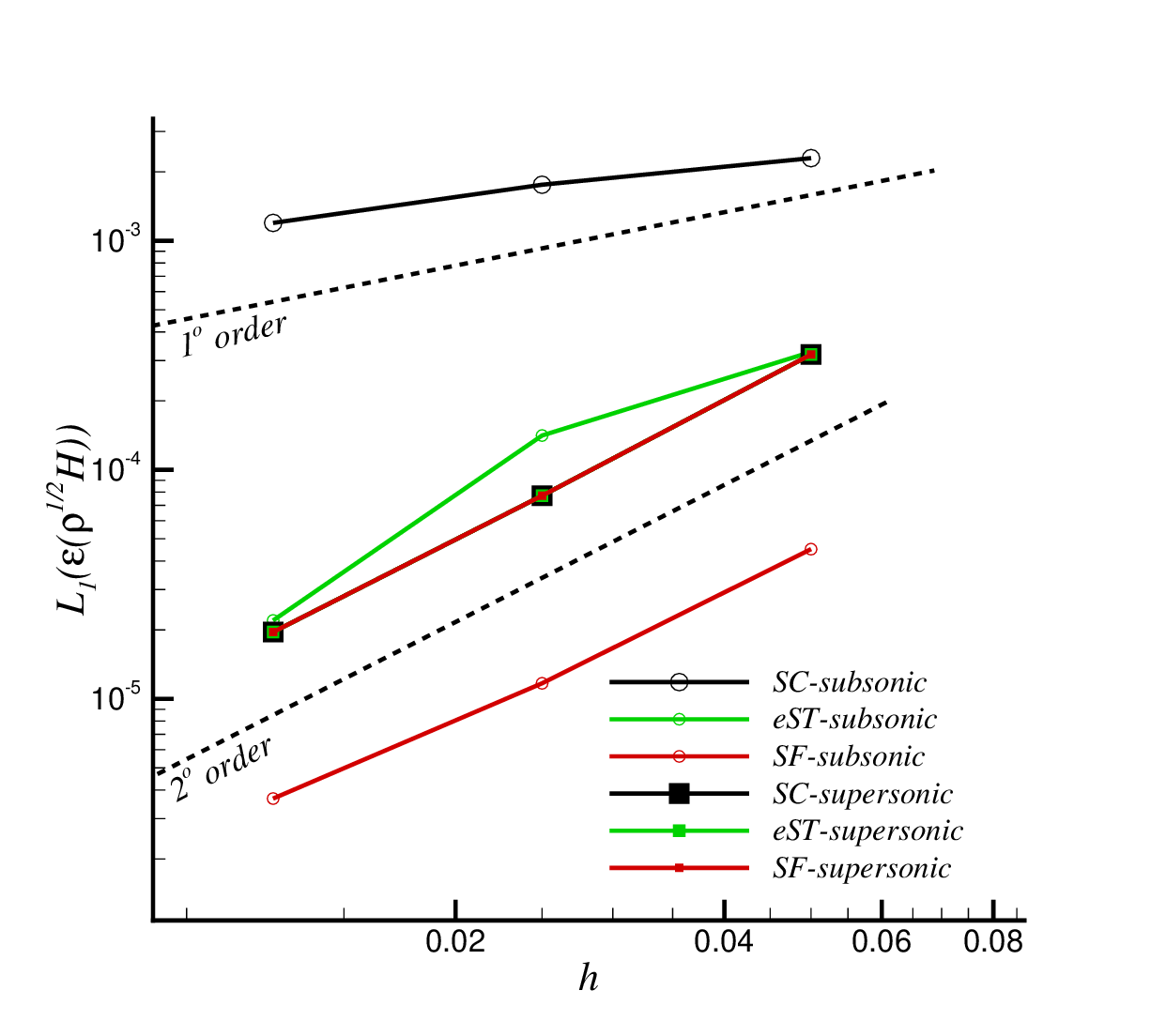}}
		\subfloat[$\sqrt{\rho}v$]{\includegraphics[width=0.5\textwidth]{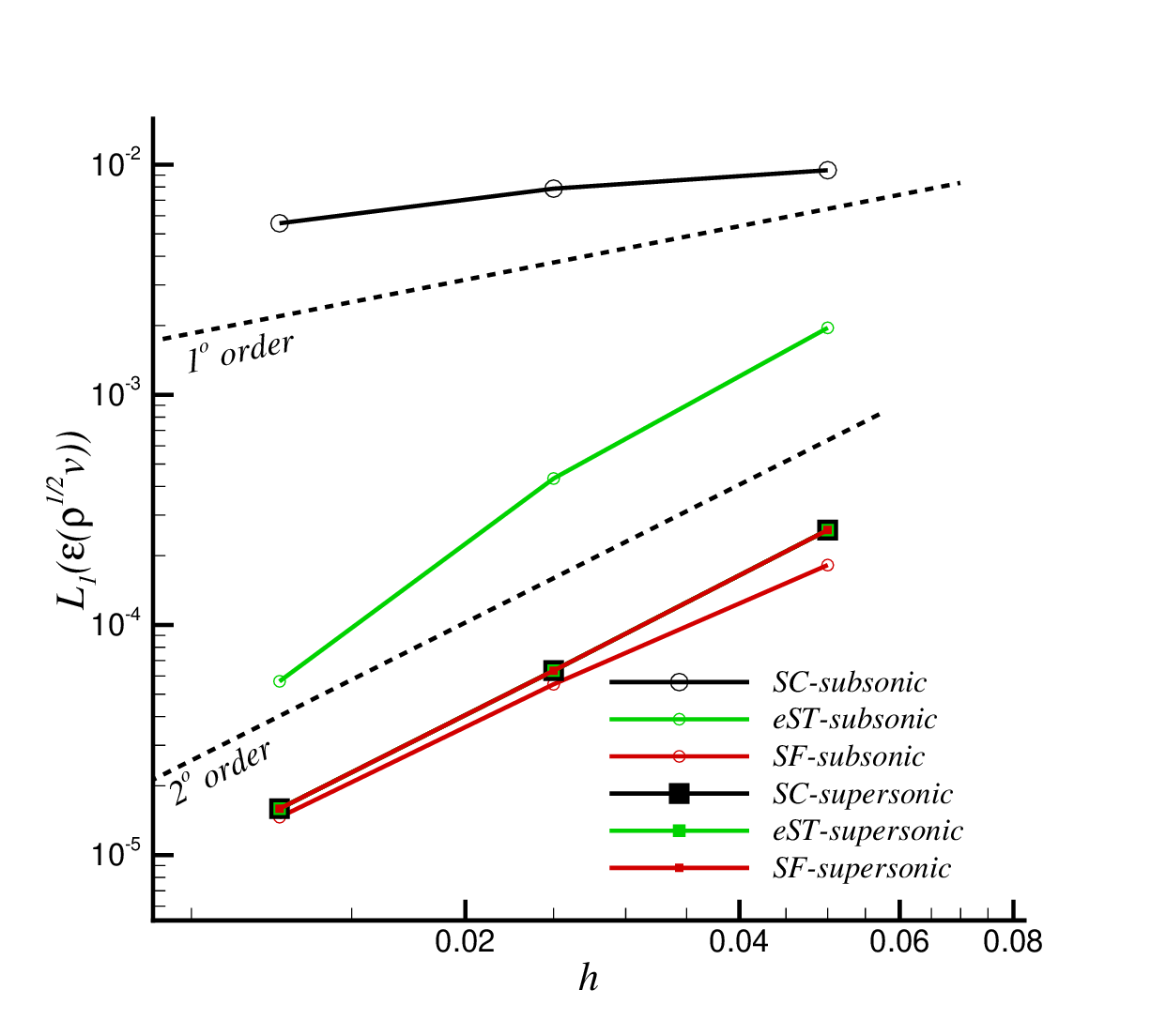}}
		\caption{Planar, transonic source flow: {\color{red}order-of-convergence comparison among SC, eST and SF w.r.t.\ the second and fourth components of the parameter vector $\mathbf{Z}$.}} \label{pic29}
	\end{center}
\end{figure}

An order-of-convergence analysis, similar to that of Sect.~\ref{Q1D}, has also been performed by repeating the same calculation 
on three nested triangulations
{\color{green}whose features are summarized in Tab.~\ref{tab2}, 
where $h$ is the mesh spacing along the inner and outer circular boundaries}.
The coarsest mesh is the one shown in Fig.~\ref{pic20b} and the two finer meshes have been obtained by recursively subdividing each
triangle of the parent mesh into four nested triangles.
\begin{table}
\centering
\caption{Planar source flow: characteristics of the background meshes used to perform the grid-convergence tests.}\label{tab2}
\begin{tabular}{crcr} \hline\hline
	Grid level & Grid-points & Triangles & $h$\\ \hline
	0& 6,916   & 13,456 & 0.05\\
	1& 27,288  & 53,824 & 0.025\\
	2& 108,400 & 215,296& 0.0125\\
\hline\hline
\end{tabular}
\end{table}

A {\em global} measure of the discretization error has been computed using the $L_1$-norm of $\epsilon\left(\sqrt{\rho}v\right)$ and $\epsilon\left(\sqrt{\rho}H\right)$,
separately within the shock-upstream and shock-downstream sub-domains; the results are displayed in Fig.~\ref{pic29}
{\color{red} and also include those published in~\cite{Paciorri6}, which have been obtained using the 
unstructured shock-fitting technique developed by some of the authors in~\cite{Paciorri1}.}
It is noted that within the supersonic, shock-upstream region, {\color{red}the three numerical solutions feature the same discretization error}
and converge to the exact solution at design (second) order as the mesh is refined.
Downstream of the shock, however, {\color{red}only the two shock-fitting techniques (SF and eST) exhibit} second-order convergence, 
whereas SC has fallen below first-order.
\textcolor{red}{Finally, the comparison between the SF technique of~\cite{Paciorri1} and the eST technique described here
reveals that the latter incurs in a slightly larger discretization error than the former within the shock-downstream region. 
This observation points to the fact that there is room for improving the data transfer algorithms described in Sect.~\ref{first transfer}
and~\ref{second transf}, which will be the subject of future work}.

	%
\subsection{\color{red}Cost vs.\ accuracy analysis}
{\color{red}
A comparative assessment of the computational cost of the shock-fitting and shock-capturing approaches
can be made by either: {\it i}) using the same meshes, or {\it ii}) estimating the (different)
mesh spacing required by the two techniques to achieve the same discretization-error level
in the shock-downstream region.

If the first standpoint is adopted, i.e.\ the same grids are used, it is clear that shock-fitting techniques,
thus including both SF and eST, incur an higher computational cost per time-step than SC.
%
%
This is because, in addition to solving the governing PDEs on the same mesh, using the same CFD solver
also used in the SC simulation,
shock-fitting techniques also have to keep track of the shock motion by solving the Rankine-Hugoniot relations
at all shock-points.
However, since the shock-mesh has a lower dimensionality ($d-1$ in the $d$-dimensional space)
than the mesh that fills the computational domain,
the overall increase in computational cost incurred by either the SF or eST techniques
amounts to a relatively small fraction of the cost per iteration of the SC technique.
The interested reader can find
a detailed analysis on the computational costs incurred by the SF technique in~\cite{Grottadaurea2011}.
%


If, on the contrary, the second standpoint is adopted, shock-fitting methods are seen to outperform
shock-capturing when it comes to achieve the same discretization error level in the shock-downstream region.
The idea here is to use  the results of the grid convergence tests of Sect.~\ref{Q1D} and~\ref{planarflow} 
to  estimate the mesh sizes required by the SC solver 
to provide error levels comparable to those of the eST approach. 
The analysis described below is based on the discretization error 
of the third component $\epsilon\left(\sqrt{\rho}u\right)$ of $\mathbf{Z}$ for the Q1D nozzle flow,
whereas the fourth component $\epsilon\left(\sqrt{\rho}v\right)$ has been considered on two-dimensional grids. 
%
%
%
%

The computations of the Q1D nozzle-flow of Sect.~\ref{Q1D} show that using the SC solver,
which has a shock-downstream convergence rate of about $1.1$, see Fig.~\ref{pic19},
a mesh size $h \simeq 3.75\; 10^{-5}$ would be required
to obtain the same error provided by the eST approach on the coarsest, level 0 mesh, see Tab.~\ref{tab1}.
This amounts to say that in 1D SC requires a mesh that is about 32 ($=h_0/h$) times finer than the level 0 mesh
to obtain the coarse-grid eST result.
Note that $h$ is even smaller than the mesh spacing $h_4 = 2.344\;10^{-4}$ (see Tab.~\ref{tab1})
of the finest mesh used in the grid-convergence study and that such a fine mesh
would be required just to compensate the error generated by capturing the shock.
Moreover, in order to obtain the same discretization error of the eST approach on the finest (level 4) mesh, 
the SC solver would need a mesh size $h' \simeq 2.73\;10^{-7}$, which is three orders of magnitude smaller than $h_4$.

The two-dimensional source flow computations of Sect.~\ref{planarflow} show that in order
to obtain the same discretization error of eST on the coarsest, level 0 mesh, SC would require a mesh spacing $h \simeq h_0/64$,
whereas to attain an error level comparable to that obtained by eST on the finest, level 2 mesh, 
SC would require a mesh spacing $h' \simeq 1.53\; 10^{-6}$, which is roughly four orders of magnitude smaller than $h_2$, see Tab.~\ref{tab2}. 

These results give indications that using uniform refinement in 2D,
the same error level of the coarse-grid, eST result would be attained with SC using 
a mesh having a number of triangular elements that is $\left(h_0/h\right)^2 \simeq 4096$ times larger
than the number of triangles of the level 0 mesh, see Tab.~\ref{tab2}. This amounts to
a number of triangular cells of the order of $10^7$.
Following the same line of reasoning,
obtaining the fine-grid, eST result 
using SC would instead require a 2D mesh with a number of elements
that is $\left(h_2/h'\right)^2 \simeq 10^8$ times larger
than the level 2 grid, which is clearly impractical. A possible solution would be to replace eST by some 
aggressive error estimation and anisotropic metric-based adaptation techniques, as e.g.\ those 
proposed in~\cite{ALAUZET201828}.
However, one should evaluate the capabilities of these techniques to  provide meshes with the 
anisotropy ratios required to drop the mesh size down several orders of magnitudes in the shocks,  
the capabilities of the flow solver involved to handle such meshes,
and, finally, the overhead of the mesh refinement itself compared to that of the eST method. 
This is perhaps  a possible avenue for future work.
}

\subsection{Blunt body problem} 
\label{subsec:blunt}

\begin{figure}
\begin{center}
\includegraphics[width=0.55\textwidth,trim={3.2cm 2.5cm 3cm 2cm},clip]{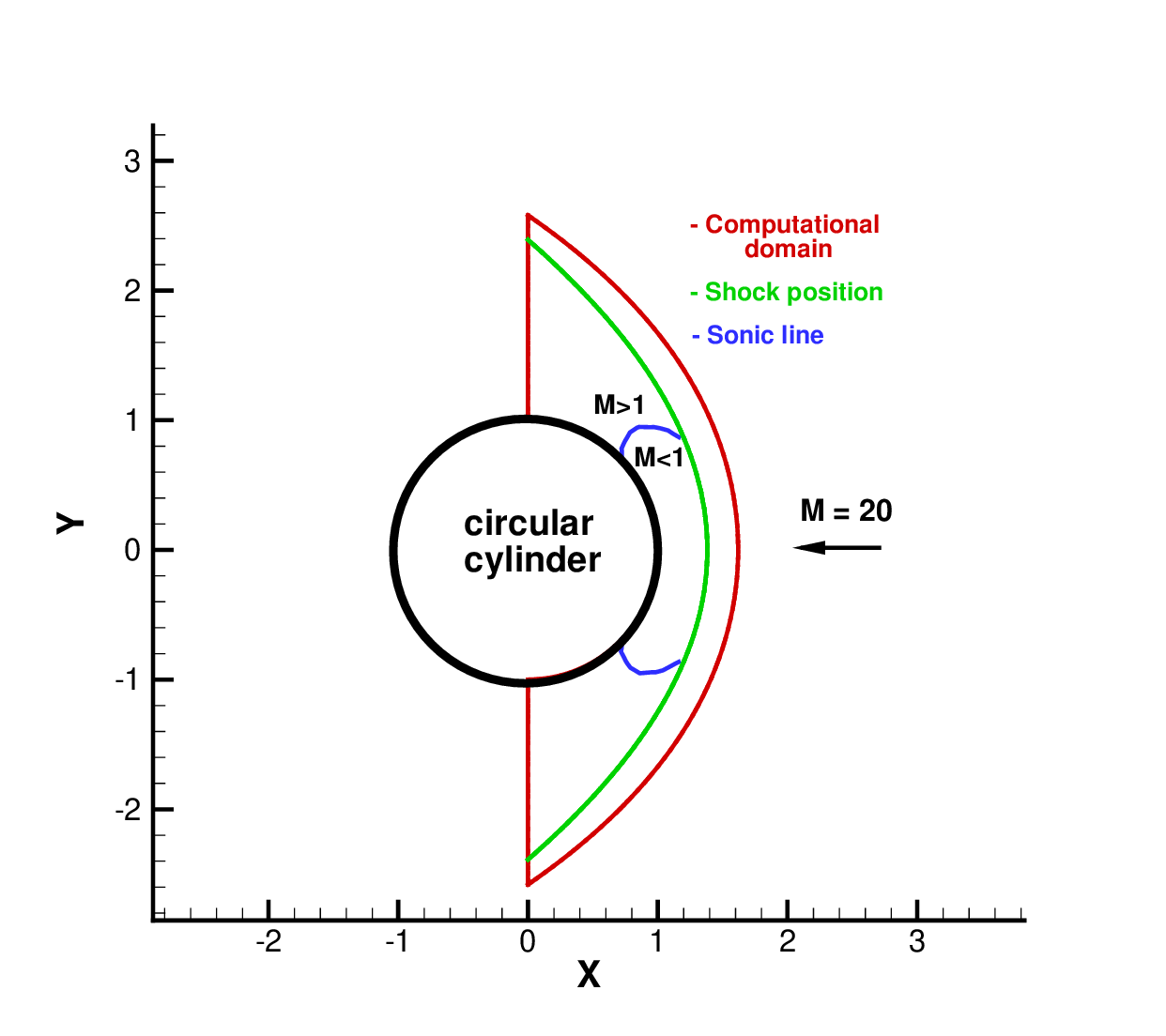}
\end{center}
\caption{Hypersonic flow past a circular cylinder: sketch of the computational domain.}\label{pic23}
\end{figure}

The hypersonic ($M_\infty$ = 20) flow past the fore-body of a circular cylinder, 
see Fig.~\ref{pic23}, is a comprehensive test-bed for the eST algorithm, because the entire shock-polar is swept whilst moving along the bow shock which stands ahead of the blunt body. 
The existence of the subsonic pocket that surrounds the stagnation point 
and the transition to supersonic flow through the sonic line may be challenging for the proposed method and, in particular, for the algorithms used to transfer data back and forth between the shock and surrogate boundaries.

\begin{figure}
\begin{center}
\includegraphics[width=0.6\textwidth]{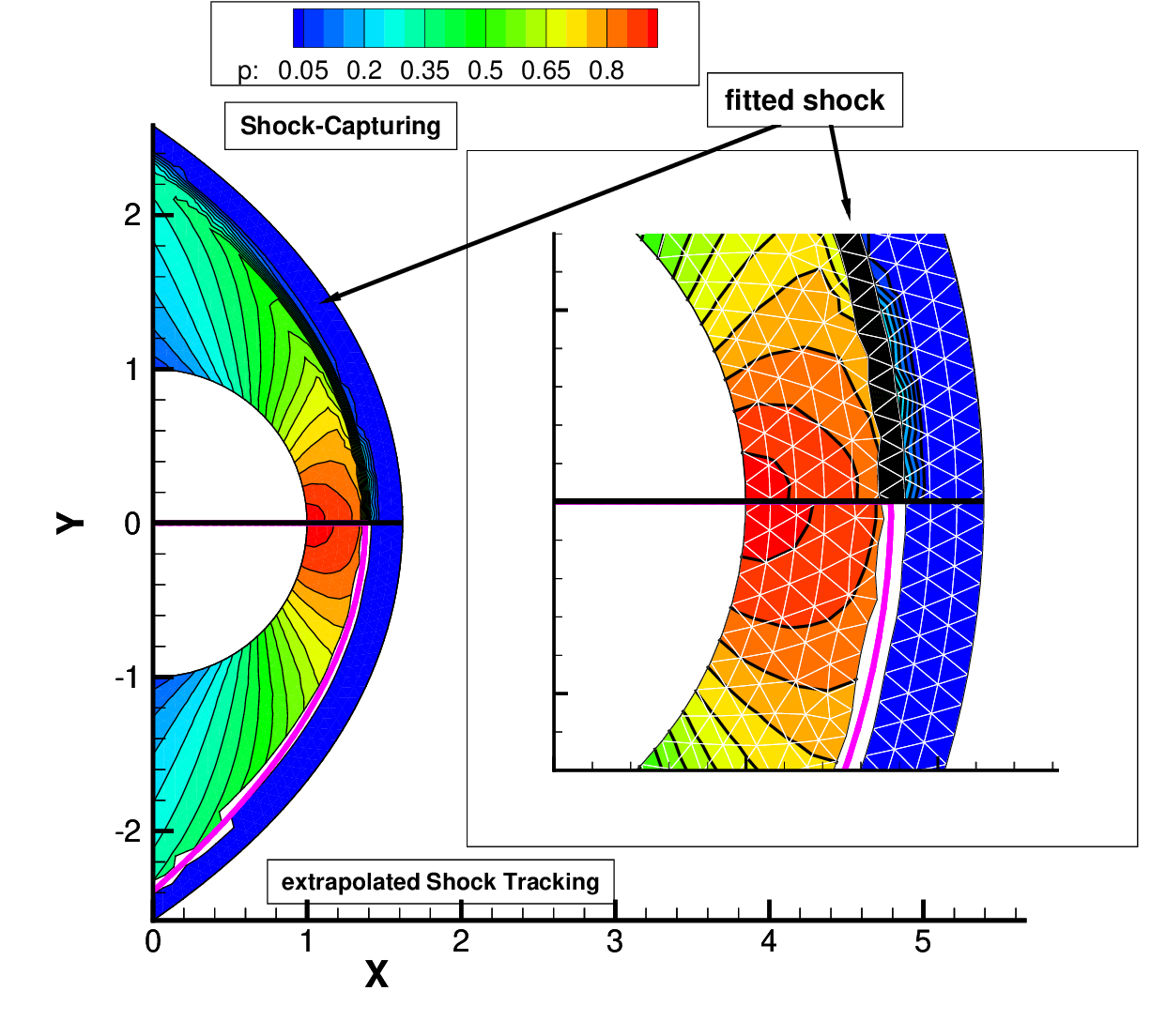}
\end{center}
\caption{Hypersonic flow past a circular cylinder: comparison between the pressure iso-contours computed by means of SC (top) and eST (bottom). Computed shock curve in pink.}\label{pic24}
\end{figure}

The mesh used as the background triangulation in the eST simulation has also been used to run the SC simulation. The computations have been run on a Delaunay mesh containing 808 points and 1458 triangles, generated using \texttt{\ttfamily delaundo}~\cite{Muller1,Muller2}.
A close up view of the mesh is reported on Fig.~\ref{pic24}: it can be seen that the only difference between the eST and SC computation consists in the removal of the triangles crossed by the fitted-shock in the eST case.

Pressure iso-contour lines are shown in Fig.~\ref{pic24}: the SC calculation is shown in the upper half of both frames
and the eST one in the lower half. The steady location of the fitted bow shock (shown using a solid bold line) has also been superimposed on the eST results. The comparison clearly reveals that the differences between the solutions obtained using the two different shock-modeling practices are remarkable within the entire shock-layer.

\begin{figure}
\begin{center}
\includegraphics[width=0.6\textwidth]{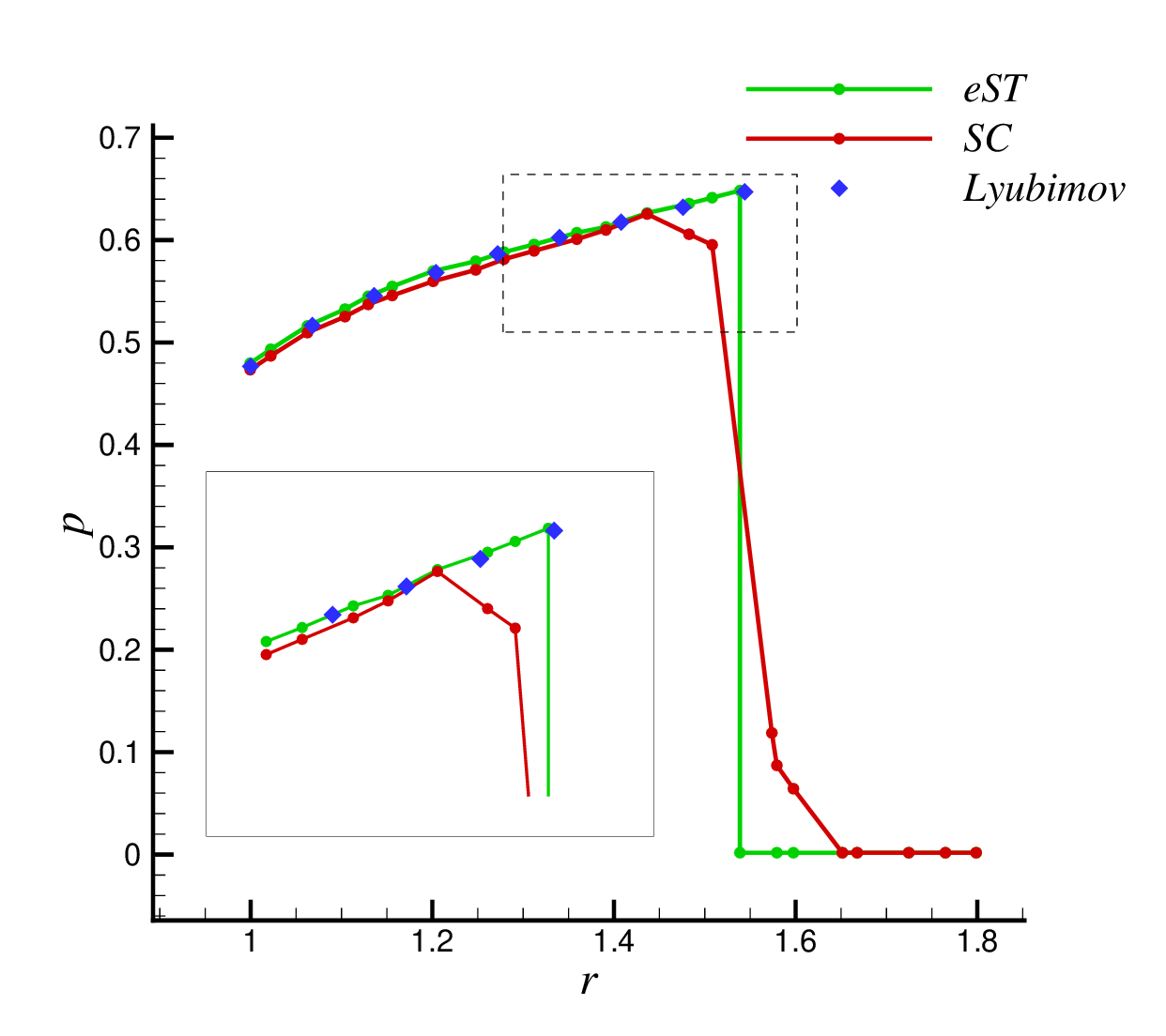}
\end{center}
	\caption{Hypersonic flow past a circular cylinder: SC, eST and reference~\cite{Rusanov} pressure distribution within the shock-layer.}\label{pic25}
\end{figure}

Figure~\ref{pic25} shows the pressure $p$ profile probed along a line
that makes a 45$^{\circ}$ angle w.r.t.\ the centerline. The SC and eST results have been compared with the
reference solution computed in~\cite{Rusanov}.
The comparison shows that the finite shock-width of the SC solution is replaced by a discontinuity
in the eST result (which also includes the shock-upstream and shock-downstream values) and that the 
shock stand-off distance computed by eST agrees very well with the reference solution.

\begin{figure*}[t]
\centering
	\subfloat[SF]{%
		\includegraphics[width=0.40\textwidth,trim={0 0 9cm 0},clip]{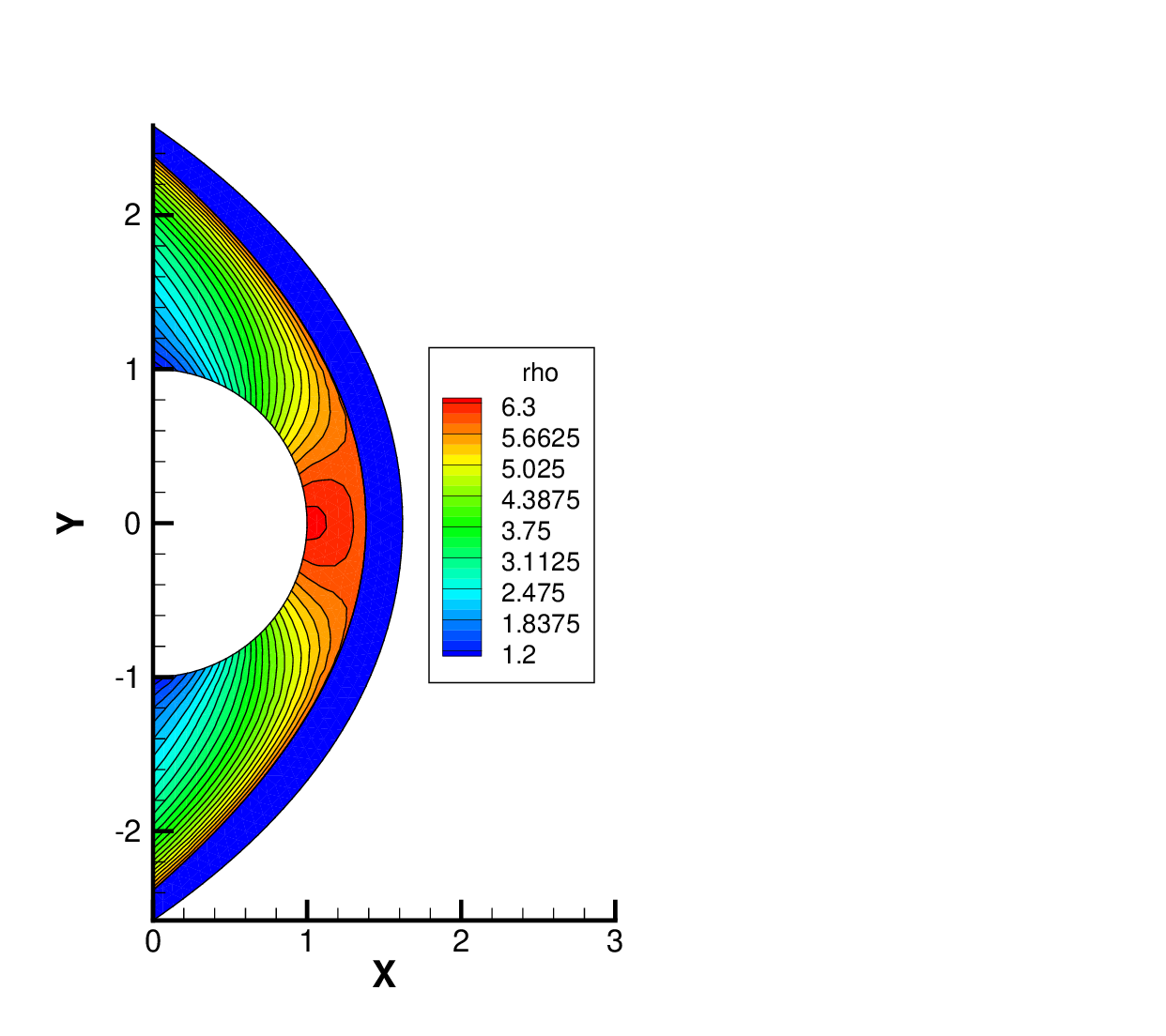}\label{pic50a}}
\hspace{\fill}
	\subfloat[eST]{%
		\includegraphics[width=0.40\textwidth,trim={0 0 9cm 0},clip]{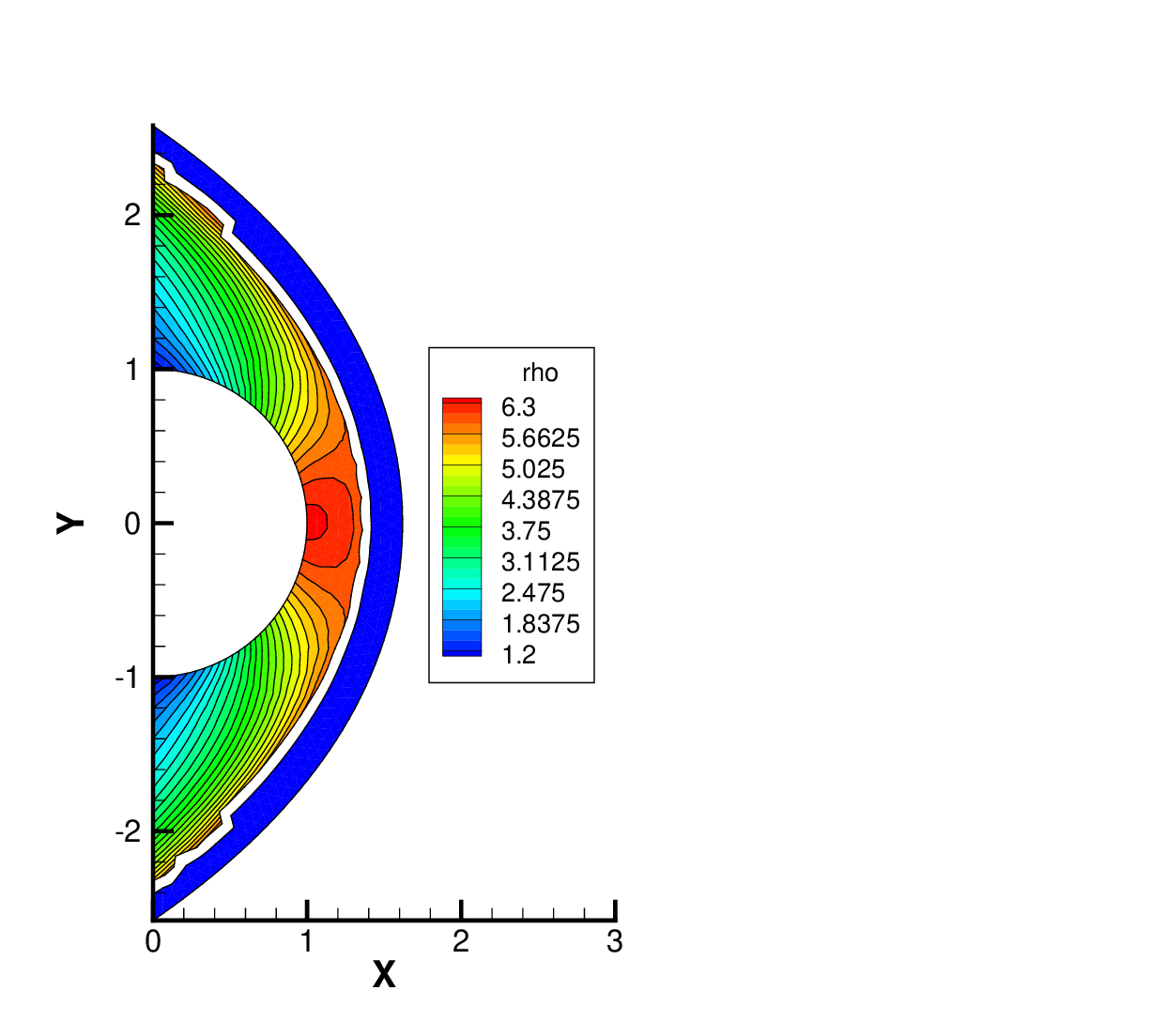}\label{pic50b}}
	\caption{\textcolor{red}{Hypersonic flow past a circular cylinder: comparison between SF and eST in terms of density iso-contour lines.}}\label{pic50}
\end{figure*}

\textcolor{red}{Finally, Fig.~\ref{pic50}, which compares the density iso-contour lines computed by 
the SF technique of Ref.~\cite{Paciorri1} and the eST technique described here,
turns out to be very useful to pinpoint the methodological
differences between the two different shock-fitting approaches.
Observe, in particular, that in the eST simulation the solution has not been computed within the blank
region surrounding the fitted shock (for clarity, the fitted-shock has not been drawn in Fig.~\ref{pic50b}). 
Even so, the eST solution within the shock-layer is
as smooth as it is the one computed by SF}. 

\subsection{Hybrid  computations of interactions}
\label{subsec:hybrid}
In its current implementation, the eST method cannot explicitly track shock interactions. However, in this section we will show that it can be applied without any problem to this type of flows by means of a hybrid fit-capture approach.
We will in particular consider two applications: a steady Mach reflection in a channel with a ramp, and a type IV shock-shock interaction arising in a supersonic flow around a circular cylinder. \\

{\it Steady Mach Reflection.} This case is quite useful as it involves a relatively simple flow pattern, but allows to clearly visualize the advantage brought by the eST approach w.r.t. SC.

\begin{figure}[h!]
\begin{center}
\includegraphics[width=0.8\textwidth]{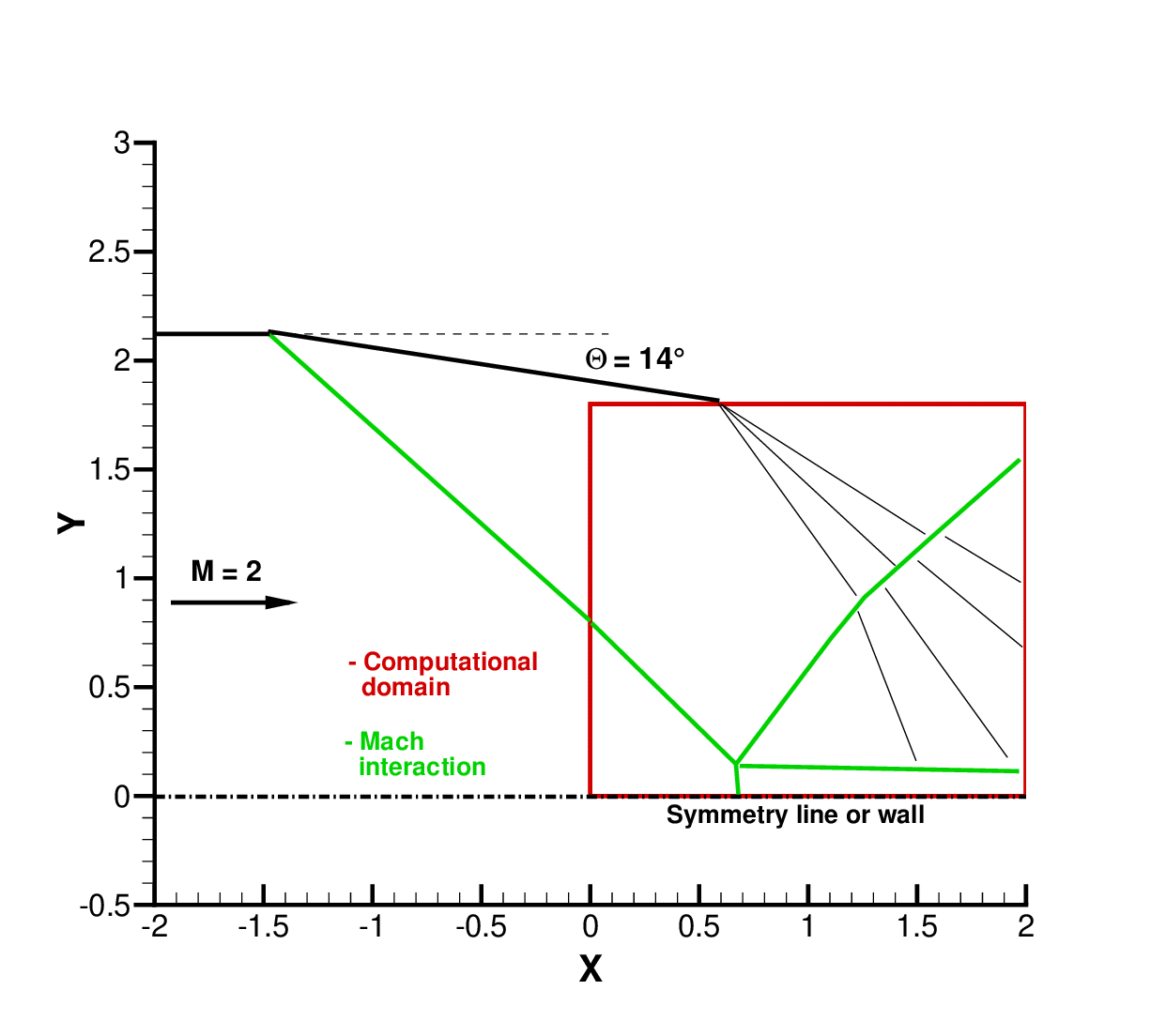}
\end{center}
\caption{Steady Mach reflection: sketch of the computational domain.}\label{pic26}
\end{figure}

The setup of the test case  is the same as in~\cite{Paciorri1} and sketched on Fig.~\ref{pic26}: it involves a $M_{\infty}=2$ flow in a channel with a wall deflection of 14 degrees. The oblique shock forming due to this deflection reflects onto the channel walls. In these conditions the reflection is not a regular one, but a Mach reflection is observed with its typical lambda-shock topology. A sketch of the resulting interaction is reported in Fig.~\ref{pic26}. Note that, as a result of this interaction, a contact discontinuity emanates from the triple point.

\begin{figure*}[h!]
   \subfloat{%
     \includegraphics[width=0.333\textwidth,trim={0 0 10cm 0},clip]{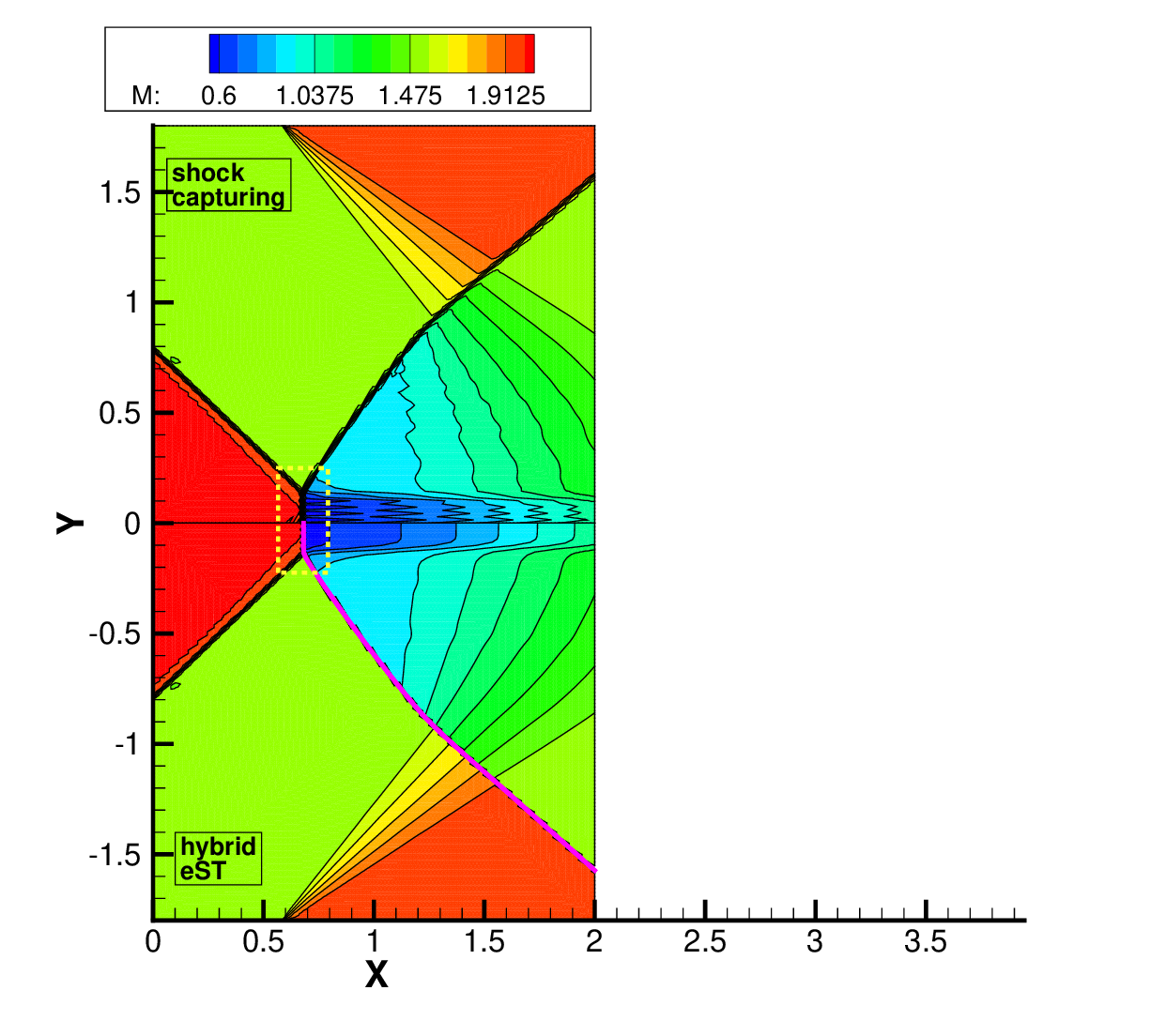}}
   \subfloat{%
	\includegraphics[width=0.333\textwidth,trim={0 0 10cm 0},clip]{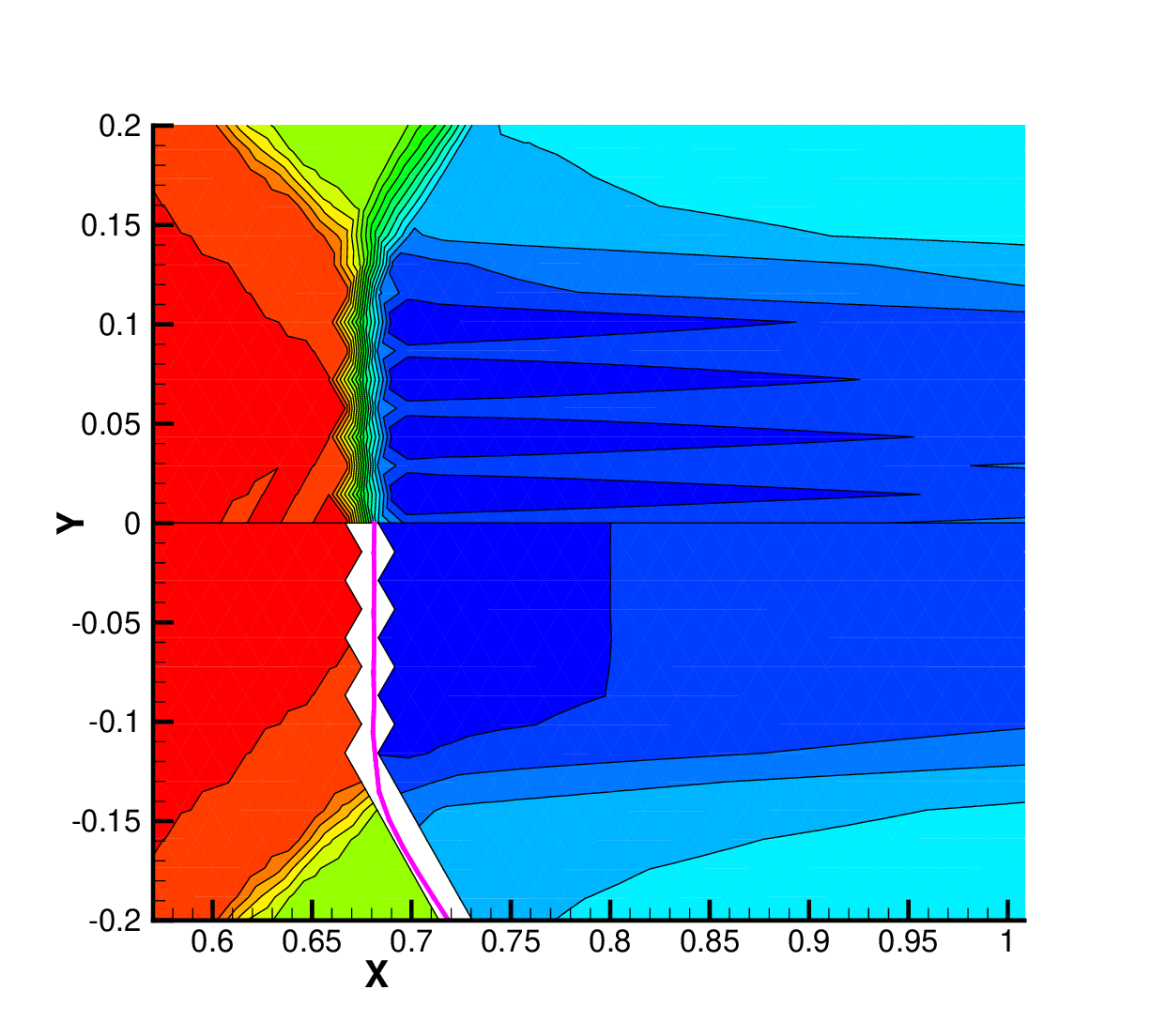}}
   \subfloat{%
	\includegraphics[width=0.333\textwidth,trim={0 0 10cm 0},clip]{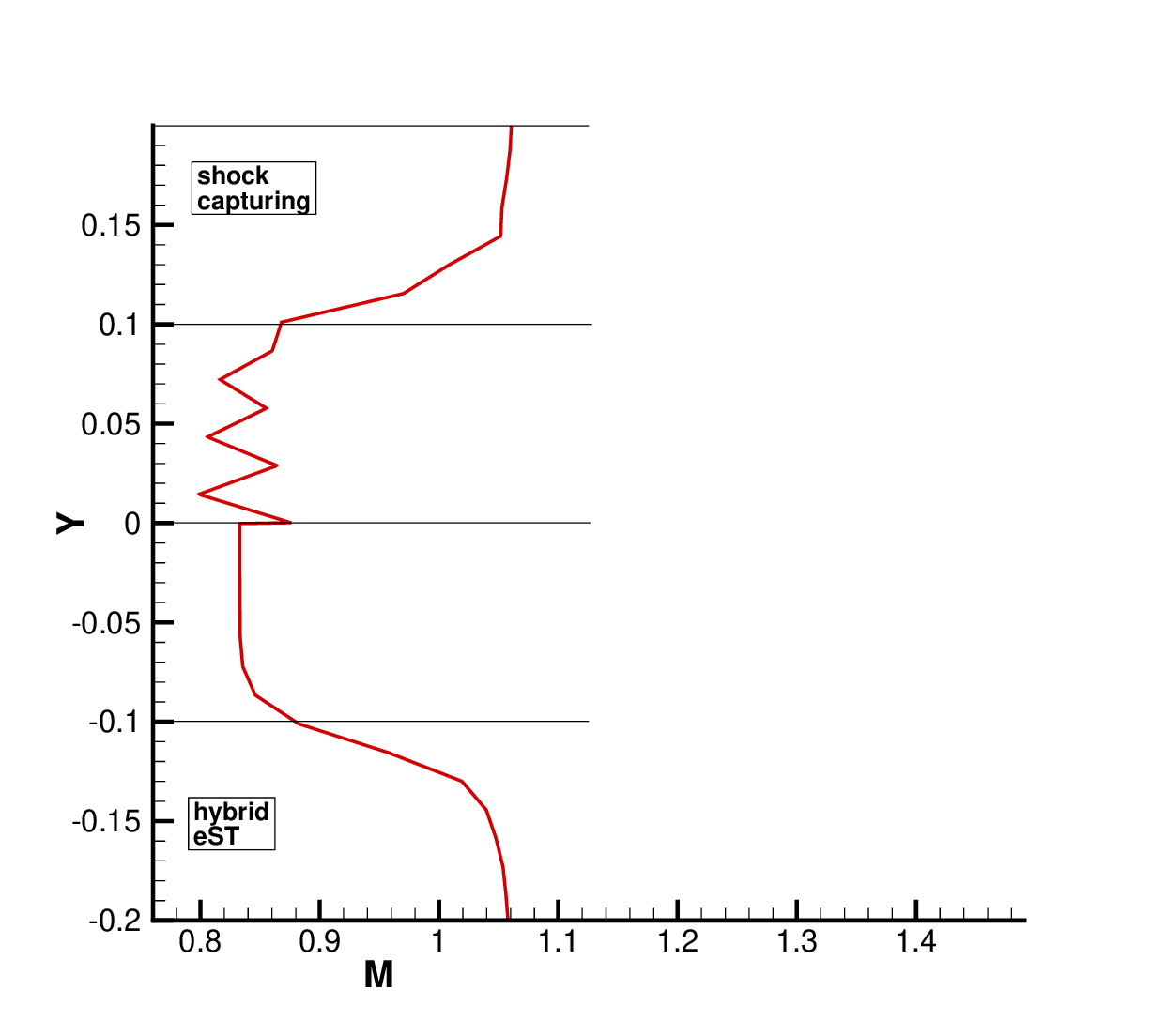}}
\caption{Steady Mach reflection: Mach number iso-contours comparison, enlargement around the triple point and Mach number distribution at \emph{x = 1.5}.}\label{pic27}
\end{figure*}

The background mesh used for this simulation contains 14833 grid-points and 29214 triangles. 
We compare on this mesh the SC solution with the hybrid result in which eST is only applied to two of the branches of the lambda shock:
the Mach stem and the reflected shock. Both the incident shock and the contact 
discontinuity are captured. Figure~\ref{pic27} displays the Mach iso-contours in the entire computational domain 
(the left frame), an enlargement of the region surrounding the triple point (the middle frame) and the Mach number distribution 
along a vertical line at \emph{x = 1.5} (the right frame) in the region downstream of the triple point. 
It can be seen in Fig.~\ref{pic27}, that the capture of the Mach stem gives rise to an unphysical behavior of the 
Mach number contour lines in the region downstream of the Mach stem. This unphysical behavior disappears in the 
hybrid solution that exhibits a smoother Mach number distribution. 
{\color{red}It must be noticed that the gradient reconstruction technique, described in~\ref{grad}, 
does not provide accurate gradient reconstruction for discontinuous solutions. Nonetheless, 
this is unlikely to significantly affect the overall quality of the eST computations because only very few grid-points are involved in the 
part of the domain where the interaction occurs.}

\begin{figure*}[tb!]
\centering
   \subfloat[SF]{%
     	\includegraphics[width=0.48\textwidth]{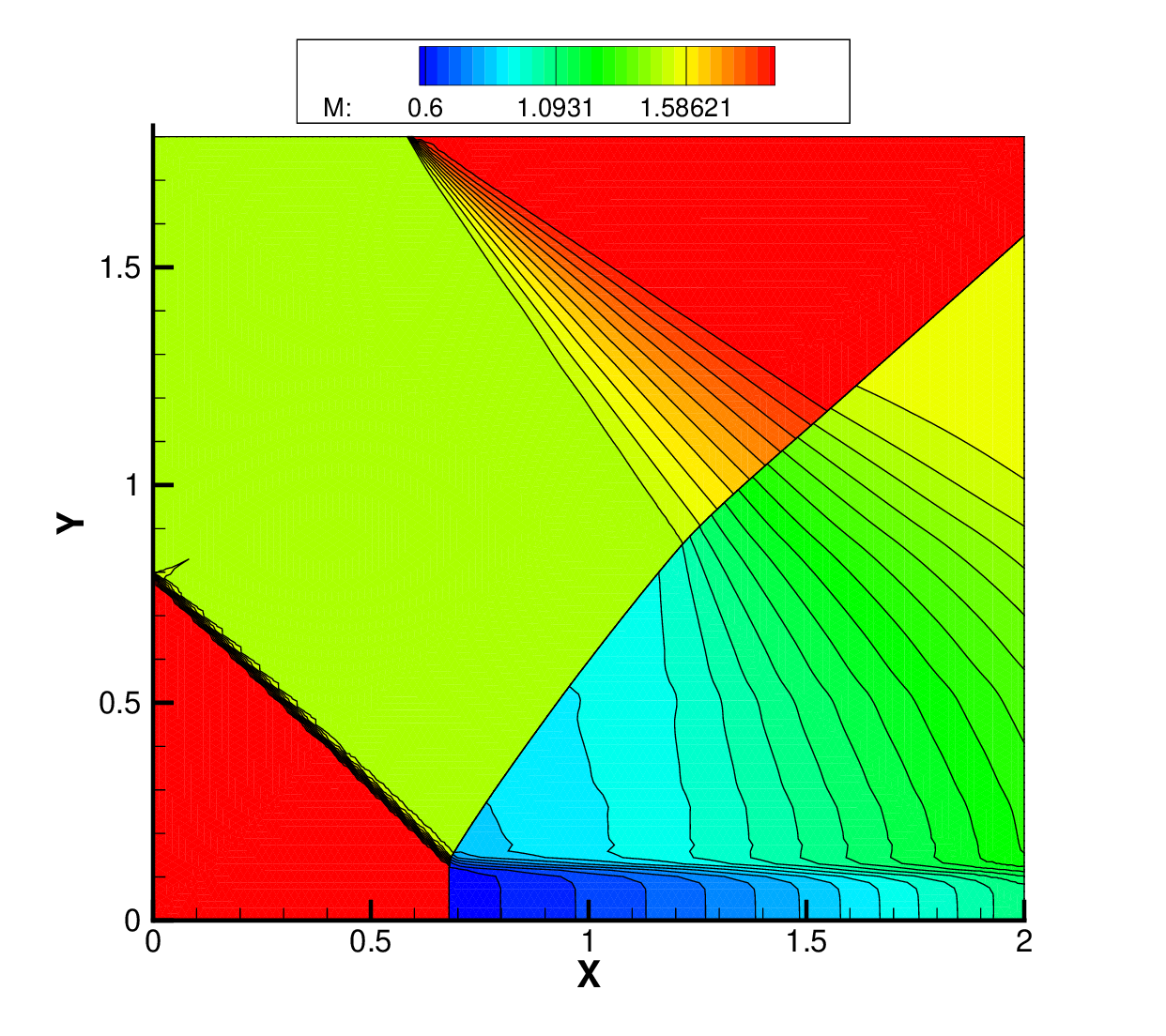}}
\hspace{\fill}
   \subfloat[eST]{%
      	\includegraphics[width=0.48\textwidth]{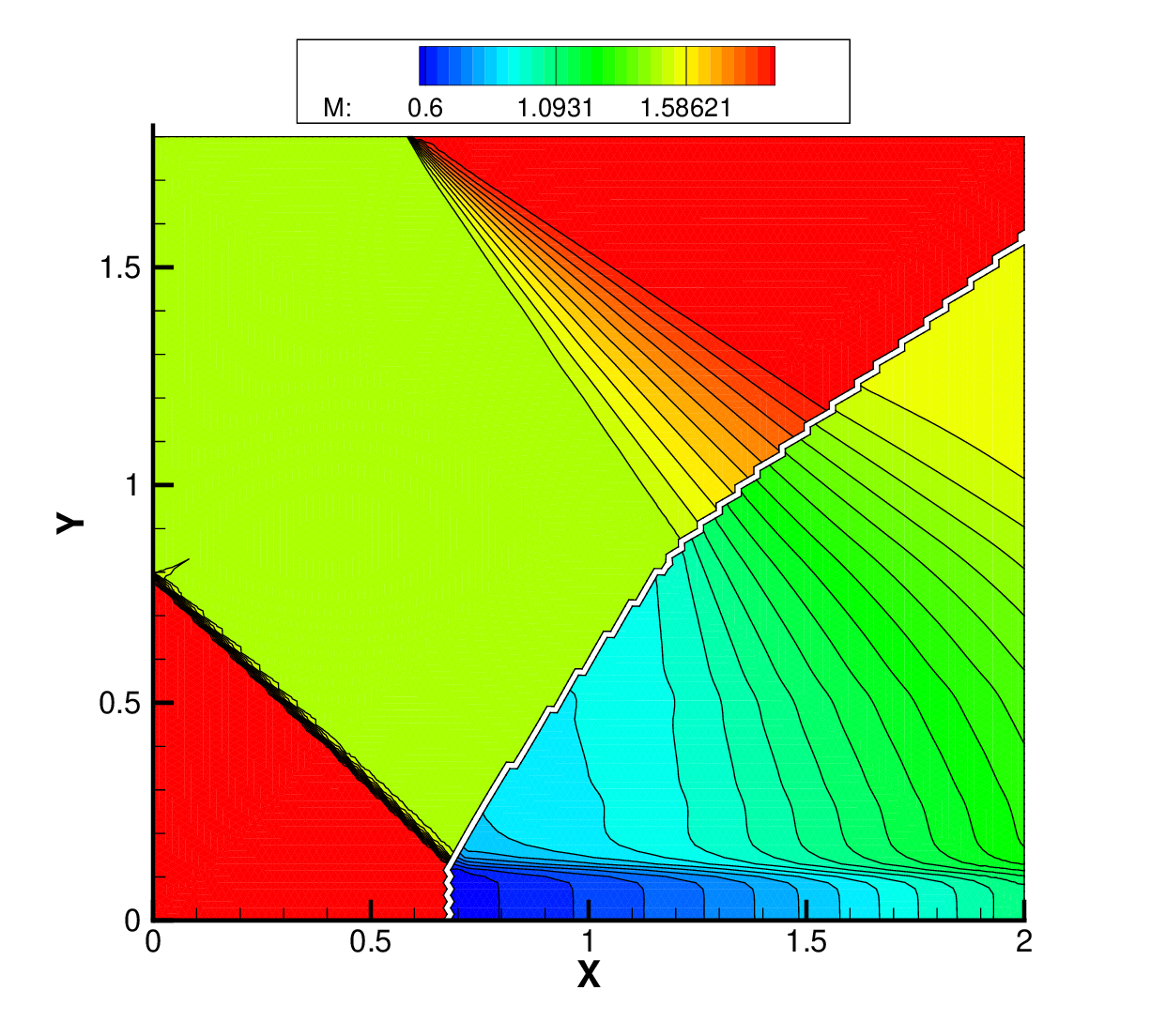}}
	\caption{\textcolor{red}{Steady Mach reflection: comparison between SF and eST in terms of Mach number iso-contour lines.}}\label{pic51}
\end{figure*}

\textcolor{red}{Finally, Fig.~\ref{pic51}, stands out, again, that the solution obtained with
the eST algorithm is notable and comparable with the one described in~\cite{Paciorri1}.}

\begin{figure}[tb!]
\begin{center}
\includegraphics[width=0.65\textwidth]{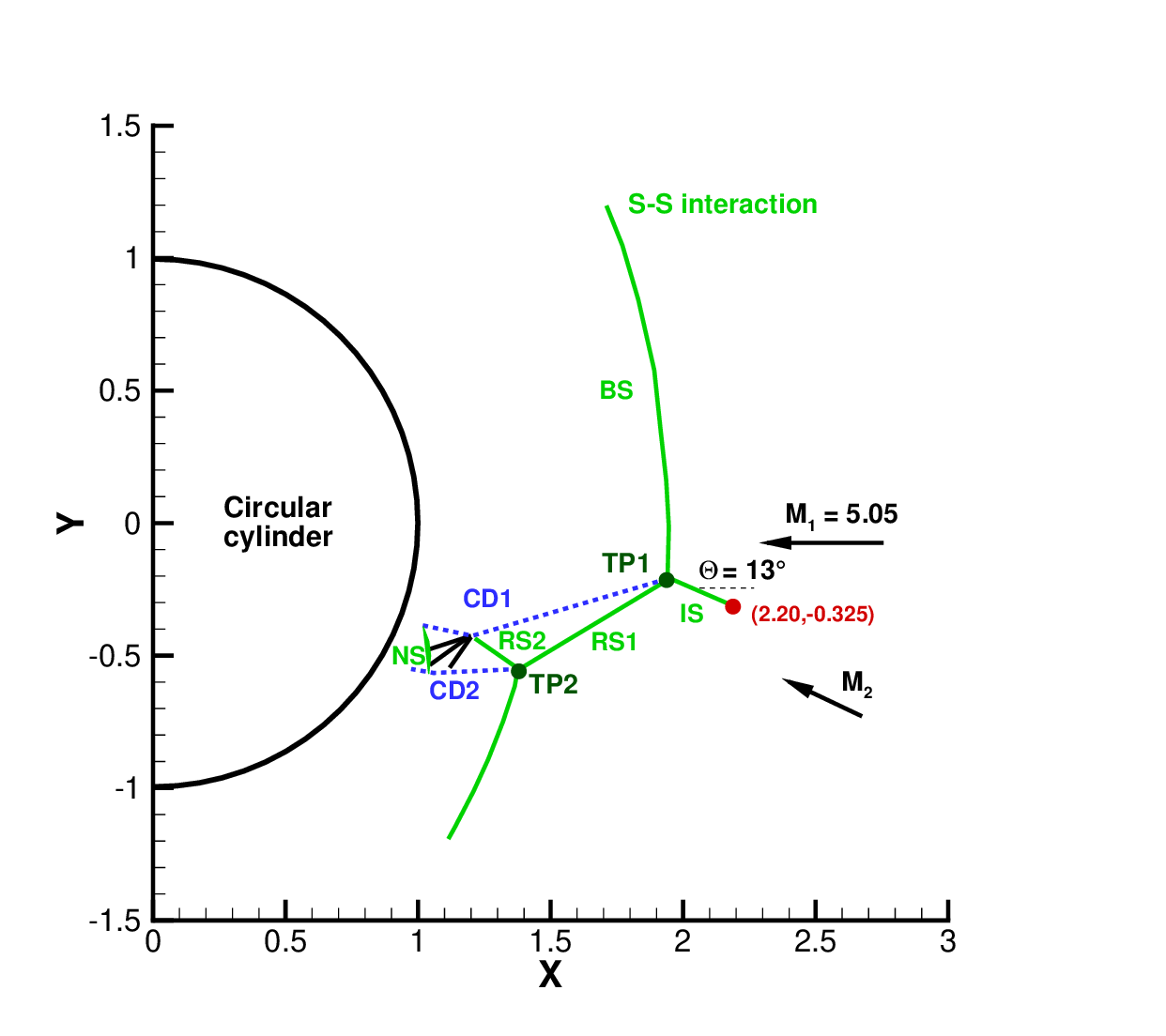}
\end{center}
\caption{Sketch of the type IV shock-shock interaction problem.}\label{pic30}
\end{figure}

{\it Type IV shock-shock interaction.} This last benchmark introduced in~\cite{SSIVpaper} involves a more complex 
pattern of  interacting discontinuities. 
An horizontal flow, characterized by a Mach number of $M=5.05$, is deflected by an oblique shock  (whose angle w.r.t.\ the horizontal direction is $\Theta=13$ degrees ) in front of a circular cylinder. The resulting flow features a bow shock, interacting with the oblique shock, and giving rise to the well known type IV interaction, which has been already studied in~\cite{Chang2019}. For clarity, we have drawn in Fig.~\ref{pic30} a sketch of this interaction.
A first triple point (TP1) occurs where the oblique shock (IS) impinges on the bow shock giving rise 
to a reflected shock (RS1) and a contact discontinuity (CD1) that move towards the stagnation point. 
The contact discontinuity separates the supersonic stream which has been deflected by the oblique shock from 
the subsonic stream downstream of the bow shock (BS). The reflected shock coming from the first triple point re-joins the 
bow shock in a second triple point (TP2) where a new reflected shock (RS2) and contact discontinuity (CD2) arise. 
The two contact discontinuities bound a supersonic jet which is directed toward the body surface. 
Within the jet the second reflected shock interacts with the first contact discontinuity giving rise to an expansion wave. 
The flow concludes his path by being decelerated by a normal shock (NS) right in front of the body surface 
causing a higher density and pressure zone on the cylinder surface. 


\begin{figure}[h!]
\centering
   \subfloat[SC]{%
	\includegraphics[width=0.5\textwidth,trim={0 0 7cm 0},clip]{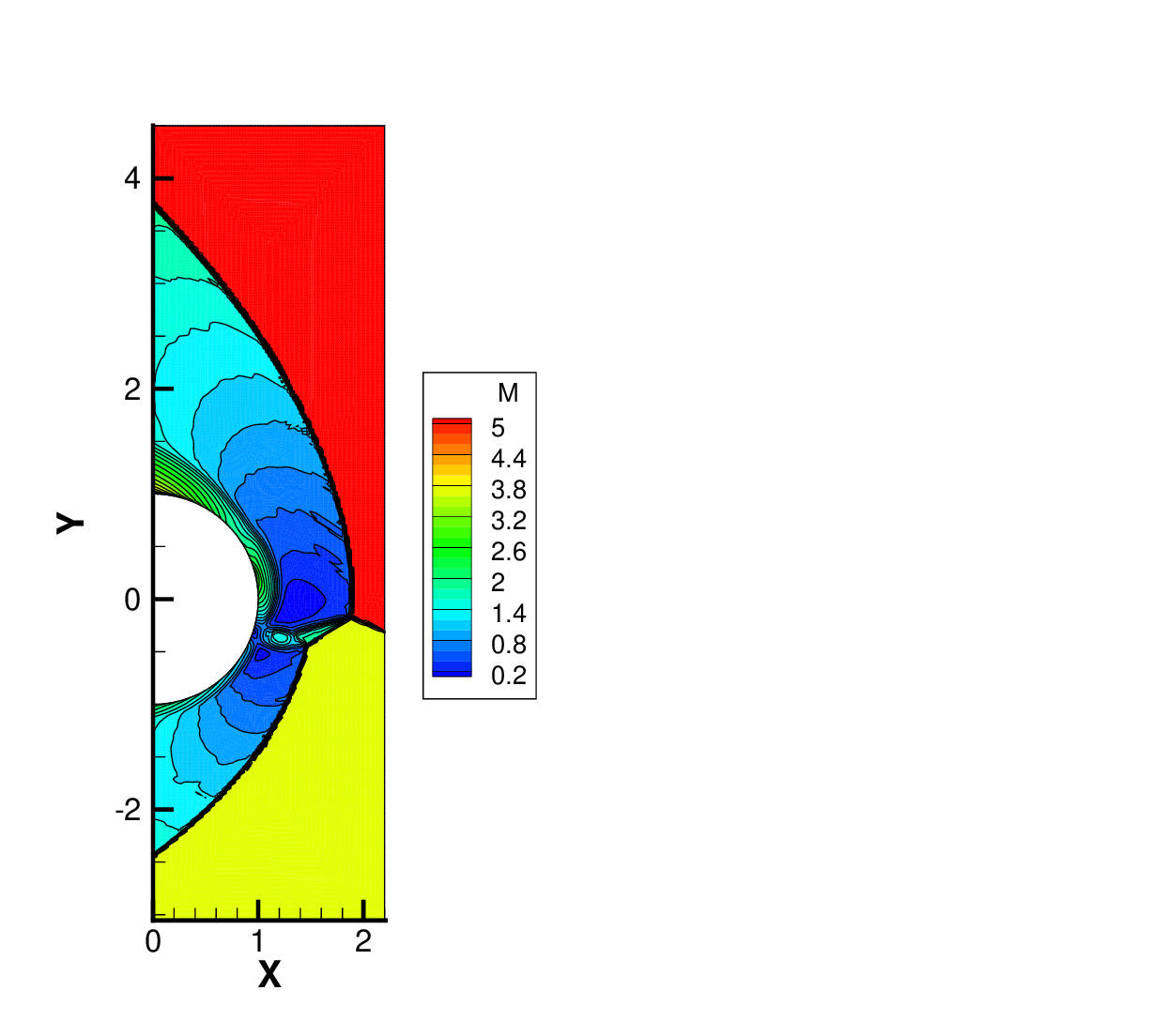}}
   \subfloat[eST]{
     	\includegraphics[width=0.5\textwidth,trim={0 0 7cm 0},clip]{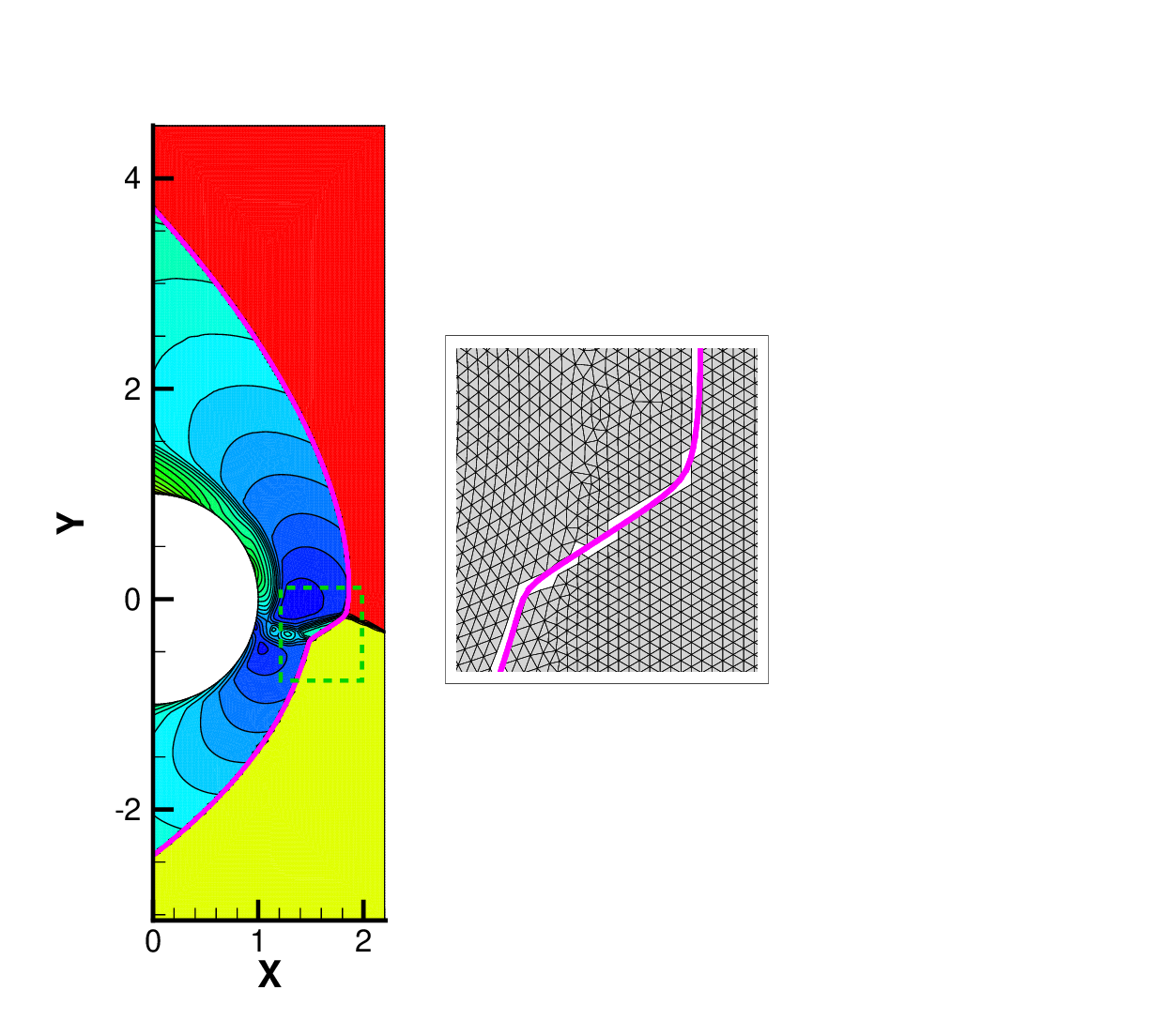}}
\caption{Type IV shock-shock interaction: Mach number iso-contours {\color{green}and a close-up on the shock-mesh}.}\label{pic32}
\end{figure}


These kinds of interactions are very difficult to study  because they require a very fine triangulation 
in order to properly describe what is going on within the flow-field. 
A Delaunay mesh containing 49660 triangles and 25231 nodes was used to compute the SC solution and also as background mesh for the hybrid computation in which eST was used to fit the entire bow-shock and the oblique shock RS1. 
Figures~\ref{pic32}.a and \ref{pic32}.b show the  
differences between the two solutions in terms of density and Mach number iso-contours. 
The pink bold line that appears in Fig.~\ref{pic32}.b represents the fitted shock, better displayed in the close-up. 
As before, the use of eST allows to obtain a considerably smoother flow-field inside the shock layer.

\section{Conclusions and perspectives}

A novel technique to simulate flows with shock waves has been illustrated and tested on several applications 
on one-dimensional and two-dimensional unstructured grids. 
The proposed extrapolated  shock tracking method borrows ideas from embedded boundary methods, combining them with a floating shock-fitting approach.
The resulting technique has been proven to be able to provide genuinely second order results for flows with very strong shocks, 
without the complexity of the re-meshing phase of the previous fitting approaches. 
The method proposed has great potential in constructing generic shock-fitting/tracking strategies,
with little dependence on the   data structure of the underlying flow solver.
As all shock/front tracking methods it has the enormous advantage of solving the exact jump conditions across the discontinuity,
which makes these methods very competitive with any kind of adaptive capturing procedure, unless these conditions are embedded in the 
discretization, as done in some DG-based recent work \cite{persson,corrigan}. This however requires to set up a dedicated solver, while 
our approach has the potential to be coupled with  several different existing CFD codes.

Indeed  one of the future challenges will be to compare its performance when coupled with different CFD solvers, not only unstructured cell-vertex, but also cell-centered
and fully structured/Cartesian codes.  Space for improvement of the method is clearly present with respect to 
its capability to handle explicitly interactions, improving the accuracy of the solution transfer to/from the shock, treating moving and complex three-dimensional discontinuities, 
and going beyond second order of accuracy.

\appendix
\section{Gradient reconstruction}\label{grad}

{\color{green}The truncated Taylor series expansion~\eqref{eq5} which is used
to transfer the dependent variable $\mathbf{Z}$ between the surrogate boundaries and the shock-mesh
relies on the availability of the gradient $\nabla \mathbf{Z}_i\left(\tilde{\mathbf{x}}\right)$ in points,
such as $A^i$ in Fig.~\ref{pic8}, and $B^i$ in Fig.~\ref{pic11},
which, respectively, belong to the surrogate boundary $\tilde{\Gamma}_U$ and $\tilde{\Gamma}_D$.
As explained in Sect.~\ref{first transfer},
the calculation of the gradient in $A^i$ or $B^i$, by means of Eq.~\eqref{eq6},
requires the knowledge of the gradient in the grid-points of the surrogate boundaries.}

\begin{figure}
\begin{center}
	\subfloat[\color{green}The gradient in grid-point \emph{i} is computed by collecting fractions of the cell-wise constant gradients of all the cells that surround $i$.]{%
\includegraphics[width=0.35\textwidth]{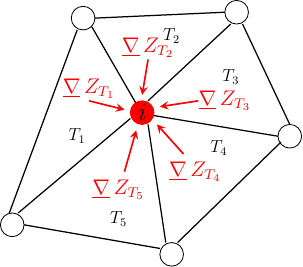}\label{pic13}}\qquad
	\subfloat[Grid elements and vectors normal to the edges of the triangle.]{%
\includegraphics[width=0.35\textwidth]{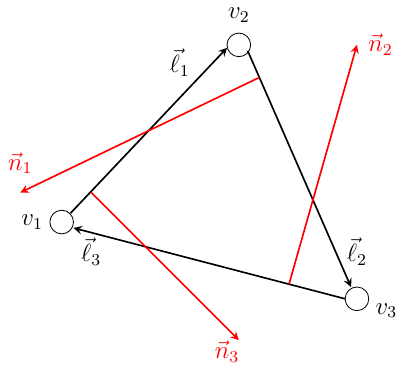}\label{pic12}}
	\caption{Reconstruction of the gradient.}
\end{center}
\end{figure}

\begin{figure}[h!]
\begin{center}
\includegraphics[width=0.5\textwidth]{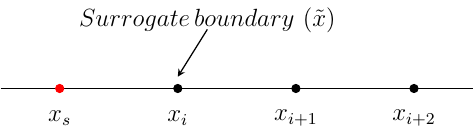}
\caption{Elements used to build the gradient on one-dimensional grids} \label{pic15}
\end{center}
\end{figure}

{\color{green}Since the dependent variable $\mathbf{Z}$ is stored in the grid-points of the
triangulation and varies linearly in space,
$\nabla \mathbf{Z}$ is not readily available within the grid-points, but it has to be reconstructed there
using the cell-wise constant gradient of the cells that surround a given grid-point,
as sketched in Fig.~\ref{pic13}.}
{\color{green}More precisely, the following} area-weighted average is used\footnote{In this Appendix the notation $\mathbf{Z}_i$ collectively refers to the four components of $\mathbf{Z}$ in grid-point $i$, rather than to the $i$th component of $\mathbf{Z}$.}:
\begin{equation}\label{eq nabla}
\nabla \mathbf{Z}_i\,=\,\frac{\sum_{i\ni T}\,(A_T\,\nabla  \mathbf{Z}_T)}{\sum_{i\ni T}\,A_T}  
\end{equation}
where the summation ranges over all the triangles that surround grid-point $i$ and $A_T$ denotes the triangle area.
{\color{green}The cell-wise constant gradient $\nabla \mathbf{Z}_T$ that appears in Eq.~(\ref{eq nabla})
can be easily computed using the values of $\mathbf{Z}$ within the vertices of triangle $T$ and
the inward normals to the edges of the triangle, scaled by the edge length $\ell$:
\begin{equation}
	\label{eq:cell_gradient}
\nabla  \mathbf{Z}_T\,=\,\frac{\sum_{i=1,3}(\mathbf{Z}_i\,\mathbf{n}_i)}{2\,\mid A_T \mid}
\end{equation}
Figure~\ref{pic12} clarifies the nomenclature used in Eq.~\eqref{eq:cell_gradient}}.
{\color{green}The gradient reconstruction described so far applies to the two-dimensional case.

In the quasi-one-dimensional framework, see Fig.~\ref{pic15} for a sketch of the 1D grid, the aforementioned approach boils down to the following 
one-sided finite-difference formula}:
\begin{equation}\label{eq11}
	\nabla \mathbf{Z}(\emph{$\tilde{x}$})\,=\, \frac{\mathbf{Z}(\emph{x}_{i+1})\,-\,\mathbf{Z}(\emph{x}_i)}{\emph{x}_{i+1} \,-\, \emph{x}_i}
\end{equation}
{\color{green}
which approximates the gradient at  the surrogate boundary, i.e.\ where $x = \tilde{x}$.}
{\color{green}
The extrapolated value of $\mathbf{Z}$ at the discontinuity is obtained from:
}
     \begin{equation}\label{eq12}
     \mathbf{Z}(\emph{x})\,=\,\mathbf{Z}(\emph{$\tilde{x}$})\,+\,\frac{\mathbf{Z}(\emph{x}_{i+1})\,-\,\mathbf{Z}(\emph{x}_i)}{\emph{x}_{i+1} \,-\, \emph{x}_i}\,(\emph{x}\,-\,\emph{$\tilde{x}$})\,+\,o(\| \emph{x}\,-\,\emph{$\tilde{x}$} \|^2)
     \end{equation}    
where \emph{$\tilde{x} = x_i$} and \emph{$x = x_s$} are, respectively, the coordinates of the surrogate boundary and of the shock-point.
	
\section{\textcolor{red}{Pseudo-temporal evolution and iterative convergence of the extrapolated Shock Tracking technique.}}
	
	\textcolor{red}{
	In this Appendix we give further insight into
	the pseudo-temporal evolution of the flow-field to show how
	the eST algorithm, starting from a converged SC solution used as initial condition, leads
	to a steady, oscillation-free, shock-fitted result}. 
	
\begin{figure*}[h!]
	\centering
		\subfloat[Iteration 1]{%
	     		\includegraphics[width=0.33\textwidth]{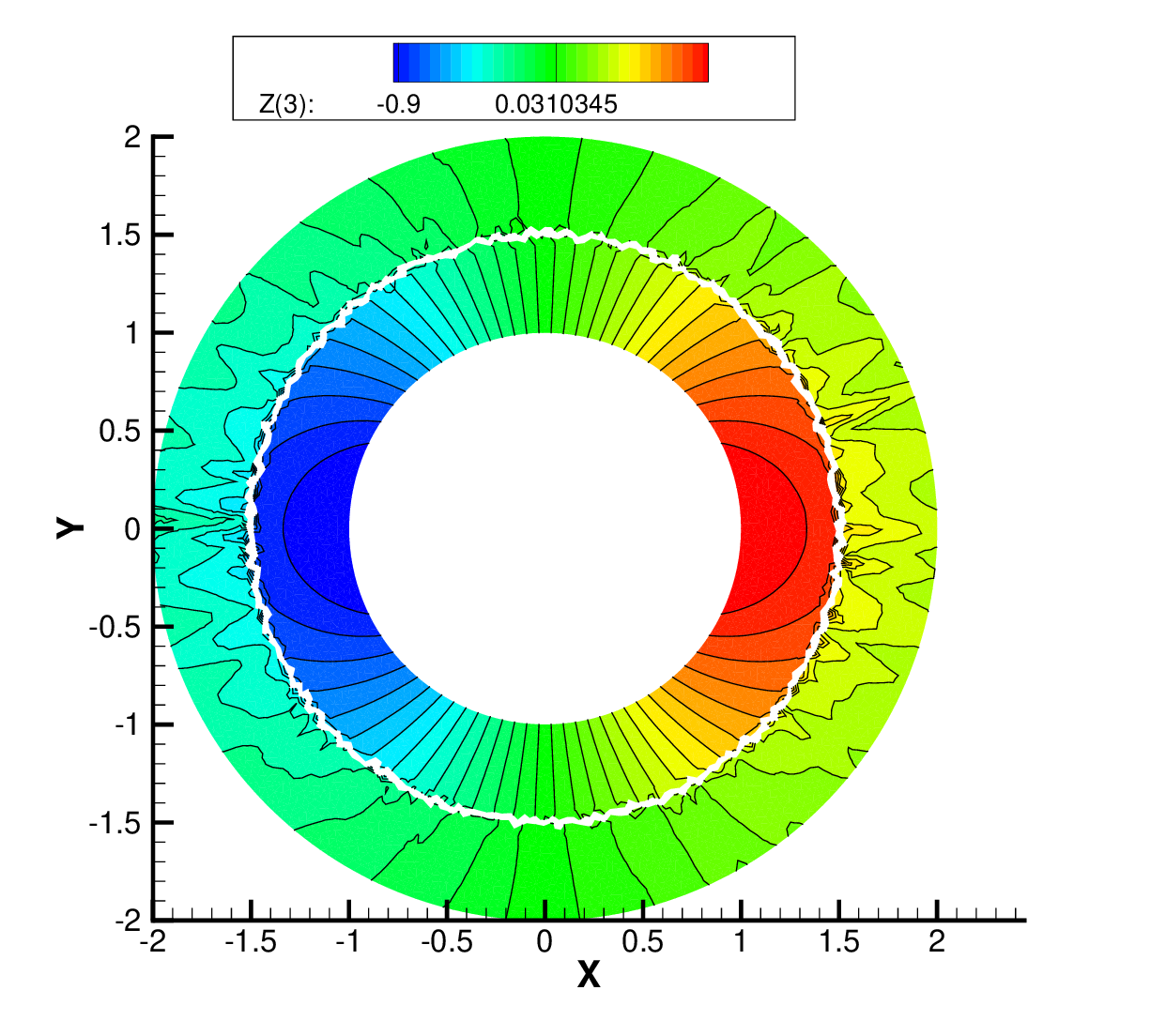}}
		\subfloat[Iteration 201]{%
	     		\includegraphics[width=0.33\textwidth]{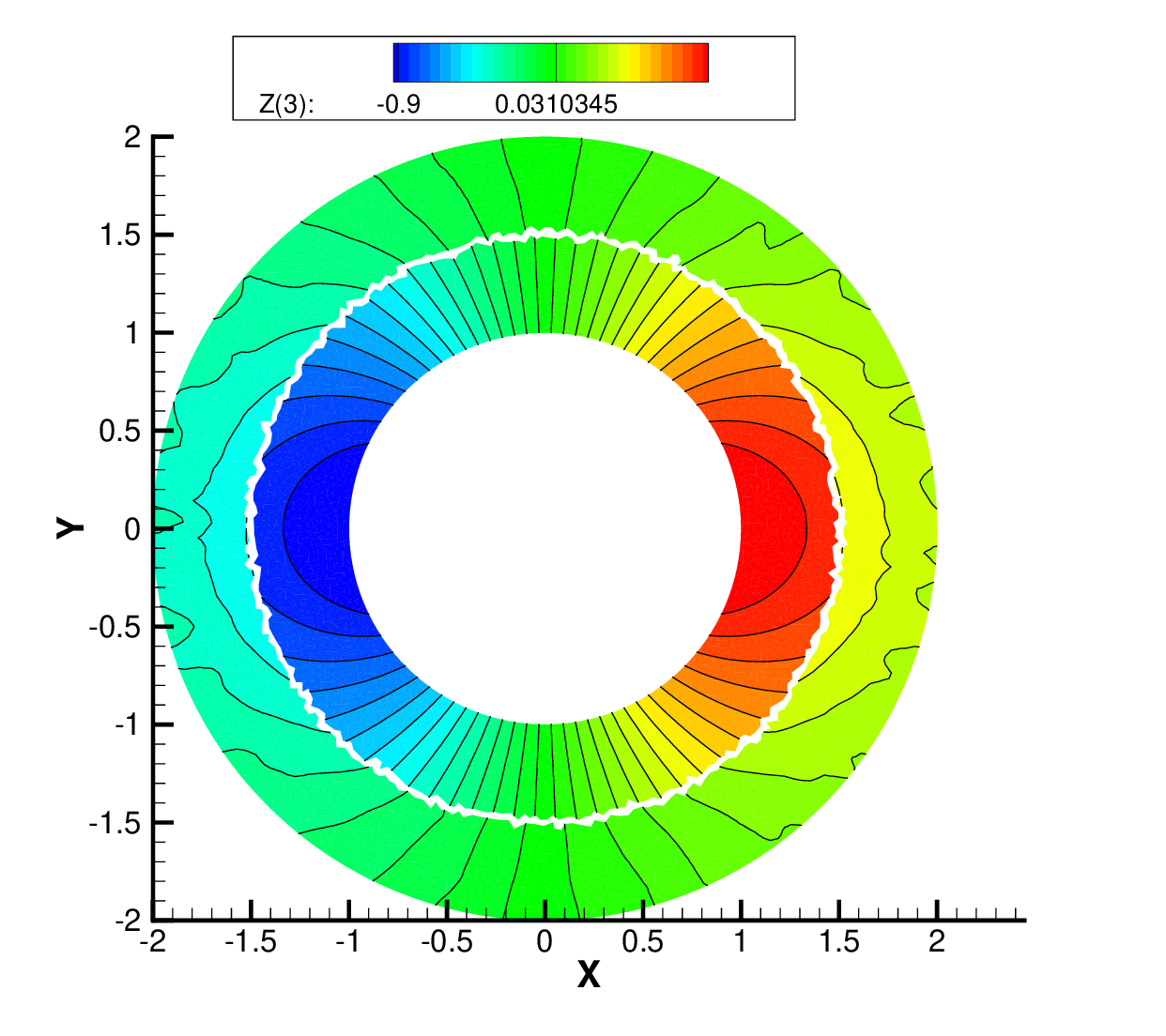}}
		\subfloat[Converged solution]{%
		\includegraphics[width=0.33\textwidth]{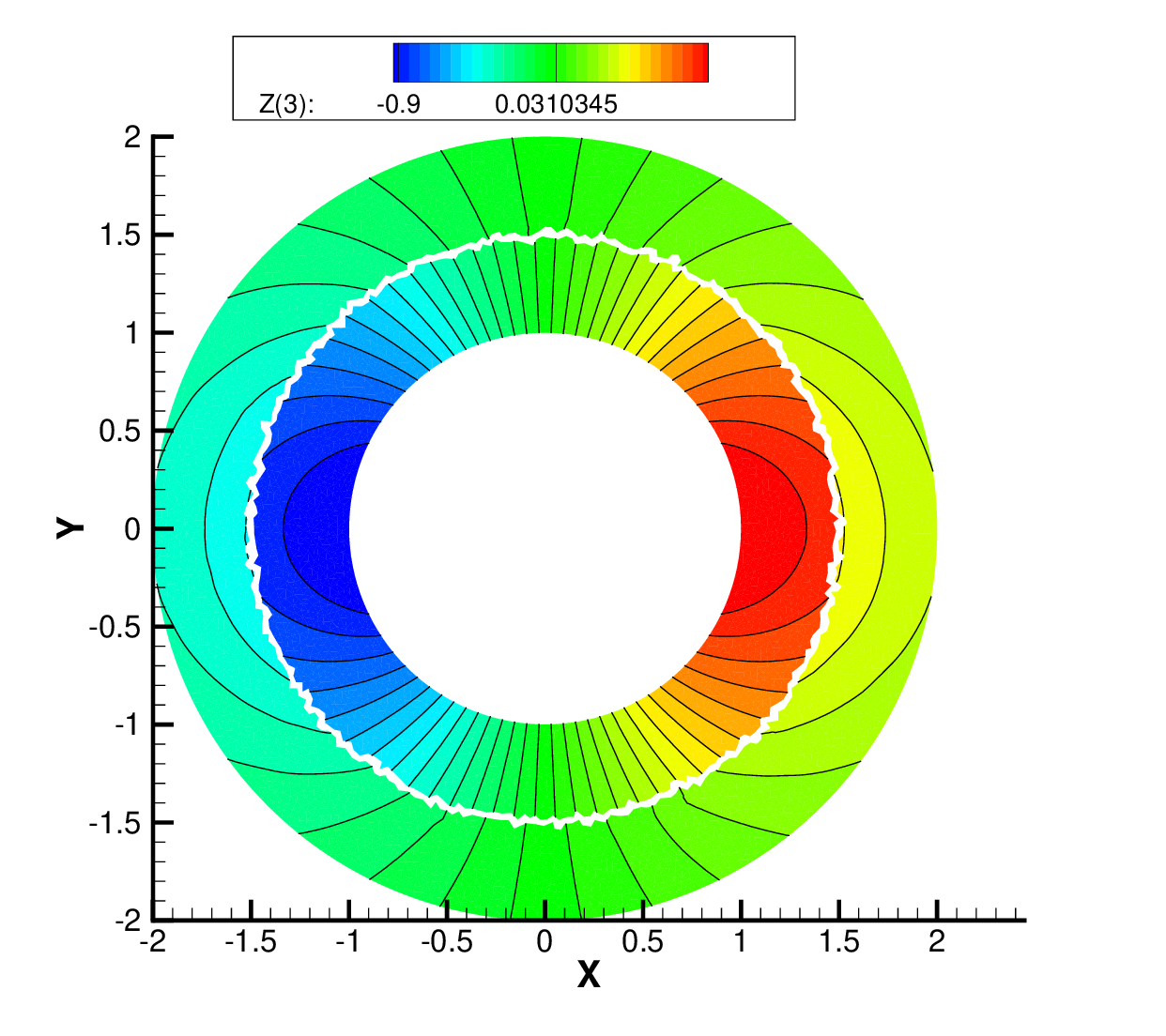}}
		\caption{\textcolor{red}{Planar, transonic source flow: {\color{green} pseudo-temporal evolution} of the eST simulation starting from a SC solution (in terms of $\sqrt{\rho}u$).}}\label{pic40}
	\end{figure*}

\begin{figure*}[h!]
	\centering
		\subfloat[Iteration 1]{%
			\includegraphics[width=0.33\textwidth,trim={0 0 9cm 0},clip]{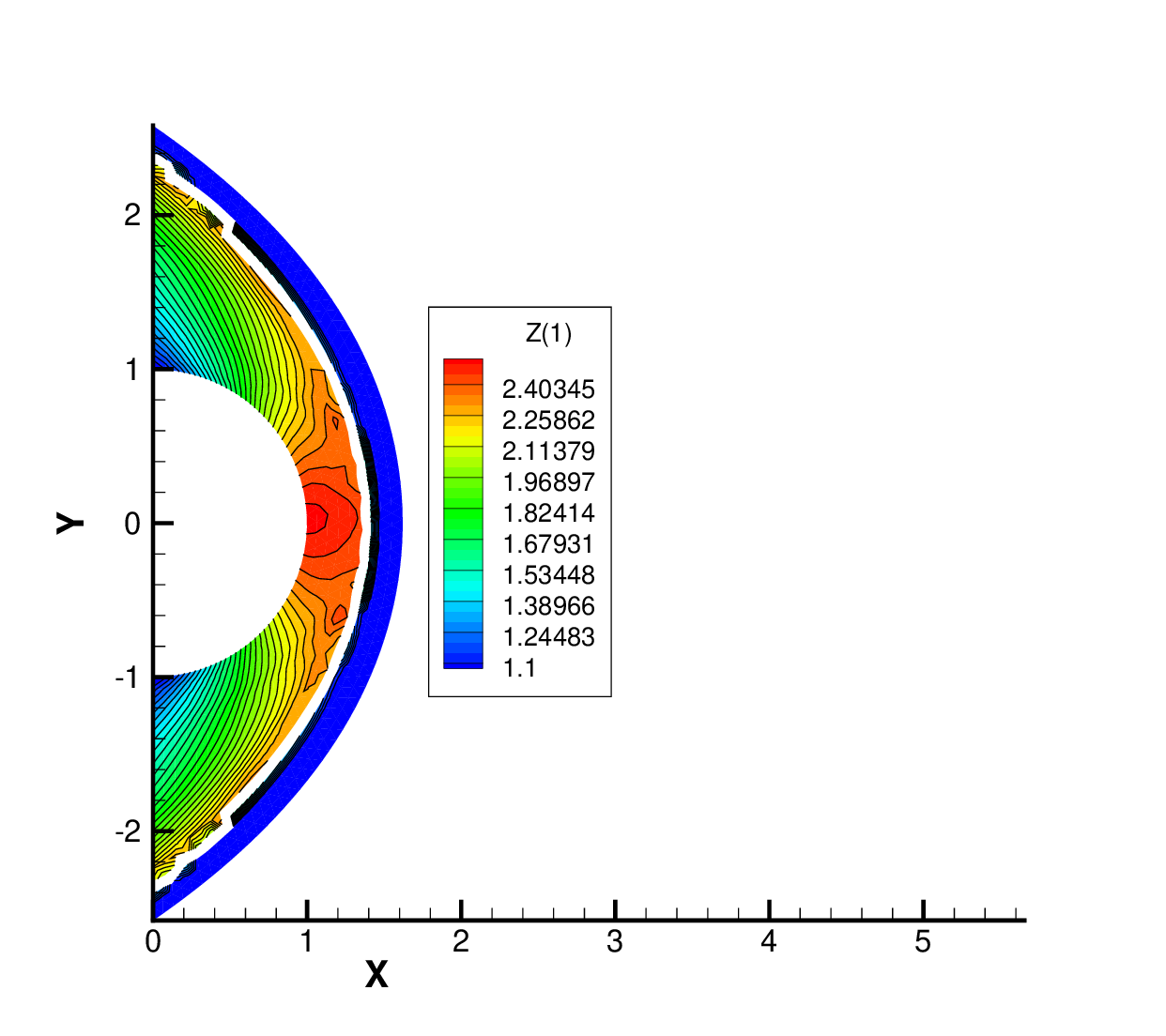}}
		\subfloat[Iteration 101]{%
			\includegraphics[width=0.33\textwidth,trim={0 0 9cm 0},clip]{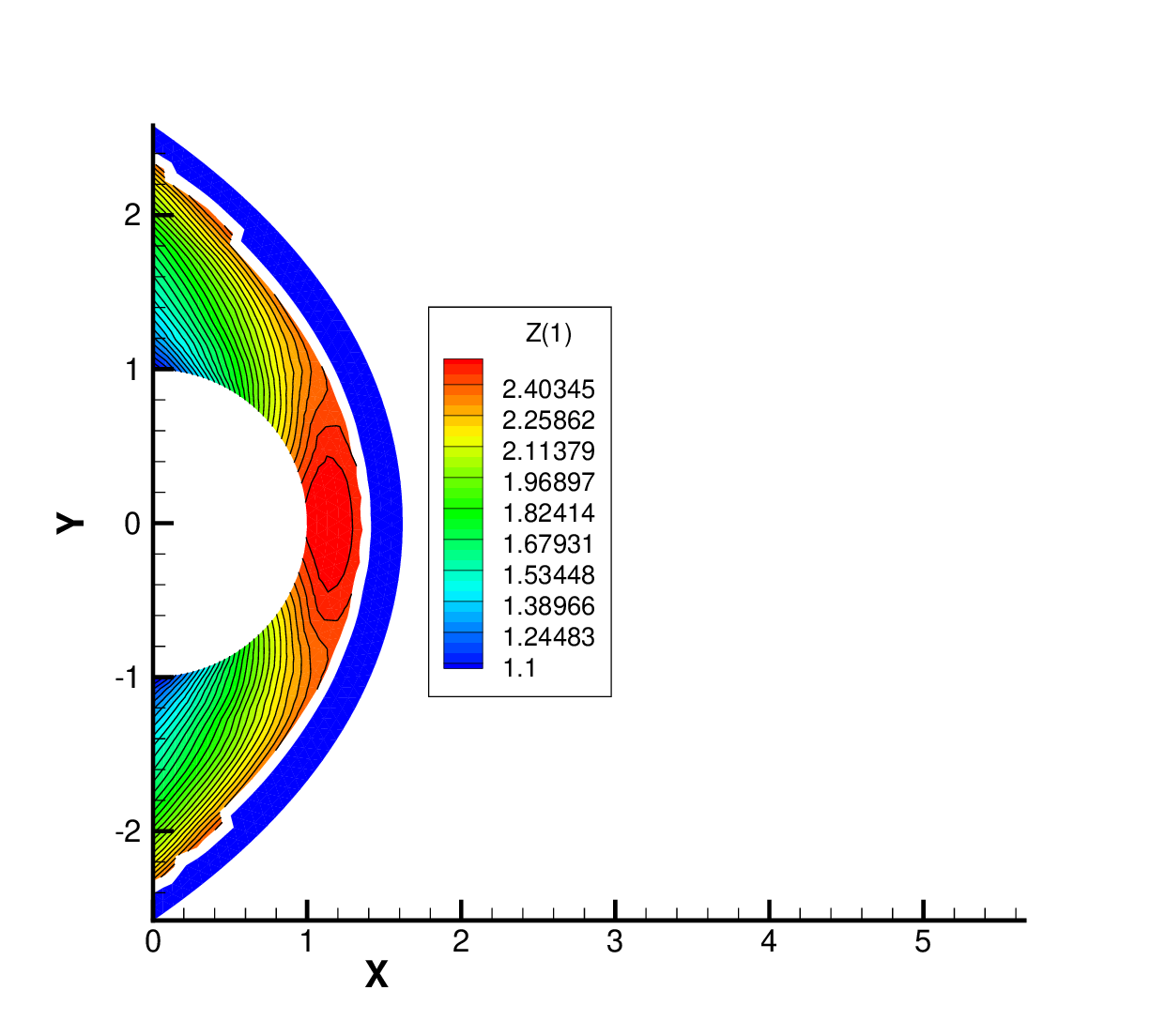}}
		\subfloat[Converged solution]{%
			\includegraphics[width=0.33\textwidth,trim={0 0 9cm 0},clip]{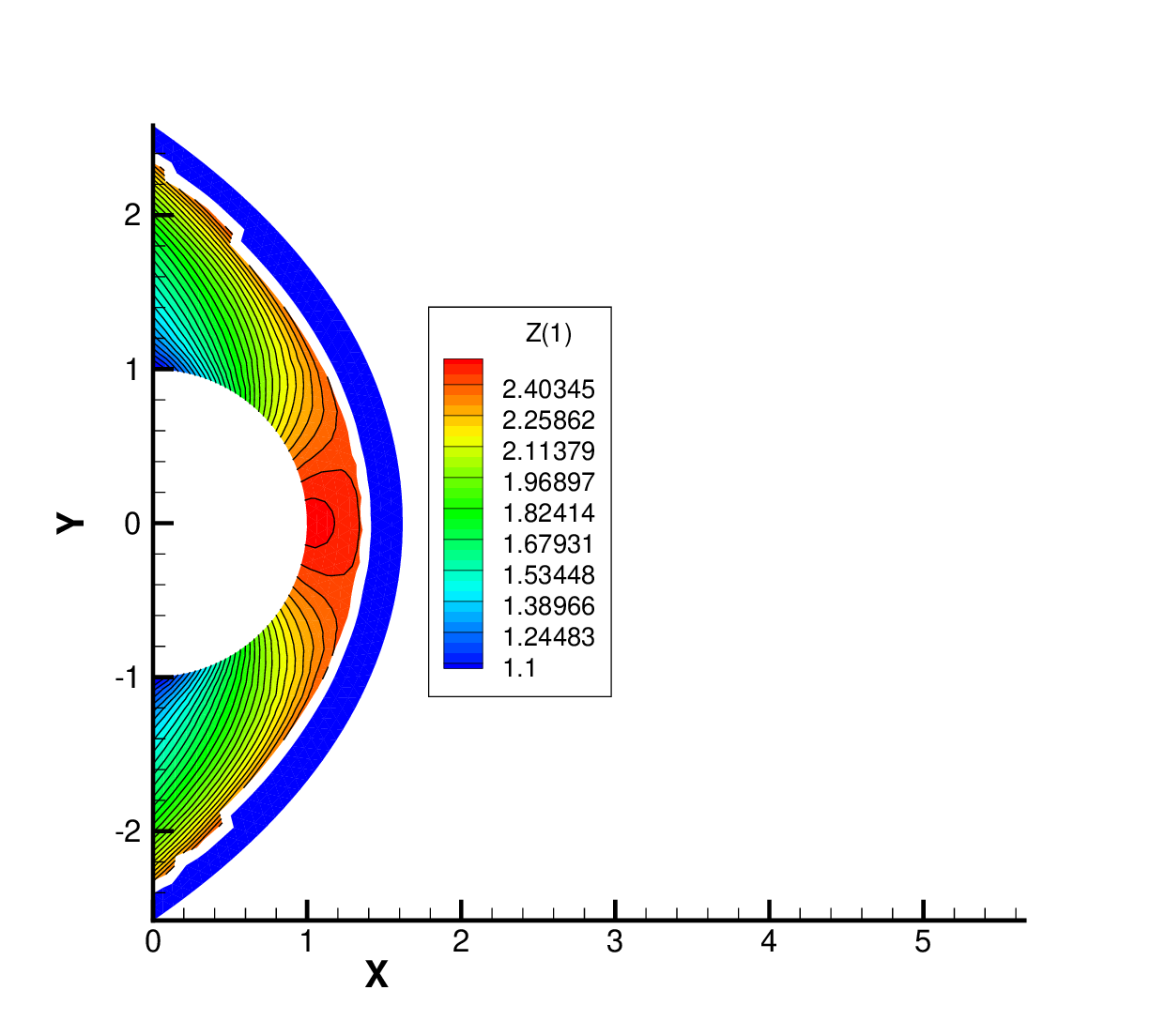}}
		\caption{\textcolor{red}{Hypersonic flow past a circular cylinder: {\color{green} pseudo-temporal evolution} of the eST simulation starting from a SC solution (in terms of $\sqrt{\rho}$).}}\label{pic43}
	\end{figure*}

\begin{figure*}[h!]
	\centering
		\subfloat[Iteration 1]{%
	     		\includegraphics[width=0.33\textwidth]{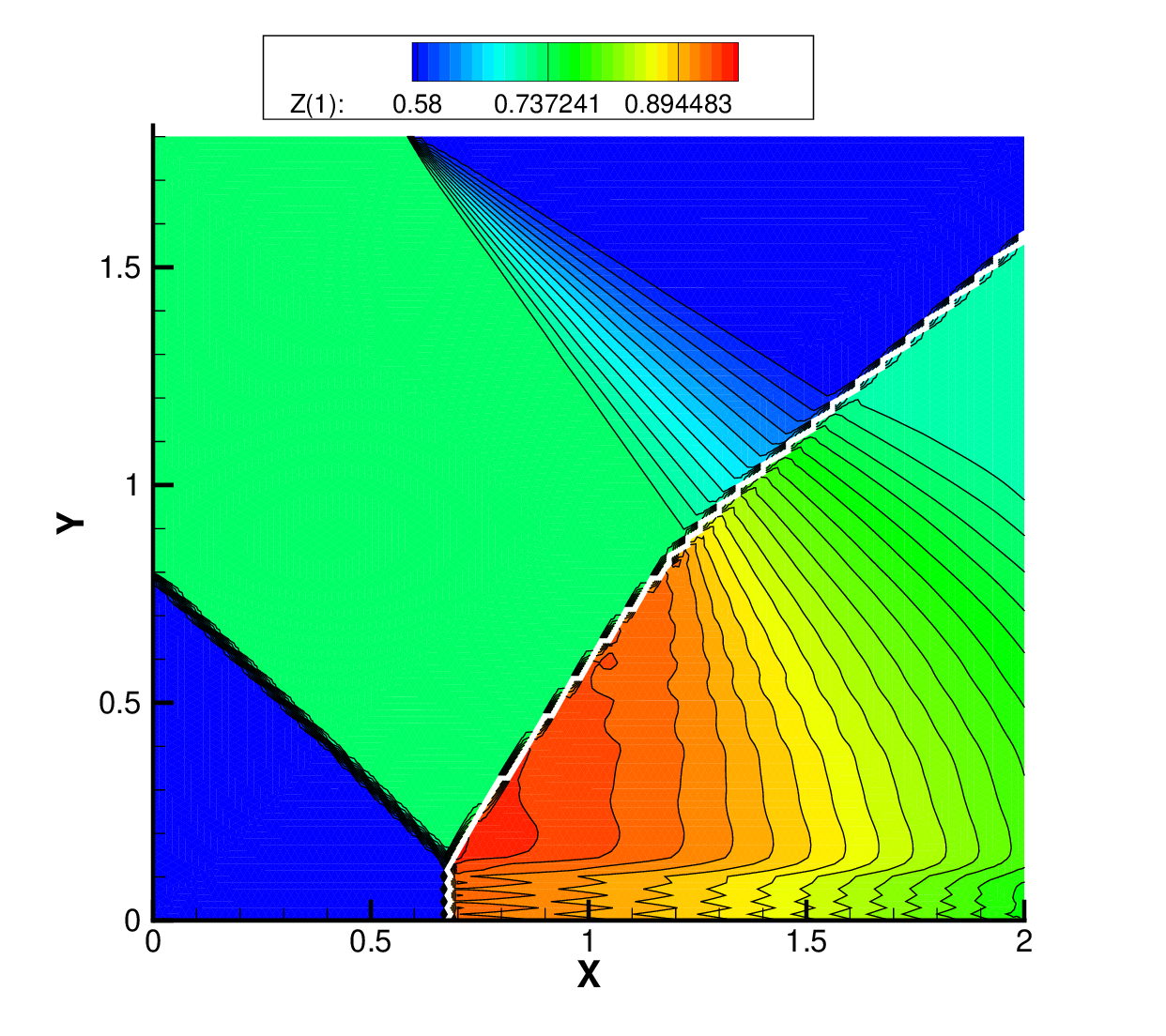}}
		\subfloat[Iteration 201]{%
	     		\includegraphics[width=0.33\textwidth]{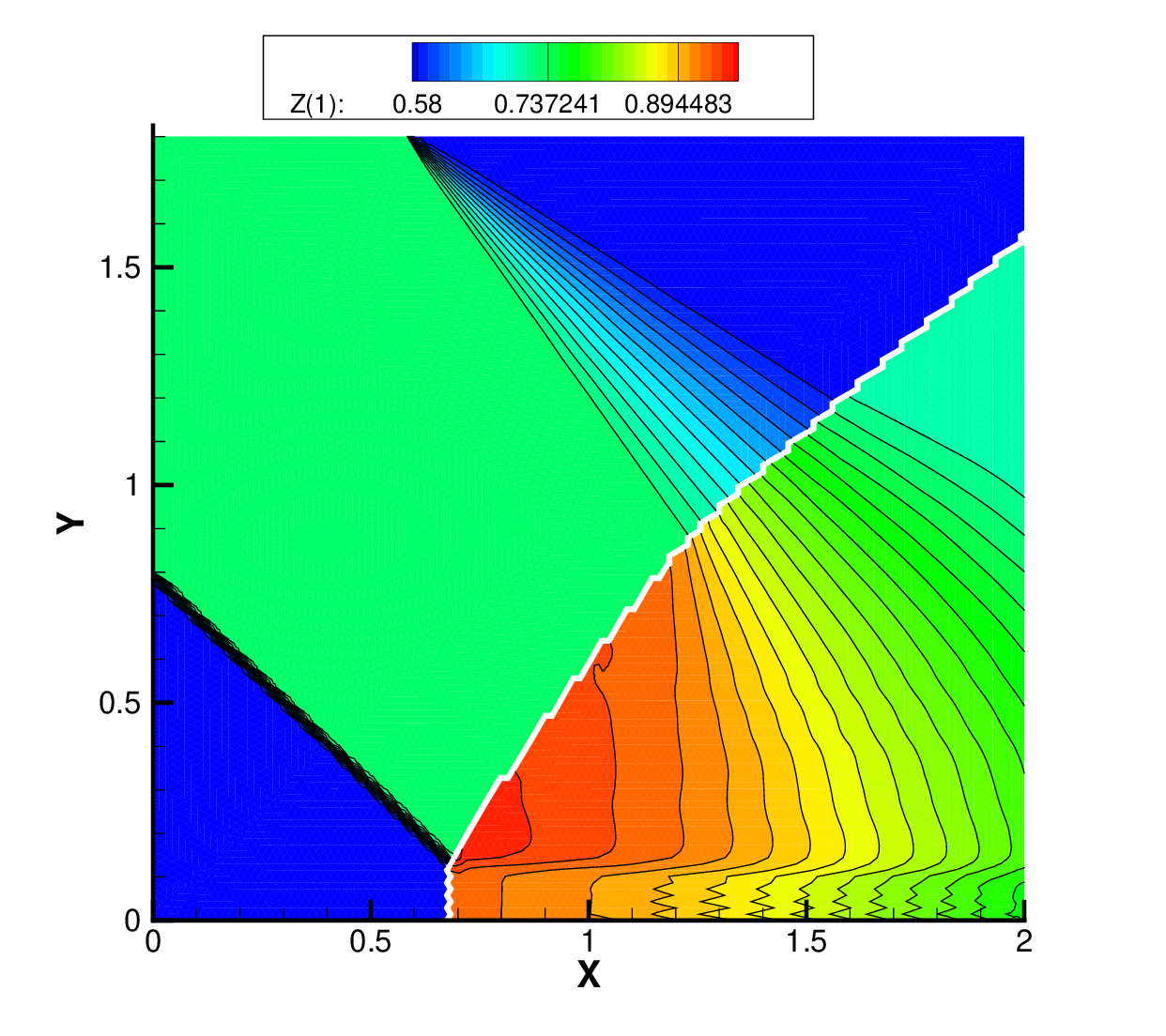}}
		\subfloat[Converged solution]{%
			\includegraphics[width=0.33\textwidth]{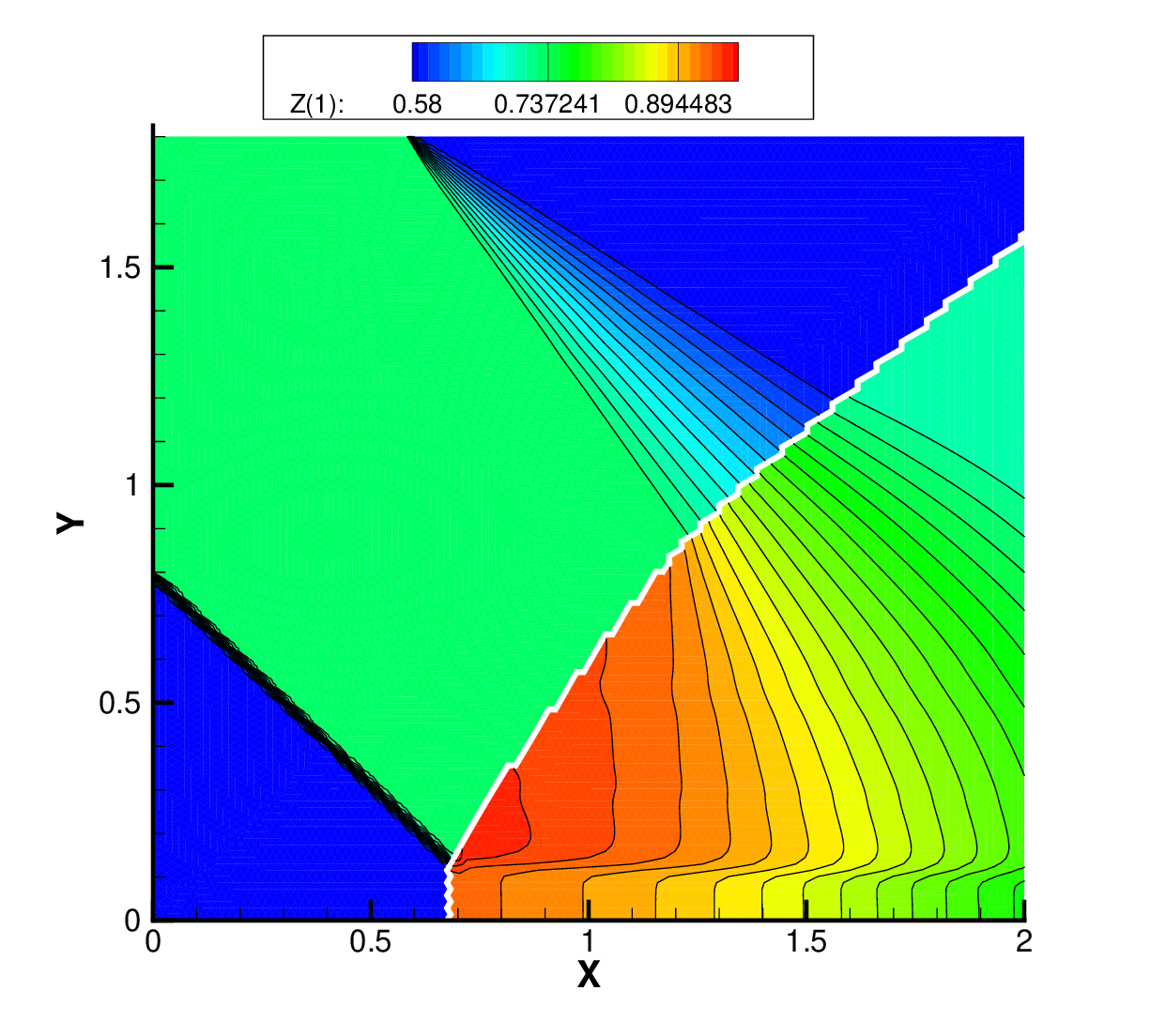}}
		\caption{\textcolor{red}{Steady Mach reflection: {\color{green} pseudo-temporal evolution} of the eST simulation starting from a SC solution (in terms of $\sqrt{\rho}$).}}\label{pic45}
	\end{figure*}

\begin{figure}[h!]
	\begin{center}
	\subfloat[]{\includegraphics[width=0.33\textwidth]{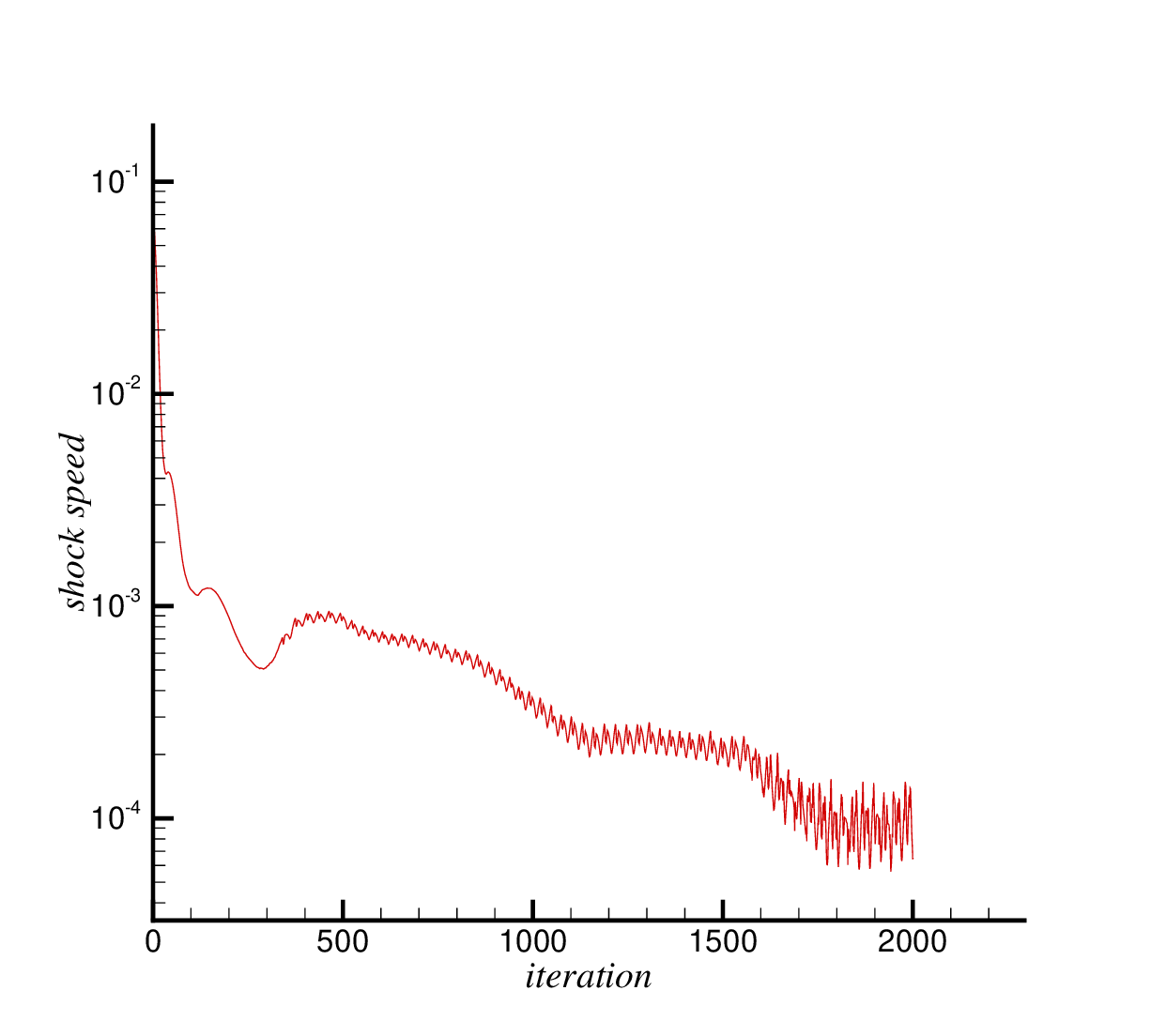}}
	\subfloat[]{\includegraphics[width=0.33\textwidth]{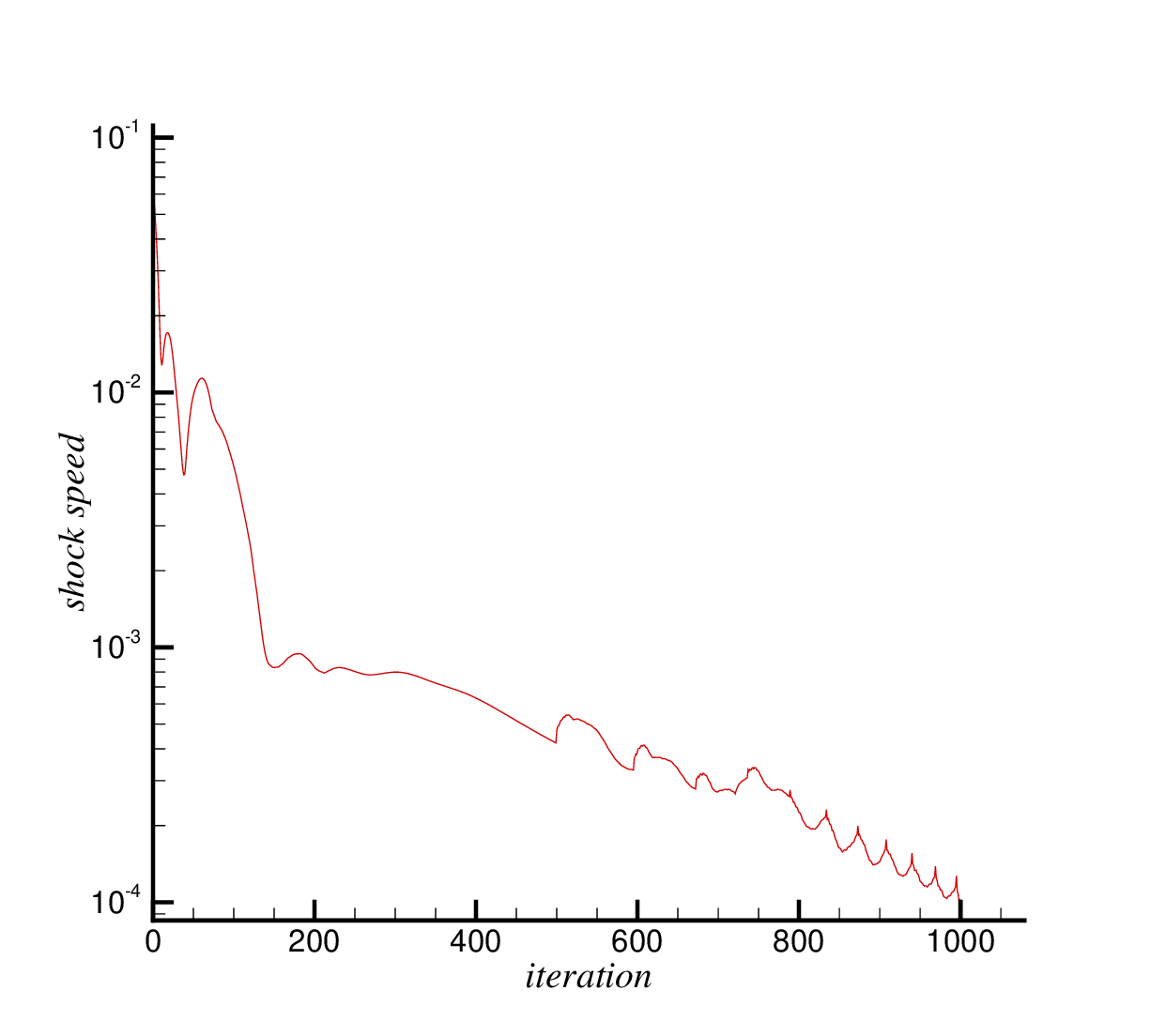}}
	\subfloat[]{\includegraphics[width=0.33\textwidth]{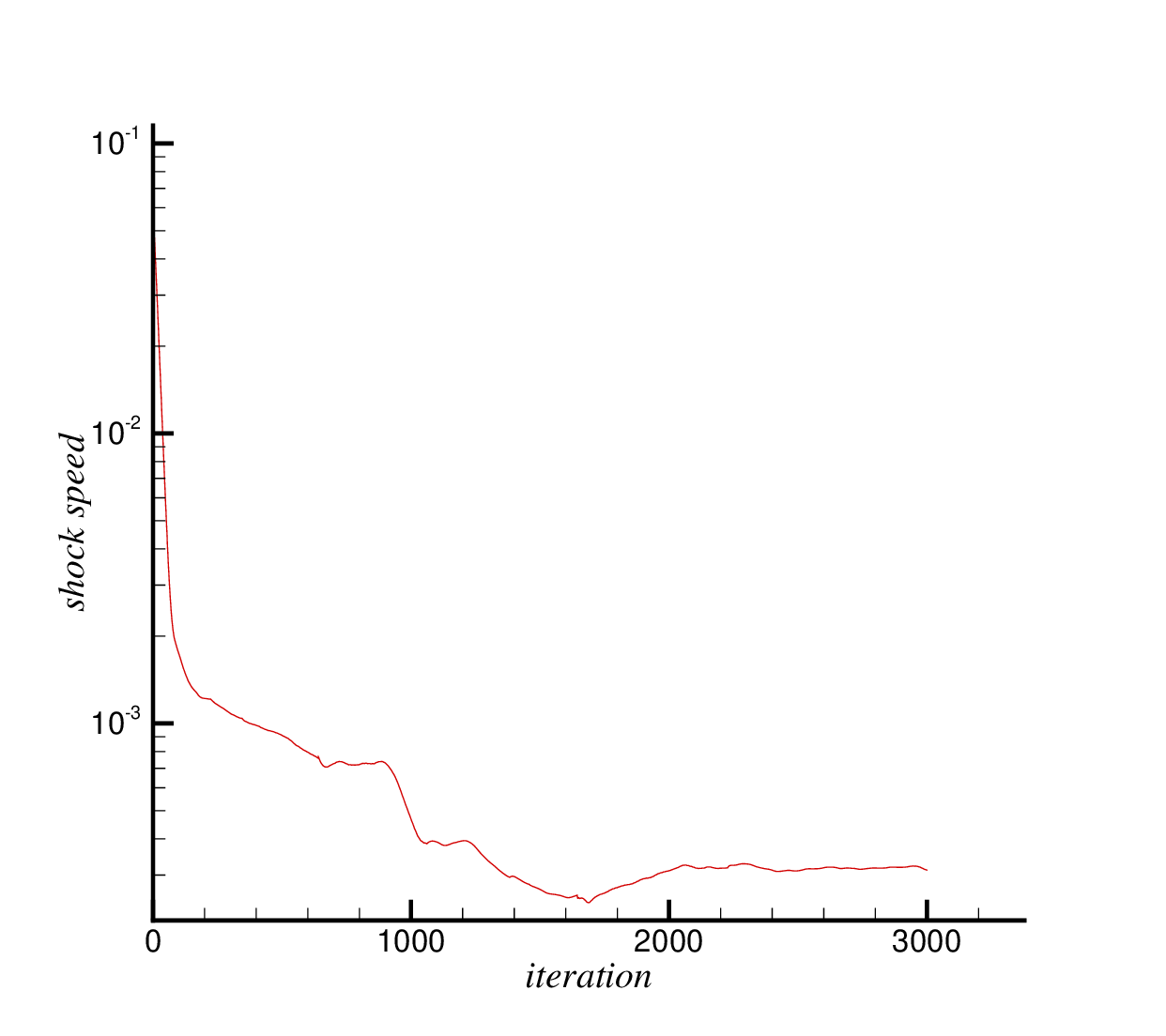}}
	\caption{\textcolor{red}{Shock speed approaching steady state: (a) Planar source flow; (b) Blunt body problem; (c) Steady Mach reflection.}}\label{pic44}
	\end{center}
	\end{figure}

        \textcolor{red}{
	Figures~\ref{pic40}, \ref{pic43} and~\ref{pic45} show a sequence of three frames that
	refer to different instances of
	the pseudo-temporal evolution of the solution for the three test-cases 
	already described in Sects.~\ref{planarflow}, \ref{subsec:blunt} and~\ref{subsec:hybrid}.
	In order to improve readability, the shock-mesh has not been plotted.
	It can be seen that the eST method requires a few hundred pseudo-time steps to get rid of the severe oscillations
        inherited by the SC calculation used to initialize the flow-field.
	Further iterations are required while the shock slows down,
	up to the point when its speed vanishes and it settles to its steady location.
	Convergence of the shock-mesh is monitored by computing, at
	each iteration, a mean shock velocity, averaged over all shock-points.
	Fig.~\ref{pic44}
	shows the pseudo-temporal evolution of the mean shock velocity, plotted against the iteration counter, 
	for all three test-cases. The solution is considered to be converged when this 
	parameter experiences a notable drop, that might also be of several orders of magnitude depending on 
	the shock initial position, as the ones shown in Fig.~\ref{pic44}.   
}

	\clearpage
          \section*{References}
	  \bibliographystyle{elsarticle-num} 
	  \bibliography{biblio}

\begin{thebibliography}{10}
\expandafter\ifx\csname url\endcsname\relax
  \def\url#1{\texttt{#1}}\fi
\expandafter\ifx\csname urlprefix\endcsname\relax\def\urlprefix{URL }\fi
\expandafter\ifx\csname href\endcsname\relax
  \def\href#1#2{#2} \def\path#1{#1}\fi

\bibitem{zaide2011shock}
D.~Zaide, P.~Roe, Shock capturing anomalies and the jump conditions in one
  dimension, in: 20th AIAA Computational Fluid Dynamics Conference, 2011.
\newblock \href {http://dx.doi.org/10.2514/6.2011-3686}
  {\path{doi:10.2514/6.2011-3686}}.

\bibitem{Quirk1994555}
J.~Quirk, A contribution to the great {Riemann }solver debate, International
  Journal for Numerical Methods in Fluids 18~(6) (1994) 555--574.
\newblock \href {http://dx.doi.org/10.1002/fld.1650180603}
  {\path{doi:10.1002/fld.1650180603}}.

\bibitem{Moretti1}
G.~Moretti, Three-dimensional, supersonic, steady flows with any number of
  embedded shocks, in: 12th Aerospace Sciences Meeting, 1974.
\newblock \href {http://dx.doi.org/https://doi.org/10.2514/6.1974-10}
  {\path{doi:https://doi.org/10.2514/6.1974-10}}.

\bibitem{Moretti2}
R.~Marsilio, G.~Moretti, Shock-fitting method for two-dimensional inviscid,
  steady supersonic flows in ducts, Meccanica 24~(4) (1989) 216--222.
\newblock \href {http://dx.doi.org/10.1007/BF01556453}
  {\path{doi:10.1007/BF01556453}}.

\bibitem{Nasuti1}
F.~Nasuti, M.~Onofri, Analysis of unsteady supersonic viscous flows by a
  shock-fitting technique, AIAA journal 34~(7) (1996) 1428--1434.
\newblock \href {http://dx.doi.org/10.2514/6.1995-2159}
  {\path{doi:10.2514/6.1995-2159}}.

\bibitem{Nasuti2}
F.~Nasuti, A multi-block shock-fitting technique to solve steady and unsteady
  compressible flows, in: Computational Fluid Dynamics 2002, Springer, 2003,
  pp. 217--222.

\bibitem{Zhong1}
A.~Prakash, N.~Parsons, X.~Wang, X.~Zhong, High-order shock-fitting methods for
  direct numerical simulation of hypersonic flow with chemical and thermal
  nonequilibrium, Journal of Computational Physics 230~(23) (2011) 8474--8507.
\newblock \href {http://dx.doi.org/10.1016/j.jcp.2011.08.001}
  {\path{doi:10.1016/j.jcp.2011.08.001}}.

\bibitem{Paciorri1}
R.~Paciorri, A.~Bonfiglioli, A shock-fitting technique for {2D} unstructured
  grids, Computers \& Fluids 38~(3) (2009) 715--726.
\newblock \href {http://dx.doi.org/10.1016/j.compfluid.2008.07.007}
  {\path{doi:10.1016/j.compfluid.2008.07.007}}.

\bibitem{Paciorri2}
R.~Paciorri, A.~Bonfiglioli, Shock interaction computations on unstructured,
  two-dimensional grids using a shock-fitting technique, Journal of
  Computational Physics 230~(8) (2011) 3155--3177.
\newblock \href {http://dx.doi.org/10.1016/j.jcp.2011.01.018}
  {\path{doi:10.1016/j.jcp.2011.01.018}}.

\bibitem{Paciorri3}
A.~Bonfiglioli, M.~Grottadaurea, R.~Paciorri, F.~Sabetta, An unstructured,
  three-dimensional, shock-fitting solver for hypersonic flows, Computers \&
  Fluids 73 (2013) 162--174.
\newblock \href {http://dx.doi.org/10.1016/j.compfluid.2012.12.022}
  {\path{doi:10.1016/j.compfluid.2012.12.022}}.

\bibitem{Paciorri4}
A.~Bonfiglioli, R.~Paciorri, L.~Campoli, Unsteady shock-fitting for
  unstructured grids, International Journal for Numerical Methods in Fluids
  81~(4) (2016) 245--261.
\newblock \href {http://dx.doi.org/10.1002/fld.4183}
  {\path{doi:10.1002/fld.4183}}.

\bibitem{Campoli2017}
L.~Campoli, P.~Quemar, A.~Bonfiglioli, M.~Ricchiuto, Shock-fitting and
  predictor-corrector explicit ale residual distribution, in: Shock Fitting:
  Classical Techniques, Recent Developments, and Memoirs of Gino Moretti,
  Springer, 2017, pp. 113--129.

\bibitem{Paciorri5}
M.~S. Ivanov, A.~Bonfiglioli, R.~Paciorri, F.~Sabetta, Computation of weak
  steady shock reflections by means of an unstructured shock-fitting solver,
  Shock Waves 20~(4) (2010) 271--284.
\newblock \href {http://dx.doi.org/10.1007/s00193-010-0266-y}
  {\path{doi:10.1007/s00193-010-0266-y}}.

\bibitem{zou2017shock}
D.~Zou, C.~Xu, H.~Dong, J.~Liu, A shock-fitting technique for cell-centered
  finite volume methods on unstructured dynamic meshes, Journal of
  Computational Physics 345 (2017) 866--882.
\newblock \href {http://dx.doi.org/10.1016/j.jcp.2017.05.047}
  {\path{doi:10.1016/j.jcp.2017.05.047}}.

\bibitem{Chang2019}
S.~Chang, X.~Bai, D.~Zou, Z.~Chen, J.~Liu, An adaptive discontinuity fitting
  technique on unstructured dynamic grids, Shock Waves 29~(8) (2019)
  1103--1115.
\newblock \href {http://dx.doi.org/10.1007/s00193-019-00913-3}
  {\path{doi:10.1007/s00193-019-00913-3}}.

\bibitem{IB0}
C.~S. Peskin, Flow patterns around heart valves: a numerical method, Journal of
  computational physics 10~(2) (1972) 252--271.
\newblock \href {http://dx.doi.org/10.1016/0021-9991(72)90065-4}
  {\path{doi:10.1016/0021-9991(72)90065-4}}.

\bibitem{IB1}
D.~Boffi, L.~Gastaldi, A finite element approach for the immersed boundary
  method, Computers \& Structures 81~(8-11) (2003) 491--501.
\newblock \href {http://dx.doi.org/10.1016/S0045-7949(02)00404-2}
  {\path{doi:10.1016/S0045-7949(02)00404-2}}.

\bibitem{IB2}
R.~Abgrall, H.~Beaugendre, C.~Dobrzynski, An immersed boundary method using
  unstructured anisotropic mesh adaptation combined with level-sets and
  penalization techniques, Journal of Computational Physics 257 (2014) 83--101.
\newblock \href {http://dx.doi.org/10.1016/j.jcp.2013.08.052}
  {\path{doi:10.1016/j.jcp.2013.08.052}}.

\bibitem{IB3}
L.~Nouveau, H.~Beaugendre, C.~Dobrzynski, R.~Abgrall, M.~Ricchiuto, An
  adaptive, residual based, splitting approach for the penalized
  {Navier-Stokes} equations, Computer Methods in Applied Mechanics and
  Engineering 303 (2016) 208--230.
\newblock \href {http://dx.doi.org/10.1016/j.cma.2016.01.009}
  {\path{doi:10.1016/j.cma.2016.01.009}}.

\bibitem{EM0}
A.~Hansbo, P.~Hansbo, An unfitted finite element method, based on {Nitsche}'s
  method, for elliptic interface problems, Computer methods in applied
  mechanics and engineering 191~(47-48) (2002) 5537--5552.
\newblock \href {http://dx.doi.org/10.1016/S0045-7825(02)00524-8}
  {\path{doi:10.1016/S0045-7825(02)00524-8}}.

\bibitem{EM1}
E.~Burman, Ghost penalty, Comptes Rendus Mathematique 348~(21-22) (2010)
  1217--1220.
\newblock \href {http://dx.doi.org/10.1016/j.crma.2010.10.006}
  {\path{doi:10.1016/j.crma.2010.10.006}}.

\bibitem{EM2}
E.~Burman, P.~Hansbo, Fictitious domain methods using cut elements: Iii. a
  stabilized {Nitsche }method for {Stokes}' problem, ESAIM: Mathematical
  Modelling and Numerical Analysis 48~(3) (2014) 859--874.
\newblock \href {http://dx.doi.org/10.1051/m2an/2013123}
  {\path{doi:10.1051/m2an/2013123}}.

\bibitem{EM3}
R.~P. Fedkiw, T.~Aslam, B.~Merriman, S.~Osher, A non-oscillatory eulerian
  approach to interfaces in multimaterial flows (the ghost fluid method),
  Journal of computational physics 152~(2) (1999) 457--492.
\newblock \href {http://dx.doi.org/10.1006/jcph.1999.6236}
  {\path{doi:10.1006/jcph.1999.6236}}.

\bibitem{EM4}
C.~Farhat, A.~Rallu, S.~Shankaran, A higher-order generalized ghost fluid
  method for the poor for the three-dimensional two-phase flow computation of
  underwater implosions, Journal of Computational Physics 227~(16) (2008)
  7674--7700.
\newblock \href {http://dx.doi.org/10.1016/j.jcp.2008.04.032}
  {\path{doi:10.1016/j.jcp.2008.04.032}}.

\bibitem{Scovazzi1}
A.~Main, G.~Scovazzi, The shifted boundary method for embedded domain
  computations. part {I}: {Poisson and Stokes }problems, Journal of
  Computational Physics 372 (2018) 972--995.
\newblock \href {http://dx.doi.org/10.1016/j.jcp.2017.10.026}
  {\path{doi:10.1016/j.jcp.2017.10.026}}.

\bibitem{Scovazzi2}
A.~Main, G.~Scovazzi, The shifted boundary method for embedded domain
  computations. part {II}: {Linear advection--diffusion and incompressible
  Navier--Stokes} equations, Journal of Computational Physics 372 (2018)
  996--1026.
\newblock \href {http://dx.doi.org/10.1016/j.jcp.2018.01.023}
  {\path{doi:10.1016/j.jcp.2018.01.023}}.

\bibitem{Scovazzi3}
T.~Song, A.~Main, G.~Scovazzi, M.~Ricchiuto, The shifted boundary method for
  hyperbolic systems: Embedded domain computations of linear waves and shallow
  water flows, Journal of Computational Physics 369 (2018) 45--79.
\newblock \href {http://dx.doi.org/10.1016/j.jcp.2018.04.052}
  {\path{doi:10.1016/j.jcp.2018.04.052}}.

\bibitem{Glimm2016}
D.~She, R.~Kaufman, H.~Lim, J.~Melvin, A.~Hsu, J.~Glimm, Front tracking
  methods, in: Handbook of Numerical Methods for Hyperbolic Problems, Vol.~17,
  Elsevier, 2016, pp. 383---402.
\newblock \href {http://dx.doi.org/https://doi.org/10.1016/bs.hna.2016.07.004}
  {\path{doi:https://doi.org/10.1016/bs.hna.2016.07.004}}.

\bibitem{lishifted}
K.~Li, N.~M. Atallah, G.~A. Main, G.~Scovazzi, The shifted interface method: A
  flexible approach to embedded interface computations, International Journal
  for Numerical Methods in Engineering 121~(3) (2020) 492--518.
\newblock \href {http://dx.doi.org/10.1002/nme.6231}
  {\path{doi:10.1002/nme.6231}}.

\bibitem{pepe2015unstructured}
R.~Pepe, A.~Bonfiglioli, A.~D’Angola, G.~Colonna, R.~Paciorri, An
  unstructured shock-fitting solver for hypersonic plasma flows in chemical
  non-equilibrium, Computer Physics Communications 196 (2015) 179--193.
\newblock \href {http://dx.doi.org/10.1016/j.cpc.2015.06.005}
  {\path{doi:10.1016/j.cpc.2015.06.005}}.

\bibitem{SFB:Andrea}
A.~Lani, V.~De~Amicis, {SF}: An open source object-oriented platform for
  unstructured shock-fitting methods, in: M.~Onofri, R.~Paciorri (Eds.), Shock
  Fitting: Gino Moretti's Memoires, Fundamentals and Recent Evolution of his
  Shock-Fitting Technique, Springer International Publishing, 2017, pp.
  127--153.

\bibitem{PACIORRI2020109196}
R.~Paciorri, A.~Bonfiglioli, Accurate detection of shock waves and shock
  interactions in two-dimensional shock-capturing solutions, Journal of
  Computational Physics 406 (2020) 109196.
\newblock \href {http://dx.doi.org/https://doi.org/10.1016/j.jcp.2019.109196}
  {\path{doi:https://doi.org/10.1016/j.jcp.2019.109196}}.

\bibitem{Scovazzi4}
L.~Nouveau, M.~Ricchiuto, G.~Scovazzi, High-order gradients with the shifted
  boundary method: An embedded enriched mixed formulation for elliptic pdes,
  Journal of Computational Physics 398 (2019) 108898.
\newblock \href {http://dx.doi.org/10.1016/j.jcp.2019.108898}
  {\path{doi:10.1016/j.jcp.2019.108898}}.

\bibitem{persson}
M.~Zahr, A.~Shi, P.-O. Persson, Implicit shock tracking using an
  optimization-based high-order discontinuous galerkin method, Journal of
  Computational Physics 410 (2020) 109385.
\newblock \href {http://dx.doi.org/{10.1016/j.jcp.2020.109385}}
  {\path{doi:{10.1016/j.jcp.2020.109385}}}.

\bibitem{corrigan}
A.~Corrigan, A.~D. Kercher, D.~A. Kessler, A moving discontinuous galerkin
  finite element method for flows with interfaces, International Journal for
  Numerical Methods in Fluids 89~(9) (2019) 362--406.
\newblock \href {http://dx.doi.org/10.1002/fld.4697}
  {\path{doi:10.1002/fld.4697}}.

\bibitem{Bonfiglioli1}
A.~Bonfiglioli, {Fluctuation splitting schemes for the compressible and
  incompressible Euler and Navier-Stokes equations}, International Journal of
  Computational Fluid Dynamics 14~(1) (2000) 21--39.
\newblock \href {http://dx.doi.org/10.1080/10618560008940713}
  {\path{doi:10.1080/10618560008940713}}.

\bibitem{dr2017}
H.~Deconinck, M.~Ricchiuto, Residual Distribution Schemes: Foundations and
  Analysis, John Wiley \& Sons, Ltd, 2017, pp. 1--53.
\newblock \href {http://dx.doi.org/10.1002/0470091355.ecm054}
  {\path{doi:10.1002/0470091355.ecm054}}.

\bibitem{ar2017}
R.~Abgrall, M.~Ricchiuto, High-Order Methods for CFD, John Wiley \& Sons, Ltd,
  2017, pp. 1--54.
\newblock \href {http://dx.doi.org/10.1002/9781119176817.ecm2112}
  {\path{doi:10.1002/9781119176817.ecm2112}}.

\bibitem{Bonfiglioli2}
A.~Bonfiglioli, R.~Paciorri, A mass-matrix formulation of unsteady fluctuation
  splitting schemes consistent with {Roe}'s parameter vector, International
  Journal of Computational Fluid Dynamics 27~(4-5) (2013) 210--227.
\newblock \href {http://dx.doi.org/10.1080/10618562.2013.813491}
  {\path{doi:10.1080/10618562.2013.813491}}.

\bibitem{Roe1}
P.~L. Roe, Approximate riemann solvers, parameter vectors, and difference
  schemes, Journal of computational physics 43~(2) (1981) 357--372.
\newblock \href {http://dx.doi.org/10.1006/jcph.1997.5705}
  {\path{doi:10.1006/jcph.1997.5705}}.

\bibitem{bonhaus1998higher}
D.~Bonhaus, A higher order accurate finite element method for viscous
  compressible flows, Ph.D. thesis, Virginia Polytechnic Institute and State
  University (1998).

\bibitem{Paciorri6}
A.~Bonfiglioli, R.~Paciorri, Convergence analysis of shock-capturing and
  shock-fitting solutions on unstructured grids, AIAA journal 52~(7) (2014)
  1404--1416.
\newblock \href {http://dx.doi.org/10.2514/1.J052567}
  {\path{doi:10.2514/1.J052567}}.

\bibitem{Campobasso2011}
M.~S. Campobasso, M.~H. Baba-Ahmadi, Ad-hoc boundary conditions for cfd
  analyses of turbomachinery problems with strong flow gradients at farfield
  boundaries, Journal of Turbomachinery 133~(4).
\newblock \href {http://dx.doi.org/10.1115/1.4002985}
  {\path{doi:10.1115/1.4002985}}.

\bibitem{Triangle1}
J.~R. Shewchuk, Triangle: Engineering a 2d quality mesh generator and delaunay
  triangulator, in: Workshop on Applied Computational Geometry, 1996, pp.
  203--222.
\newblock \href {http://dx.doi.org/10.1007/BFb0014497}
  {\path{doi:10.1007/BFb0014497}}.

\bibitem{Triangle2}
J.~Shewchuk, Triangle mesh generator, Available at
  \url{https://www.cs.cmu.edu/~quake/triangle.html}.
\newblock \href {http://dx.doi.org/https://doi.org/10.1007/BFb0014497}
  {\path{doi:https://doi.org/10.1007/BFb0014497}}.

\bibitem{Grottadaurea2011}
M.~Grottadaurea, R.~Paciorri, A.~Bonfiglioli, F.~Sabetta, M.~Onofri,
  D.~Bianchi, Numerical simulation of hypersonic flows past three-dimensional
  blunt bodies through an unstructured shock-fitting solver., in: International
  Space Planes and Hypersonic Systems and Technologies Conferences, American
  Institute of Aeronautics and Astronautics, 2011.
\newblock \href {http://dx.doi.org/https://doi.org/10.2514/6.2011-2288}
  {\path{doi:https://doi.org/10.2514/6.2011-2288}}.

\bibitem{ALAUZET201828}
F.~Alauzet, A.~Loseille, G.~Olivier, Time-accurate multi-scale anisotropic mesh
  adaptation for unsteady flows in cfd, Journal of Computational Physics 373
  (2018) 28 -- 63.
\newblock \href {http://dx.doi.org/https://doi.org/10.1016/j.jcp.2018.06.043}
  {\path{doi:https://doi.org/10.1016/j.jcp.2018.06.043}}.

\bibitem{Muller1}
J.-D. M{\"u}ller, P.~L. Roe, H.~Deconinck, A frontal approach for internal node
  generation in delaunay triangulations, International Journal for Numerical
  Methods in Fluids 17~(3) (1993) 241--255.
\newblock \href {http://dx.doi.org/10.1002/fld.1650170305}
  {\path{doi:10.1002/fld.1650170305}}.

\bibitem{Muller2}
J.~M{\"u}ller, Delaundo mesh generator, Available at
  \url{http://www.ae.metu.edu.tr/tuncer/ae546/prj/delaundo/}.

\bibitem{Rusanov}
A.~Lyubimov, V.~Rusanov, {Gas Flows Past Blunt Bodies, Part II: Tables of the
  Gasdynamic Functions}, NASA TT F-715.

\bibitem{SSIVpaper}
N.~Duquesne, T.~Alziary~de Roquefort, Numerical investigation of a
  three-dimensional turbulent shock/shock interaction, in: 36th AIAA Aerospace
  Sciences Meeting and Exhibit, 1998.
\newblock \href {http://dx.doi.org/https://doi.org/10.2514/6.1998-774}
  {\path{doi:https://doi.org/10.2514/6.1998-774}}.

\end{thebibliography}
	%
	%
	%
	%
	%
	\end{document}